\documentclass[12pt]{imsart}

\RequirePackage[OT1]{fontenc}
\RequirePackage{amsthm,amsmath}
\RequirePackage[numbers]{natbib}
\RequirePackage[colorlinks,citecolor=blue,urlcolor=blue]{hyperref}

\usepackage{amsthm,amsmath,amssymb}
\usepackage[numbers]{natbib}
\usepackage{color}

\usepackage{hyperref}
\hypersetup{
    colorlinks,%
    citecolor=blue,%
    filecolor=red,%
    linkcolor=red,%
    urlcolor=black,
}
\usepackage{hypernat}
\usepackage[top=2.5cm, bottom=2.5cm, left=2.8cm,
right=2.8cm]{geometry}

\newtheorem{theorem}{Theorem}[section]

\newtheorem{lemma}[theorem]{Lemma}

\newtheorem{corollary}[theorem]{Corollary}
\newtheorem{remark}{Remark}[section]

\newcommand{\cqfd}{\hfill{$\Box$}}
\newcommand{\eps}{{\epsilon}}
\newcommand{\E}{{\mathbb E}}

\newcommand{\R}{{\mathbb R}}

\renewcommand{\P}{{\mathbb P}}

\newcommand{\F}{{\mathcal F}}

\newcommand{\density} {(R5)}
\newcommand{\sigmamin}{(R3)}
\newcommand{\momentq}{(R4)}
\newcommand{\monorandom}{(R1)}
\newcommand{\momentexpo}{(R4)}
\newcommand{\densityUI}{(R2)}

\newdimen\AAdi%
\newbox\AAbo%
\def\AArm{\fam0 }%\tenrm}%
\def\AAk#1#2{\setbox\AAbo=\hbox{#2}\AAdi=\wd\AAbo\kern#1\AAdi{}}%
\def\BBone{{\AArm 1\AAk{-.8}{I}I}}%

\newcounter{rcnt}[section]

\renewcommand{\d}{\mbox{ \sl\em d}}
\def\argmin{\mathop{\rm argmin}}
\def\argmax{\mathop{\rm argmax}}

\arxiv{arXiv:0000.0000}

\startlocaldefs
\numberwithin{equation}{section}
\theoremstyle{plain}

\endlocaldefs

\begin{document}

\begin{frontmatter}
\title{Divide and Conquer in Non-standard Problems and the Super-efficiency Phenomenon}
\runtitle{Divide and Conquer in Non-standard Problems}
%\thankstext{T1}{Footnote to the title with the ``thankstext'' command.}

\begin{aug}
\author{\fnms{Moulinath} \snm{Banerjee}\thanksref{t1}\ead[label=e1]{moulib@umich.edu}},
\author{\fnms{Cecile} \snm{Durot}%\thanksref{t3,m1,m2}
\ead[label=e2]{cecile.durot@gmail.com}}
\and
\author{\fnms{Bodhisattva} \snm{Sen}\thanksref{t3}
\ead[label=e3]{bodhi@stat.columbia.edu}
%\ead[label=u1,url]{http://www.foo.com}
}

\thankstext{t1}{Supported by NSF Grant DMS-1308890}
%\thankstext{t2}{First supporter of the project}
\thankstext{t3}{Supported by NSF CAREER Grant DMS-1150435}
\runauthor{Banerjee, M., Durot, C. and Sen, B.}

\affiliation{University on Michigan, Universit\'{e} Paris Ouest Nanterre La D\'{e}fense and Columbia University}

\address{University of Michigan \\
451 West Hall, 1085 South University \\
Ann Arbor, MI 48109 \\
\printead{e1}
%\phantom{E-mail:\ }\printead*{e2}
}

\address{Universit\'{e} Paris Ouest Nanterre La D\'{e}fense \\
200 avenue de la r\'{e}publique\\
92001 Nanterre Cedex, France\\
\printead{e2}}

\address{Columbia University \\
1255 Amsterdam Av., Room \# 1032 SSW\\
New York, NY 10027\\
\printead{e3}}
\end{aug}

\begin{abstract}
We study how the divide and conquer principle --- partition the available data into subsamples, compute an estimate from each subsample and combine these appropriately to form the final estimator --- works in non-standard problems where rates of convergence are typically slower than $\sqrt{n}$ and limit distributions are non-Gaussian, with a special emphasis on the least squares estimator (and its inverse) of a monotone regression function. We find that the pooled estimator, obtained by averaging non-standard estimates across the mutually exclusive subsamples, outperforms the non-standard estimator based on the entire sample in the sense of \emph{pointwise inference}. We also show that, under appropriate conditions, if the number of subsamples is allowed to increase at appropriate rates, the pooled estimator is asymptotically normally distributed with a variance that is empirically estimable from the subsample-level estimates. Further, in the context of monotone function estimation we show that this gain in pointwise efficiency comes at a price ---  the pooled estimator's performance, in a \emph{uniform sense} (maximal risk) over a class of models worsens as the number of subsamples increases, leading to a version of the super-efficiency phenomenon. In the process, we develop analytical results for the order of the bias in isotonic regression, which are of independent interest. 

\end{abstract}



\end{frontmatter}

\section{Introduction}

Suppose that $W_1,\ldots, W_N$ are i.i.d.~random elements having a common distribution $P$. We assume that $P$ is unknown and $\theta_{0} \equiv \theta_{0}(P)$ is a finite dimensional parameter of interest. In this paper we focus on non-standard statistical problems where a natural estimator $\hat \theta$ (of $\theta_{0}$) converges in distribution to a non-normal limit at a rate slower than $N^{1/2}$, i.e., 
\begin{equation}\label{eq:LimDist}
	r_N(\hat \theta -  \theta_0) \stackrel{d}{\to} G,
\end{equation}
where $r_N = o(\sqrt{N})$ and $G$ is non-normal, has mean zero and finite variance $\sigma^2$. However, $\sigma^2$, the variance of $G$, can depend on $P$ in a complicated fashion which often makes it difficult to use~\eqref{eq:LimDist} to construct confidence intervals (CIs) and hypothesis tests for $\theta_{0}$.  Such non-standard limits primarily arise due to the inherent lack of smoothness in the underlying estimation procedure. Also, in many such scenarios the computation of $\hat \theta$ is complicated, requiring computationally intensive algorithms. Thus, in the face of a humongous sample size $N$ --- quite common with present-day data --- these problems present a significant challenge both in computation and inference. 

In this paper, our primary goal is to investigate how such non-standard estimates behave under a {\it sample-splitting} strategy, the so-called ``divide-and-conquer'' method that has been much used in the analysis of massive data sets; see e.g.,~\cite{LiEtAl13},~\cite{Zhang13} and~\cite{Zhao14}. In divide and conquer, the available data is partitioned into subsamples, an estimate of $\theta_{0}$ is computed from each subsample, and finally the subsample level estimates are combined appropriately to form the final estimator. Our combining/pooling strategy will be based on averaging the estimators obtained from the different subsamples.

A rich class of such problems arises in the world of ``cube-root asymptotics'' (see~\cite{KP90}), which include, e.g., estimation of the mode (see~\cite{C64}), Manski's maximum score estimator (see~\cite{M75}), change-point estimation under smooth mis-specification (see~\cite{BM07}), least absolute median of squares (see~\cite{R84}), shorth estimation (see e.g.,~\cite{DT06}), and last but not least, isotonic regression (see e.g.,~\cite{B55},\cite{G56}). We elaborate below on the last of the aforementioned examples: the estimation of a monotone function. 

Consider  i.i.d.~data $\{W_i:=(X_i,Y_i):i =1,\ldots, N\}$ from the regression model 
\begin{equation}\label{eq:IsoReg}
	Y = \mu(X) + \eps
\end{equation}
where $Y \in \R$ is the response variable, $X \in [0,1]$ (with density $f$) is the covariate, $\mu$ is the unknown {\it nonincreasing} regression function, and $\eps$ is independent of $X$ and has mean 0 and variance $v^2>0$. The goal is to estimate $\mu:[0,1] \to \R$ nonparametrically, under the known constraint of monotonicity. We will consider the least squares estimator (LSE) $\hat \mu$ defined as 
\begin{equation}\label{eq:LSE}
\hat \mu \in \arg \min_{\psi \downarrow} \sum_{i=1}^n  (Y_i- \psi(X_i))^2,
\end{equation} 
 where the minimization is over all nonincreasing functions $\psi:[0,1] \to \R$. We know that $\hat \mu$ is unique at the data points $X_i$'s and can be connected  to  the left-hand slope of the least concave majorant of the cumulative sum diagram (see e.g., \cite[Chapter 1]{RWD88}). If $\mu'(t_0) \ne 0$, where $t_0$ is an interior point in the support of $X$, then
\begin{equation}
\label{eq:Chernoff-forward} 
	N^{1/3}(\hat \mu(t_0) - \mu(t_0)) \stackrel{d}{\to} \kappa \mathbb{Z},
\end{equation}
with $\kappa := |4 v^2 \mu'(t_0)/f(t_0)|^{1/3}$ and $\mathbb{Z} := \argmin_{s \in \R} \{W(s) + s^2\}$ (where $W$ is a standard two-sided Brownian motion starting at 0)  has the so-called Chernoff's distribution; see e.g., Theorem 1 in~\cite{Wright81}. It is known that $\mathbb{Z}$ is symmetric (around 0) and has mean zero. Lastly $\sigma^2 = \mathrm{Var}(\kappa \mathbb{Z})$, the variance of the limiting distribution, is difficult to estimate as it involves the derivative of $\mu$, the estimation of which is well-known to be a challenging problem (see e.g.,~\cite{BW05}). 

A closely related problem is the estimation of the inverse isotonic function at a point. If $a$ is an interior point in the range of $\mu$ and $t_0 = \mu^{-1}(a) \in (0,1)$ satisfies $\mu'(t_0) \ne 0$, then
\begin{equation}
\label{eq:Chernoff-backward} 
	N^{1/3}(\hat \mu^{-1}(a) - \mu^{-1}(a)) \stackrel{d}{\to} \tilde{\kappa} \mathbb{Z},
\end{equation}
where $\tilde{\kappa} := |4 v^2/\mu'(t_0)^2 f(t_0)|^{1/3}$; this can be derived, e.g., from the arguments in~\cite{Durot08}. Similar results hold across a vast array of monotone function problems: in particular, in the heteroscedastic regression model where $\eps$ is no longer independent of $X$, in Grenander's problem (\cite{G56}) on the estimation of a monotone density, and monotone response models as considered in~\cite{B07}.

We now formally introduce the sample-splitting idea. Assume that $N$ is large and write $N = n \times m$, where $n$ is still large and $m$ relatively smaller (e.g., $n = 1000$, $m = 50$, so that $N = 50000$). We define our new ``averaged'' estimator as follows:
\begin{enumerate}
	\item Divide the set of samples $W_1,\ldots, W_N$ into $m$ disjoint subsets $S_1,\ldots, S_m$.
	
	\item For each $j = 1,\ldots, m$, compute the estimator $\hat \theta_j$ based on the data points in $S_j$.
	
	\item Average together these estimators to obtain the final `pooled' estimator: 
\begin{equation}\label{eq:GlEst}
	\bar \theta = \frac{1}{m} \sum_{j=1}^m \hat \theta_j.
\end{equation}
\end{enumerate}
Observe that if the computation of $\hat \theta$, the global estimator based on all $N$ observations, is of super-linear computational complexity in the sample size, computing $\bar \theta$ saves resources compared to $\hat \theta$. Further, the computation of $\bar \theta$ can be readily  parallelized, using $m$ CPU's. This idea of averaging estimators based on disjoint subsets of the data has been used by many authors recently to estimate nonparametric functions, but typically  under smoothness constraints; see e.g.,~\cite{Zhang13},~\cite{Zhao14}, and also~\cite{LiEtAl13} for a discussion with a broader scope.  The above papers illustrate that the sample-splitting approach significantly reduces the required amount of primary memory and computation time in a variety of cases, yet statistical optimality  --- in the sense that the resulting estimator is as efficient (or minimax rate optimal) as the global estimate based on applying the estimation algorithm to the entire data set  --- is retained. 

We show in this paper that in certain non-standard problems, by sample-splitting not only do we have computational gains, but the resulting estimator $\bar \theta$  acquires a \emph{faster rate} of convergence to a normal limit. This is quite interesting, and to the best of our knowledge, hitherto unobserved in the statistical literature. However, this faster rate of convergence of the pooled estimator at a point $\theta_0$ is typically accompanied with an inferior performance in the sense of the {\it maximal} risk over a suitably large class of models in a neighborhood of a fixed model, leading to a version of the {\it super-efficiency} phenomenon. We lay down our contributions below.
\begin{enumerate}
	\item We present general results on the asymptotic distribution of the averaged (pooled) estimator $\bar \theta$, both when $m$ is fixed and when allowed to increase as $N$ increases, in which case a normal distribution arises in the limit. Furthermore, in the latter case, the order of $m$, which affects the rate of convergence of $\bar \theta$, crucially depends on the bias of $\hat \theta_j$.  Pooling provides us with a novel way to construct a CI for $\theta_{0}$ whose length is shorter than that using $\hat \theta$ owing to the faster convergence rate involved: in fact, the ratio of the lengths of the CIs shrinks to 0. The calibration of the new CI involves normal quantiles, instead of quantiles of those of the non-standard limits that describe $\hat{\theta}$ asymptotically. Moreover, the variance $\sigma^2$ can be estimated empirically using the subsample-level estimates, whereas in the method involving $\hat{\theta}$, one is typically forced to impute values of several nuisance parameters that arise in the expression for $\sigma^2$ using estimates that can be quite unreliable. \newline

	\item The quantity that drives the possible gain by sample-splitting is the bias of the non-standard estimator. Hence, to obtain results on the rates of convergence (and asymptotic distribution) of the averaged estimator $\bar \theta$, we study the bias in prototypical non-standard problems: the LSE of a monotone regression function and its inverse (exhibiting cube-root asymptotics). The bias of the LSE or the maximum likelihood estimator in non-standard problems is \emph{hard} to compute because the usual Taylor expansion arguments that work in smooth function estimation fail in most non-standard problems. In particular, almost nothing seems to be known about the bias of the isotonic regression in the statistical literature. For the first time, we provide a non-trivial bound on the order of the bias of the monotone LSE under mild regularity assumptions. \newline

Furthermore, establishing the asymptotic normality of the pooled estimator in the monotone regression model, requires showing uniform integrability of certain powers of the normalized LSE as well as its inverse, pointwise. We are able to establish this property for all powers $p \geq 1$ in the general monotone regression model under a suitable `light-tail' assumption on the errors. As a consequence of these results we obtain upper bounds on the maximal  risk of the isotonic LSE and its inverse, over suitable classes of monotone functions. Although such bounds on the maximal risk are known for most nonparametric function estimators, this is the first instance of such a result in the general isotonic regression problem\footnote{Similar risk-bounds are presented in the special case of current status data in Theorem 11.3 in \cite{groeneboom2014nonparametric}; however, their derivation uses a special feature of the isotonic estimator in that particular model which is not true in the general scenario we consider, as discussed later in Remark~\ref{Groen-remark}.}. \newline
	%We show that (see Theorem~\ref{theo: biasDirect}) for the isotonic estimator, $$\mu(t_0) - \E[\hat \mu(t_0)] = O(N^{-7/15+\xi})$$ for any $\xi>0$, if the underlying regression function $\mu$ is assumed to have a  H\"olderian first derivative. We further study the bias of the inverse regression estimator and show that (see Theorem~\ref{theorem: BiasInverse}) $$\mu^{-1}(a) - \E [\hat \mu^{-1}(a)] = o(N^{-1/2}).$$
%\newline	
\item We present a rigorous study of a super-efficiency phenomenon that comes into play when using the pooled estimator in the context of estimating the inverse of an isotonic function. Let $\overline{\theta}$ denote the average of the $\hat{\mu}_{n,j}^{-1}(a)$'s, where $\hat{\mu}_{n,j}$ is the isotonic LSE from the $j$'th subsample and let $\theta_0 := t_0 \equiv \mu^{-1} (a)$ (see the discussion around~\eqref{eq:Chernoff-backward}). For a fixed $m$, we establish that for a \emph{fixed} $\mu$, 
\[ \E_{\mu}\left[N^{2/3}(\overline{\theta}  - \theta_0)^2\right] \rightarrow m^{-1/3}\,\mbox{Var}(\tilde{\kappa} \mathbb{Z}),\qquad \mbox{as } N \to \infty\,;\] see~\eqref{eq:Chernoff-backward}. On the other hand, we also show that for a suitably chosen (large enough)  class of models $\mathcal{M}_0$,  when $m \equiv m_n \rightarrow \infty$, 
$$\liminf_{N \to \infty}\,\sup_{\mu \in \mathcal{M}_0}\E_{\mu}\left[ N^{2/3}(\overline{\theta}  - \theta_0)^2 \right] = \infty,$$ 
whereas, for the global estimator $\hat{\theta}  \equiv \hat{\mu}_N^{-1}(a)$, 
$$\limsup_{N \to \infty}\,\sup_{\mu \in \mathcal{M}_0}\E_{\mu}\left[ N^{2/3}(\hat{\theta} - \theta_0)^2 \right] < \infty.$$ 
Thus, while the pooled estimator $\overline{\theta}$ can outperform the LSE under any \emph{fixed} model, its performance over a class of models is compromised relative to the isotonic LSE. The larger the number of splits ($m$), the better the performance under a fixed model, but the worse the performance over appropriately chosen classes of models. While we study the super-efficiency problem in a specific context, it is fairly clear (from computational evidence as well as heuristic considerations) that it will arise in other non-standard problems as well, e.g., the pointwise estimation of the monotone function itself. Our discoveries therefore serve as a cautionary tale that illustrates the potential pitfalls of using sample-splitting: the benefits from sample-splitting, both computational and in the sense of \emph{pointwise inference} may come at subtle costs. 
	
\end{enumerate}

The paper is organized as follows. In Section~\ref{sec:Fixed_m} we consider the case when $m$ is fixed and $n$ grows to infinity and study the behavior of the pooled estimator $\bar \theta$, while in Section~\ref{sec:Grow_m} we allow $m$ to grow with $n$. Section~\ref{sec:MonoReg} deals with a general monotone regression model where we derive bounds on the pointwise bias of both the isotonic LSE and its inverse, as well as $L_p$-risks. We use these results to study sample-splitting in various monotone function models in Section \ref{sec: appli}. Section \ref{sec: superEf} studies in some detail the super-efficiency phenomenon that comes into play in the isotonic regression setting under sample-splitting and compares and contrasts it with what transpires in kernel density estimation. The proofs of some of the main results are presented in Section~\ref{proof-section} and Appendix~\ref{Appendix} provides detailed coverage of additional technical material. 

Before we move on to the rest of the paper there is one point on which some clarity needs to be provided: in subsequent sections, the total sample size $N$ will be written as $m \times n$. Now, if $m$ is a fixed number, not all sample-sizes $N$ can be represented as a product of that form. To get around this difficulty, one can work with the understanding that we reduce our sample size from $N$ to $\tilde{N}:= m \times \lfloor N/m \rfloor$ (which is then renamed $N$) with the last few samples being discarded. Since finitely many are discarded, the resulting pooled estimate will be as precise in an asymptotic sense as the one based on the original $N$: in this (latter) case, one of the subsamples will have size less than $m$ but the contribution of the estimate from that subsample to the behavior of the pooled estimate is negligible in the long run. Similar considerations can be applied to the case of a growing $m$, so long as it is of a smaller order than $N$ which will always be the case in the sequel.  In this paper we work with the $\tilde{N}$ interpretation. 

\section{Fixed $m$ and growing $n$}\label{sec:Fixed_m}
Consider the setup of~\eqref{eq:LimDist}, where $\theta_0$ is the parameter of interest and let $\bar \theta$ be the pooled estimator as defined in~\eqref{eq:GlEst}. We start with a simple lemma that illustrates the statistical benefits of sample-splitting in the setting of~\eqref{eq:LimDist} when $n$ is large and $m$ is held fixed. 
\begin{lemma}\label{lem:m_Finite}
Suppose that~\eqref{eq:LimDist} holds where $G$ has mean zero and variance $\sigma^2 >0$. For $m$ fixed and $N = m \times n$, 
	\begin{equation}\label{eq:mFixedLimit}
		\sqrt{m} r_n (\bar \theta - \theta_0)  \stackrel{d}{\to}  H := m^{-1/2} (G_1 + G_2 + \ldots + G_m), \qquad \mbox{as } n \to \infty,
	\end{equation}
where $G_1, G_2, \ldots, G_m$ are i.i.d. $G$. Note that the limiting random variable $H$ has mean zero and variance $\sigma^2$. 
\end{lemma}
Compare the above result with the fact that if all $N$ data points were used together to obtain $\hat \theta$ we would have the limiting distribution given in~\eqref{eq:LimDist}. In particular, if $\{[r_n(\hat{\theta} - \theta_0)]^2\}_{n \ge 1}$ is uniformly integrable (which we will prove later for certain problems), we conclude that 
\[ \E [r_N^{2}(\hat \theta - \theta_0)^2] \rightarrow  \mathrm{Var}(G) \,,\qquad \mbox{as } N \to \infty,\]
while
\begin{equation}\label{eq:VarReduc}
\E \left[\frac{mr_n^2}{r_N^2}\,r_N^{2}(\bar \theta - \theta_0)^2 \right] \rightarrow \mathrm{Var}(G) \,,\qquad \mbox{as } N \to \infty,
\end{equation}
noting that $G$ and $H$ have the same variance. Thus, the asymptotic relative efficiency of $\bar \theta$ with respect to $\hat \theta$ is $m r_n^2/r_N^2$. For example, if $r_N = N^{\gamma}, \gamma < 1/2$, then using $\bar \theta$ gives us a reduction in asymptotic variance by a factor of $m^{1-2\gamma}$.  Hence, for estimating $\theta_0$, the pooled estimator $\bar \theta$ \emph{outperforms} $\hat \theta$.

\begin{remark}
If $\{[r_n(\hat{\theta}_{j} - \theta_0)]^2\}_{n \ge 1}$ is uniformly integrable then, $$\sigma_n^2 := \mathrm{Var}[r_n(\hat{\theta}_{j} - \theta_0)]  = r_n^2 \mathrm{Var}(\hat{\theta}_{j}) \rightarrow \sigma^2,\qquad \mbox{ as } n \to \infty,$$ for every $j = 1,\ldots,m$. As we have $m$ independent replicates from the distribution of $\hat \theta_j$, $\sigma^2$ can be approximated by 
\begin{equation}\label{eq:sigmahat}
\hat \sigma^2 := \frac{r_n^2}{m-1} \sum_{j=1}^m (\hat \theta_j - \bar \theta)^2.
\end{equation}
\end{remark}

\begin{remark}
For moderately large $m$ (e.g., $m \ge 30$) the $m$-fold convolution $H$ in~\eqref{eq:mFixedLimit} maybe well approximated by $N(0, \sigma^2).$ This yields a simple and natural way to construct an approximate CI for $\theta$:  
$$\left[\bar \theta - \frac{\hat \sigma}{r_n \sqrt{m}}  z_{\alpha/2}, \bar \theta +\frac{\hat \sigma}{r_n \sqrt{m}}  z_{\alpha/2}\right]$$ is an approximate $(1-\alpha)$ CI for $\theta$, where $z_\alpha$ is the $(1-\alpha)$'th quantile of the standard normal distribution. Note that, we have completely by-passed the direct estimation of the problematic nuisance parameter $\sigma^2$.
\end{remark}

\begin{remark}
The normal approximation can be avoided at the expense of  simulating the distribution of the limiting $m$-fold convolution $H$ in~\eqref{eq:mFixedLimit} and estimating the appropriate quantiles. For example, when $G$ is a scaled Chernoff's distribution, i.e., $G \equiv_d \sigma\,\tilde{Z}$, where $\tilde{Z}$ is the Chernoff random variable scaled by its standard deviation, we would simulate the distribution of $\sigma \sum_{j=1}^m \tilde{Z}_j/\sqrt{m}$, where $\tilde{Z}_1,\ldots, \tilde{Z}_m$ are i.i.d.~$\tilde{Z}$. Since it is easy to generate from $\tilde{Z}_j$ (see e.g.,~\cite{GW01}), and $\sigma^2$ can be estimated as shown in~\eqref{eq:sigmahat}, fairly accurate empirical quantiles of the exact limit can be generated. 
\end{remark}

\section{Letting $m$ grow with $n$: asymptotic considerations}\label{sec:Grow_m}
In this section, we derive the asymptotic distribution of $\sqrt{m} r_n(\bar{\theta} - \theta_0)$ under certain conditions, as $m \rightarrow \infty$. We first introduce some notation. To highlight the dependence on $n$, we write $m \equiv m_n$, $\hat{\theta}_{j} \equiv \hat{\theta}_{n,j}$ and $\bar{\theta} = \bar{\theta}_{m_n}$. Consider the triangular array of i.i.d.~random variables $\{\xi_{n,1}, \xi_{n,2}, \ldots, \xi_{n,m_n}\}_{n\ge1}$ where 
$\xi_{n,j} := r_n(\hat{\theta}_{n,j} - \theta_0).$ Let $b_{n} := \E(\xi_{n,1}) = r_n(\theta_n - \theta_0)$ where $ \theta_n := \E(\hat{\theta}_{n,1})$ is assumed to be well-defined. The following theorem is proved in Section~\ref{proof:thm:m_n_infty}.
\begin{theorem}\label{thm:m_n_infty}
Suppose that~\eqref{eq:LimDist} holds where $G$ has mean zero and variance $\sigma^2 >0$.  Also, suppose that $b_n = O(c_n^{-1}),$ where $c_n \to \infty$ as $n \to \infty,$ and that the sequence $\{\xi_{n,1}^2\}$ is uniformly integrable. Then,  as $n\to\infty$, 
\begin{enumerate}
	\item[(i)] for any $m_n \rightarrow \infty$ such that $m_n = o(c_n^2)$, $$\sqrt{m_n}r_n (\bar{\theta}_{m_n} - \theta_0) \stackrel{d}{\to} N(0, \sigma^2);$$
	 \item[(ii)] if $m_n \sim O(c_n^2)$, and furthermore $\sqrt{m_n}\,b_n \rightarrow \tau$, then 
\begin{equation}\label{eq:Bias}
\sqrt{m_n}r_n (\overline{\theta}_{m_n} - \theta_0) \stackrel{d}{\to} N(\tau, \sigma^2). \vspace{0.15in}
\end{equation}
\end{enumerate}
\end{theorem}

\begin{remark}[Gains from sample-splitting: ``divide to conquer'']
The pooled estimator $\bar{\theta}_{m_n}$ is more effective than $\hat{\theta}_N$, when its convergence rate exceeds that of the latter, i.e., $$ \frac{r_N}{\sqrt{m_n} r_n} \rightarrow 0 \Leftrightarrow \frac{r_N/r_n}{m_n^{1/2}} \rightarrow 0;$$ thus, if $r_N = N^{\alpha}$, using $N = n \times m_n$, this requires $\alpha < 1/2$. In other words, acceleration is only possible if the initial estimator has a slower convergence rate  than the parametric rate. 
\end{remark}
\begin{remark}[Choice of $m_n$] 
\label{optimal-mn}
As above, let $r_N = N^{\alpha}$ with $\alpha < 1/2$, and let $c_n = n^{\phi}$. Choosing $m_n = n^{2\phi - \delta}$, with $0< \delta < 2\phi$, so that $m_n = o(c_n^2)$, we have $\sqrt{m_n}\,r_n = n^{\phi - \delta/2 + \alpha}$. Using $m_n \times n = N$, we get $n = N^{1/(2\phi - \delta + 1)}$. The convergence rate of the pooled estimator in terms of the total sample size is therefore $N^{( \phi - \delta/2+\alpha )/2(\phi- \delta/2 + 1/2)}$. Since $\alpha < 1/2$, this rate is strictly less than $N^{1/2}$. Next, the improvement in the convergence rate is given by 
\[  \frac{ \phi - \delta/2+\alpha }{2(\phi- \delta/2 + 1/2)} - \alpha = 2\left(\frac{1}{2} - \alpha \right)\,\frac{\phi - \delta/2}{\phi - \delta/2 + 1/2}\,,\]
which is monotone decreasing in $\delta$. This means that smaller values of $\delta$, corresponding to larger values of $m_n = N^{(2\phi - \delta)/(2\phi - \delta + 1)}$ give greater improvements in the convergence rate. In the situation of conclusion (ii) of the above theorem, when $\delta = 0$ and $m_n = O(c_n^2)$, we get the maximal convergence rate: $N^{(\alpha + \phi)/2(\phi + 1/2)}$.

To get the best possible rate out of sample-splitting, ideally, we would like to get hold of the \emph{optimal} value of $c_n$, i.e., we would want $b_n = O(c_n^{-1})$ but not $o(c_n^{-1})$. The optimal $c_n$ might, of course, be difficult to obtain in a particular application; however, sub-optimal $c_n$'s will also improve the rate of convergence, albeit not to the best possible extent. 
\end{remark}
From Theorem \ref{thm:m_n_infty} we see that the two key challenges to establishing the asymptotic normality of the pooled estimator are: (a) establishing uniform integrability as desired above, and, (b) determining an order for the bias $b_n$. In the following sections we consider the example of monotone regression and address (a) and (b) for the isotonic MLE and its inverse.  

\section{Isotonic regression: uniform integrability and bias}\label{sec:MonoReg}
In this section, we consider the example of monotone regression (a prototypical example of non-standard asymptotics) and establish the uniform integrability of the isotonic LSE and its inverse as well as an order for the point-wise bias. While this section is useful in applying Theorem~\ref{thm:m_n_infty} to derive the rate of convergence and asymptotic distribution of the pooled estimator in isotonic regression (Section \ref{sec: appli}), we believe it is also of independent interest. To the best of our knowledge, Theorems \ref{theorem: BiasInverse} and \ref{theo: biasDirect} constitute the first attempts in the literature to study the order of the bias of the monotone LSE under mild regularity assumptions and may well have implications beyond the sample-splitting methodology considered in this paper. Our formal treatment is developed in the framework of~\cite{Durot08} which considers a general monotone non-increasing regression model described below. The results, of course, extend immediately to the nondecreasing case.

We observe independent copies $\{W_i \equiv (X_i,Y_i): i=1,\ldots, n\}$ of a bivariate random variable $(X,Y)\in[0,1]\times\R$. We aim at estimating the regression function $\mu$ defined by $\mu(x)=\E(Y|X=x)$ for $x\in[0,1]$, under the constraint that  it is nonincreasing on $[0,1]$. Alternatively, we may be interested in estimating the inverse function $\mu^{-1}$. With $\eps=Y-\mu(X)$ we define $v(x):=\E(\eps_i^2|X_i=x)$ for all $x\in[0,1]$ and we make the following assumptions.

\begin{description}\item{\monorandom} $\mu$ is differentiable and decreasing on $[0,1]$ 
with $\inf_t| \mu'(t)|>0$ and $\sup_t| \mu '(t)|<\infty.$ 
\item{\densityUI} $X$ has a density $f$ which is bounded and
bounded away from zero. 
\item{\sigmamin} There exists $c_0>0$ such that $v^2(t)\geq c_0(t\wedge(1-t))$ 
for all $t\in[0,1]$.
\item{\momentq} There exist $K>0$ and $\alpha>0$ such that $\E\left(e^{\theta \eps}|X\right)\leq K\exp(\theta^2\alpha)\ a.e.$ for all $\theta\in\R$,
\end{description}
Assumption {\sigmamin} is slightly less restrictive than the usual assumption of a bounded variance function $v$ away from zero and allows us to handle, for example, the current status model in Subsection \ref{currstat}. Assumption {\momentq} is fulfilled for instance if the conditional distribution of $\epsilon$ given $X$ is sub-Gaussian and the variance function $v$ is bounded, or if  $\epsilon$ is bounded.

\subsection{The isotonic LSE of $\mu$ and the inverse estimator}\label{sec: LSE}
We start with an exposition of the characterization of the least-squares estimator (LSE) of $\mu$ and its inverse under the monotonicity constraint. With  $X_{(1)}<\dots<X_{(n)}$ the order statistics corresponding to $X_1,\dots,X_n$, and $Y_{(i)}$ the observation corresponding to $X_{(i)}$,  let $\Lambda_n$ be the piecewise-linear process on $[0,1]$ such that 
\begin{equation}\label{Lambdan}
\Lambda_n\left(\frac{i}{n}\right)=\frac{1}{n}\sum_{j\leq i}Y_{(j)}
\end{equation}
for all $i\in\{0,\dots,n\}$, where we set $\sum_{j\leq 0}Y_{(j)}=0$. 
%and linearly interpolate between those points. 
Let $\hat\lambda_n$ be the left-hand slope of the least 
concave majorant of $\Lambda_n$ with $\hat\lambda_{n}(0)=\lim_{t\downarrow 0}\hat\lambda_{n}(t)$. It is well known that a monotone $\hat\mu_{n}$ is an LSE if and only if it satisfies
\begin{equation}\label{brunk0}
\hat \mu_n(X_{(i)})=\hat\lambda_n(i/n)
\end{equation}
for all $ i=1,\dots,n.$ In the sequel, we consider the piecewise-constant left-continuous LSE $\hat\mu_{n}$ that is constant on  the intervals $[0, X_{(1)}]$, $(X_{(n)},1]$ and $(X_{(i-1)},X_{(i)}]$ for all $i=2,\dots,n-1$.

Now, recall that for every nonincreasing left-continuous function $h:[0,1]\to\R$, the generalized inverse of $h$  is defined as: for every $a\in\R$, $h^{-1}(a)$ is the greatest  $t\in[0,1]$ that satisfies $h(t)\geq a$, with the convention that the supremum of an empty set is zero. In the sequel, we consider the generalized inverse  $\hat\mu_{n}^{-1}$ of  $\hat\mu_{n}$ as an estimator for  $\mu^{-1}$.

%%%%%%%%%%%%%%%%%%%%%%%%%%%%%%
\subsection{Uniform integrability and bias of the direct and inverse estimators}\label{sec: intdirect}
%%%%%%%%%%%%%%%%%%%%%%%%%%%%%%
In this subsection we provide bounds on the absolute centered moments of the isotonic LSE and its inverse. These results will imply uniform integrability of the corresponding estimates. We also establish the order of the bias for both the LSE and its inverse. 
First, we consider the absolute centered moments. The proofs of the main results in this subsection are given in Section~\ref{sec: randomU}.
\begin{theorem}\label{theo: UnifIntRD}
Assume {\momentexpo}, $X$ has a density function $f$, $\mu$ is nonincreasing, and there exist positive numbers $A_{1},\dots,A_{5}$ such that 
$A_{1}<|\mu'(t)|<A_{2}\mbox{, }A_{3}<f(t)<A_{4}\mbox{ and } |\mu(t)|\leq A_{5}$
for all $t\in[0,1]$. Then, for any $p\geq 1$, there exists $K_{p}>0$  that depends only on $p,A_{1},\dots,A_{5},K,\alpha$, where $K$ and $\alpha$ are taken from  {\momentexpo}, such that for all $n$,
\begin{enumerate}
\item $\E\left(|\hat \mu_{n}^{-1}(a) -\mu^{-1}(a)|^p\right)\leq K_{p}n^{-p/3}$
for all $a\in\R$,
\item
$\E\left(|\hat \mu_{n}(t) - \mu(t)|^p\right)\leq K_{p}n^{-p/3}$
for all $t\in[n^{-1/3},1-n^{-1/3}]$. $\vspace{0.15in}$
\end{enumerate}
\end{theorem}

\begin{remark}
\label{Groen-remark} 
Inequalities (11.32) and (11.33) in \cite{groeneboom2014nonparametric} are special cases of the second assertion of Theorem~\ref{theo: UnifIntRD} above since the current status model is a special case of the general regression model (see Section~\ref{currstat}). However,  the inequalities in \cite{groeneboom2014nonparametric} hold for all $t$, whereas the corresponding inequality above holds only for $t$  in a restricted interval. This is due to a very specific feature of the estimator in the current status model: it has the same range as the estimated function since both of them are distribution functions. In particular,  the estimator is consistent at the boundaries in the current status model, whereas it is not in the general regression model. Hence, the strategy of proof  in \cite{groeneboom2014nonparametric} does not extend to the general regression model. Specifically, the proof in 
\cite{groeneboom2014nonparametric} is based solely on an exponential inequality for the tail probabilities of the inverse estimator given in Theorem 11.3 of that book, whereas our proof  is based on two exponential inequalities, see Section~\ref{sec: expobounds}:  Lemma \ref{lem: expoboundinvm}  extends Theorem 11.3 in \cite{groeneboom2014nonparametric} to our general setting, and Lemma \ref{lem: expoboundU+-} gives a sharper exponential inequality  for the case when the inverse estimator is computed at some point, $a$, that does not belong to the range of $\mu$. 
\end{remark}

A direct corollary (below) to Theorem \ref{theo: UnifIntRD} is an upper  bound on the maximal  risk of the two estimators discussed above. Although such bounds on the maximal risk over  suitable classes of functions are known for most nonparametric function estimators, this is the \emph{first instance} for such a result in the context of isotonic regression.

\begin{corollary}
\label{maximal-risk} 
Let $A_{1}, A_{2}, A_{3}$ be positive numbers. Let $\F_1$  be the class of nonincreasing functions $\mu$ on $[0,1]$ such that $A_{1}<|\mu'(t)|<A_{2}$ and $|\mu(t)|\leq A_{3}$ for all $t\in[0,1]$. If {\densityUI} and {\momentexpo} hold, then for any $p\geq 1$, there exists $K_{p}>0$ such that
\begin{enumerate}
\item $\limsup \limits_{n \to \infty} \;\sup_{\mu \in  \F_1} n^{p/3}\, \E_{\mu}\left(|\hat \mu_{n}^{-1}(a) -\mu^{-1}(a)|^p\right)\leq K_{p}$
for all fixed $a\in\R$, 
\item
$\limsup\limits_{n \to \infty} \; \sup_{\mu \in \F_1} n^{p/3}\, \E_{\mu}\left(|\hat \mu_{n}(t) - \mu(t)|^p\right)\leq K_{p}$
for all fixed $t\in(0,1)$.
\end{enumerate}
\end{corollary}

We next consider the order of the bias. Tackling the bias requires imposing additional smoothness assumptions on the underlying parameters of the problem. Precisely,  we assume  for some of our results that $v^2$ has a bounded second derivative on $[0,1]$, that
\begin{equation}\label{eq: holder}
|\mu'(x)-\mu'(y)|\leq C|x-y|^s, \qquad \mbox{ for all } x,y\in[0,1],
\end{equation}
for some $C>0$ and $s>0$ (where bounds on $s$ will be specified precisely while stating the actual results); and, instead of {\densityUI}, the more restrictive assumption:
\begin{description}
\item{\density} The density $f$ of $X$  is 
bounded away from zero with a bounded first derivative on $[0,1]$. 
\end{description}
\begin{theorem}\label{theorem: BiasInverse}
Assume {\monorandom}, {\density}, {\sigmamin} and {\momentq}.  Assume, furthermore, that  $v^2$ has a bounded second derivative on $[0,1]$ and $\mu$ satisfies \eqref{eq: holder} for some $C>0$ and $s>3/4$.   For an arbitrary constant $K>0$ we then have
\begin{equation}\notag
\E\left(\hat \mu_{n}^{-1}(a)-\mu^{-1}(a)\right)=o(n^{-1/2})
\end{equation}
where the small-$o$ term is uniform in $a\in [\mu(1)+Kn^{-1/6}\log n,\mu(0)-Kn^{-1/6}\log n]$.
\end{theorem}
Now, consider the bias of the direct estimator. Ideally, one would like to prove that $\E\left(\hat \mu_{n}(t)- \mu(t)\right)=o(n^{-1/2})
$ uniformly in $t\in [Kn^{-1/6}\log n,1-Kn^{-1/6}\log n]$, with an arbitrary $K>0$ that does not depend on $n$. Unfortunately, we are only able to obtain a somewhat less precise bound. We also require a higher degree of smoothness $s=1$ on $\mu'$ than needed for dealing with the inverse function\footnote{For smaller values of $s$, we obtain an even larger bound but this is not discussed any further in the paper.}. 
\begin{theorem}\label{theo: biasDirect}
 Assume {\monorandom}, {\density}, {\sigmamin}, {\momentq},  $v^2$ has a bounded second derivative on $[0,1]$ and $\mu$ satisfies \eqref{eq: holder} for some $C>0$ and $s=1$.  For an arbitrary fixed interval $[c_{1},c_{2}]\subset(0,1)$, we have
 $$\E\left(\hat \mu_{n}(t)-\mu(t)\right)=O(n^{-7/15+\zeta})$$
with an arbitrary $\zeta>0$, where the big-$O$ term is uniform in $t\in [c_{1},c_{2}]$.
\end{theorem}

\section{Applications to sample-splitting in monotone function models} \label{sec: appli}

\subsection{Simple isotonic regression model}  
{\bf Function estimation at a point: } We consider $N$ i.i.d.~data $\{X_i, Y_i\}_{i=1}^N$ from the simple isotonic regression model with error independent of covariate as considered in (\ref{eq:IsoReg}) and recall the notation used therein. The parameter of interest is $\theta_0 \equiv \mu(t_0)$ which is estimated by  $$\bar \theta_{m_n} =  \frac{1}{m_n} \sum_{j=1}^{m_n} \hat{\mu}_{n,j}(t_0),$$ $\hat{\mu}_{n,j}$ being the isotonic LSE computed from the $j$-th split-sample. For this problem, the function $v^2(t) \equiv v^2 > 0$. Under (a subset of) the assumptions on the parameters of the model made in Theorem~\ref{theo: biasDirect}, convergence in law to 
Chernoff's distribution (recall (\ref{eq:Chernoff-forward})) holds. To apply Theorem~\ref{thm:m_n_infty}, we need to show that: $(a)$ $n^{1/3}(\theta_n - \mu(t_0)) = O(n^{-\phi})$ (here $\theta_n = \E[\hat{\mu}_{n,1}(t_0)]$) for some $\phi >0$, and $(b)$ the uniform integrability of the sequence $\{n^{2/3}(\hat{\mu}_{n,1}(t_0) - \mu(t_0))^2\}_{n \ge 1}$. 

Now, $(b)$ is a direct consequence of Theorem \ref{theo: UnifIntRD} applied with any $p > 2$. As far as $(a)$ is concerned, by Theorem \ref{theo: biasDirect}, we know that the desired condition in $(a)$ is satisfied for $s =1$ in~\eqref{eq: holder} for any fixed $t_0 \in (0,1)$, by taking $\phi = (7/15 - 1/3) - \zeta = (2/15 - \zeta)$ 
where $\zeta > 0$ can be taken to be arbitrarily small. From Remark~\ref{optimal-mn}, choosing $m_n = n^{2\phi - \delta} = n^{4/15 - 2\zeta - \delta}$ for a small enough 
$0 < \delta < 2\phi$, we conclude that with $\sigma^2 = {\kappa}^2\,\mbox{Var}(\mathbb{Z})$, we have
\begin{equation}
\label{pooled-iso-convergence}
N^{( {7/15} - \zeta - \delta/2)/( {19/15 - 2\zeta - \delta})}(\overline{\theta}_{m_n} - \theta_0) \stackrel{d}{\to} N(0, \sigma^2)\,.
\end{equation}

\noindent
{\bf Inverse function estimation at a point}:  Consider the same set-up as in the above problem. We now consider estimation of $\mu^{-1}(a)$ via the inverse isotonic LSE {under the assumptions of Theorem \ref{theorem: BiasInverse}}. The behavior of the isotonic estimator $\hat{\mu}_N$ based on the entire data of size $N$ is given in \eqref{eq:Chernoff-backward}. To apply Theorem \ref{thm:m_n_infty}, we need to show that: $(a)$ $n^{1/3}(\theta_n - \mu^{-1}(a)) = O(n^{-\phi})$ (here $\theta_n = \E[\hat{\mu}_{n,1}^{-1}(a)]$) for some $\phi >0$, and $(b)$ the uniform integrability of the sequence $\{n^{2/3}(\hat{\mu}_{n,1}^{-1}(a) - \mu^{-1}(a))^2\}_{n \ge 1}$. 

In this case, $(b)$ is a direct consequence of Theorem \ref{theo: UnifIntRD} applied with any $p > 2$. As far as $(a)$ is concerned, by Theorem \ref{theorem: BiasInverse}, we know that the desired condition in $(a)$ is satisfied for  {$s >3/4$} in~\eqref{eq: holder} for any fixed $a$ in the interior of the range of $\mu$, by taking $\phi = (1/2 - 1/3) = 1/6$. From Remark~\ref{optimal-mn}, choosing $m_n = n^{2\phi} = n^{ {1/3}}$ (for the inverse function estimation problem we are actually in the situation of conclusion (ii) of Theorem \ref{thm:m_n_infty} with $\tau = 0$), we conclude that: 
\begin{equation}
\label{pooled-iso-convergence-inv}
N^{(1/3 + 1/6)/2(1/6+1/2)}(\overline{\theta}_{m_n} - \theta_0) \equiv N^{3/8}(\overline{\theta}_{m_n} - \theta_0)  \stackrel{d}{\to} N(0, \widetilde{\sigma}^2)\,,
\end{equation}
where $\widetilde{\sigma}^2 = \tilde{\kappa}^2\,\mbox{Var}(\mathbb{Z})$. The pooled estimator, therefore, has a convergence rate of $N^{3/8}$. 
\begin{remark}
\label{integrated-bias-remark} 
Note that the order of the bias obtained in the forward problem (Theorem \ref{theo: biasDirect}) is slower than that obtained in the inverse problem 
(Theorem \ref{theorem: BiasInverse})and comes at the expense of increased smoothness ($s=1$) compared to Theorem \ref{theorem: BiasInverse} (where we assume $s > 3/4$). This seems to be, at least partly, an artifact of our approach where we start from the characterization of the inverse estimator as our starting point and derive results for the forward problem from those in the inverse problem through the switching relationship. Ideally, one would want to derive the same order for the bias in both forward and inverse problems for a fixed degree of H\"{o}lder smoothness on $\mu'$.

Next, even for the inverse problem, it is not clear at this point whether the order of the bias obtained in Theorem \ref{theorem: BiasInverse} is optimal, i.e., the best possible one under the assumed smoothness. It is conceivable that when $s > 3/4$ the exact order of the bias is smaller than the obtained $o(n^{-1/2})$ rate from Theorem \ref{theorem: BiasInverse}. A smaller bias would allow a faster rate of convergence than $N^{3/8}$ through an appropriate choice of $m_n$. A complete resolution of the bias problem would require characterizing the optimal order of the bias in the isotonic regression problem as a function of $s$ (with larger $s$'s corresponding to smaller orders), but this is outside the scope of this paper. It is, however, worth reiterating that Theorems \ref{theorem: BiasInverse} and \ref{theo: biasDirect} are the first systematic attempts in the literature to quantify the bias of isotonic estimators. 
\end{remark}

\subsection{The current status model}
\label{currstat}   Our framework covers the important case of the current status model, which has found extensive applications in epidemiology and biomedicine. The problem is to estimate the distribution function $F_T$ of a failure time $T\geq 0$ on $[0,1]$, based on observing $n$ independent copies of the censored  pair $(X,\BBone_{T\leq X})$. Here, $X\in[0,1]$ is the observation time independent of $T$, and $\BBone_{T\leq X}$ stipulates whether or not the failure has occurred before time $X$. Then,
$$F_T(x)=\P(T\leq x)=\E(\BBone_{T\leq X}|X=x)$$
for all $x\in[0,1]$. This falls in the general framework of Section~\ref{sec:MonoReg} with $Y=-\BBone_{T\leq X}$ and $\mu=-F_T$, which is nonincreasing. It turns out that the nonparametric maximum likelihood estimator (MLE) of $F_T$ is precisely $-\hat \mu_n$ where $\hat \mu_n$ is the LSE from Section~\ref{sec: LSE},  see \cite{groeneboom1992information}. We present results separately for the current status model in the following theorem, proved in Section~\ref{proof:theo: currentstatus}. 

\begin{theorem}\label{theo: currentstatus}
Assume that we observe $n$ independent copies of $(X,\BBone_{T\leq X})$, where $X\in[0,1]$ is independent of $T\geq 0$. Assume that $T$ has a density function $f_{T}$ that is bounded away from both zero and infinity on $[0,1]$, and that $X$ has a density function $f$ on $[0,1]$ that is bounded away from zero and has a bounded first derivative on $[0,1]$. With $\hat F_{Tn}$ the MLE of the distribution function $F_{T}$ of $T$, and $\hat F_{Tn}^{-1}$  the corresponding quantile function, we have:
\begin{enumerate}
\item For any $p\geq 1$, there exists $K_{p}>0$ such that for all $n$,
$$\E\left(|\hat F_{Tn}(t)-F_{T}(t)|^p\right)\leq K_{p}n^{-p/3}\mbox{ 
for all }t\in[n^{-1/3},1-n^{-1/3}]$$
and 
$$\E\left(|\hat F_{Tn}^{-1}(a)-F_{T}^{-1}(a)|^p\right)\leq K_{p}n^{-p/3}\mbox{ 
for all }a\in\R.$$
\item If moreover,  $f_T$ has a bounded first derivative, then with $K>0$, $c_1>0$, $c_{2}<1$, and $\phi>0$ arbitrary constants, 
$$\E\left(\hat F_{Tn}^{-1}(a)-F_{T}^{-1}(a)\right)=o(n^{-1/2})$$
uniformly for all $a\in [Kn^{-1/6}\log n,1-Kn^{-1/6}\log n]$ and
 $$\E\left(\hat F_{Tn}(t)-F_T(t)\right)=O(n^{-7/15+\phi})$$
 uniformly for all  $t\in [c_{1},c_{2}]$.  
 \item Now, let $\hat{F}_{TN}$ denote the MLE based on $N = m_n \times n$ observations from the current status model, $\hat{F}_{Tn}^{(j)}$ the MLE from the $j$'th subsample and $\overline{F}_{m_n}$ the pooled isotonic estimator obtained by averaging the $\hat{F}_{Tn}^{(j)}$s. If $f_T$ has a bounded first derivative, then for all $\zeta, \delta > 0$, sufficiently small, and any 
 $0 < t < 1$, 
 \[ N^{(7/15 - \zeta - \delta/2)/(19/15 - 2\zeta - \delta)}(\overline{F}_{m_n}(t) - F(t)) \stackrel{d}\rightarrow N(0, \sigma^2)\,,\]
where $\sigma^2 = \{4\,F_T(t)(1-F_T(t))f_T(t)/f(t)\}^{2/3}\,\mathrm{Var}(\mathbb{Z})$.

Moreover, for any $a \in(0,1)$, with $\overline{\theta}_{m_n} $ the pooled  estimator obtained by averaging the $ (\hat{F}_{Tn}^{(j)})^{-1}(a)$s, 
\[ N^{3/8}(\overline{\theta}_{m_n} - F_T^{-1}(a)) \stackrel{d}{\rightarrow} N(0, \tilde{\sigma}^2)\,,\]
where $\tilde\sigma^2 = \{4\,a(1-a)/f_T^{'}(t_a)^2 f(t_a)\}^{2/3}\,\mathrm{Var}(\mathbb{Z})$ with $t_a=F_T^{-1}(a)$.
\end{enumerate} 
\end{theorem}

\section{Sample-splitting and the super-efficiency phenomenon} \label{sec: superEf}
The variance reduction accomplished by sample-splitting (see~\eqref{eq:VarReduc}) for estimating a fixed monotone function, or its inverse, at a given point comes at a price. We show in this section, in the context of the inverse problem, that though a larger number of splits ($m$) brings about greater reduction in the variance for a fixed function, the performance of the pooled estimator in a \emph{uniform sense}, over an appropriately large class of functions, deteriorates in comparison to the global estimator as $m$ increases. This can be viewed as a \emph{super-efficiency} phenomenon: a trade-off between point wise performance and performance in a uniform sense. We elaborate below. 

\subsection{Super-efficiency of the pooled estimator}

Fix a continuous monotone (nonincreasing) function $\mu_0$ on $[0,1]$ that is continuously differentiable on $[0,1]$ with $0 < c < |\mu_0^{'}(t)| < d < \infty$ for all $t \in [0,1]$. Let $x_0 \in (0,1)$. Define a neighborhood $\mathcal{M}_0$ of $\mu_0$ as the class of all continuous non-increasing functions $\mu$ on $[0,1]$  that are continuously differentiable on $[0,1]$, that coincide with $\mu_0$ outside of $(x_0 - \epsilon_0, x_0 + \epsilon_0)$ for some (small) $\epsilon_0>0$, and such that $0 < c < |\mu'(t)| < d < \infty$ for all $t \in [0,1]$. Now, consider $N$ i.i.d.~observations $\{Y_i, X_i\}_{i=1}^N$ from the model: $$Y = \mu_0(X) + \epsilon,$$ where $X \sim \mbox{Uniform}(0,1)$ is independent of $\epsilon \sim N(0, v^2)$. Let $\hat{\theta}_N$ denote the isotonic estimate of $\theta_0: = 
\mu_0^{-1}(a)$ as considered before. We know that as $N \rightarrow \infty$, 
\begin{equation}\label{eq:IsoRegLimDist}
N^{1/3}\,(\hat{\theta}_N - \theta_0) \stackrel{d}{\rightarrow} G,
\end{equation}
where $G =_d \tilde{\kappa}\,\mathbb{Z}$, $\mathbb{Z}$ being the Chernoff random variable, and $\tilde{\kappa} > 0$ being a constant. If we split $N$ as $m \times n$, where $m$ is a fixed integer, then as $N \rightarrow 
\infty$, Lemma \ref{lem:m_Finite} tells us that  
\[ N^{1/3}(\overline{\theta}_{m} - \theta_0) \stackrel{d}{\rightarrow} m^{-1/6} H,\]
where $\overline{\theta}_{m}$ is the pooled estimator and $H$ has the same variance as $G$. By Theorem \ref {theo: UnifIntRD} we have uniform integrability under $\mu_0$ and conclude that: 
\begin{equation}\label{eq:MSEGlobalEst}
\E_{\mu_0}\left[N^{2/3}(\hat{\theta}_N - \theta_0)^2\right] \rightarrow \mathrm{Var}(G), \qquad \mbox{as } N \to \infty,
\end{equation}
while
\begin{equation}\label{eq:MSEPooledEst} 
\E_{\mu_0}\left[N^{2/3}(\overline{\theta}_{m} - \theta_0)^2\right] \rightarrow m^{-1/3}\,\mathrm{Var}(G),\qquad \mbox{as } N \to \infty\,.
\end{equation}
Hence, for estimating $\theta_0 = \mu_0^{-1}(a)$, the pooled estimator \emph{outperforms} the isotonic regression estimator.

We now focus on comparing the performance of the two estimators over the class $\mathcal{M}_0$. In this regard we have the following theorem, proved in Section \ref{supeff-proof}.
\begin{theorem}\label{superefficient}
Let
\begin{equation}\label{eq:IsoRegUnifBd}
E := \limsup_{N \to \infty}\,\sup_{\mu \in \mathcal{M}_0}\E_{\mu}\left[N^{2/3}(\hat{\theta}_N - \theta_0)^2 \right],
\end{equation}
and 
\[  E_m: =  \liminf_{N \to \infty}\,\sup_{\mu \in \mathcal{M}_0}\E_{\mu}\left[ N^{2/3}(\overline{\theta}_{m} - \theta_0)^2 \right],\]
where the subscript $m$ indicates that the maximal risk of the $m$-fold pooled estimator ($m$ fixed) is being considered. Then $E < \infty$ while $E_m \geq m^{2/3}\,c_0,$ for some $c_0 > 0$. When $m = m_n$ diverges to infinity, 
$$\liminf_{N \to \infty}\,\sup_{\mu \in \mathcal{M}_0}\E_{\mu}\left[ N^{2/3}(\overline{\theta}_{m_n} - \theta_0)^2 \right] = \infty \,.$$ 
\end{theorem} 
Therefore, from Theorem~\ref{superefficient} it follows that the asymptotic maximal risk of the pooled estimator diverges to $\infty$ (at least) at rate $m^{2/3}$. Thus, the better off we are in a pointwise sense with the pooled estimator, the worse off we are in the uniform sense over the class of functions~$\mathcal{M}_0$.

{\small 
\begin{table}
\begin{tabular}{|r|rrrrrrr|}
\hline
$(n,m)$ &  5  & 10& 15 & 30 & 45 & 60 & 90 \\
\hline
%    50      &  1.23 & 1.58 & 1.52 & 1.53 & 1.38 & 1.48 & 1.37 \\
 50&        1.67  &  1.71  &  1.90&    1.66 &   1.57 &   1.65  &  1.17\\
100  &    1.31  &  1.76  &  2.21  &  2.29  &  2.16   & 2.46 &   2.33\\
200  &    1.75  &  2.06  &  2.42   & 2.81  &  2.58  &  3.16 &   3.39\\
500   &   1.70 &   2.13   & 2.12 &   2.80  &  3.16  &  3.59  &  4.11\\
1000  &  1.46  &  2.04  &  2.46  &  2.88  &  3.60  &  3.51 &   4.31\\
3000 &   1.63  &  2.12  &  2.33  &  3.11  &  4.15  &  3.84  &  3.69\\
10000  & 1.75  &  2.11  &  2.70  &  2.86  &  3.31 &   5.08  &  5.18\\ \hline
\end{tabular}
\begin{tabular}{|rrrrrrr|}
\hline
5 & 10 & 15 & 30 & 45 & 60 & 90 \\
\hline
       1.47  &  1.21   & 0.94 &   0.70  &  0.55 &   0.54 &   0.39\\
   1.04  &  0.97  &  0.90 &   0.59  &  0.47  &  0.40  &  0.31\\
    1.03  &  0.94  &  0.76  &  0.68  &  0.42   & 0.38  &  0.29\\
     1.01  &  0.90 &   0.69  &  0.54  &  0.44  &  0.34  &  0.24\\
    1.16 &   0.88  &  0.66  &  0.52 &   0.36 &   0.34  &  0.24\\
   1.09  &  0.87  &  0.75   & 0.43 &   0.40   & 0.31  &  0.21\\
   0.94  &  0.79 &   0.80  &  0.43  &  0.33   & 0.31  &  0.23\\          
     \hline
\end{tabular}
\caption{Ratios of the (estimated) mean squared errors $\frac{\E \left[ (\hat{\mu}_N^{-1}(a) - \theta_0)^2\right]}{\E \left[ (\overline{\theta}_m - \theta_0)^2\right]}$ comparing the performance of the pooled estimator $\overline{\theta}_m$ with the global estimator $\hat{\mu}_N^{-1}$ as $n$ and $m$ change for the model: $Y = \mu(X) + \eps$, $X \sim $ Unif$(0,1)$, $\eps \sim N(0,0.2^2)$,  and $a = 0.5$, with (i) $\mu(x) = x$,  and (ii) $\mu(x) = \mu_n(x) = x + n^{-1/3}B(n^{1/3}(x- x_0))$ with $B(u) = 2^{-1}(1-(|u|-1)^2)^2 \BBone_{\{|u| \le 2\}}$. For both (i) and (ii), $\theta_0 \equiv \mu^{-1}(a) = 0.5$.}
\label{table:Ratio}
\end{table}
}
Table~\ref{table:Ratio} gives the ratios of the (estimated) mean squared errors ${\E \left[ (\hat{\mu}_N^{-1}(a) - \theta_0)^2\right]}/{\E \left[ (\overline{\theta}_m - \theta_0)^2\right]}$  comparing the performance of the pooled estimator $\overline{\theta}_m$ with the global estimator $\hat{\mu}_N^{-1}(a) $ as $n$ and $m$ change for two different models, which are described in the caption to the table. For the first model (left table) we fix the regression function at $\mu(x) = x$ and let $N \to \infty$ and find that the pooled estimator has superior performance to the global estimator as $m$ (and $n$) grows. The ratio of the mean squared errors is generally 
close to $m^{1/3}$, as per~\eqref{eq:MSEGlobalEst} and~\eqref{eq:MSEPooledEst}. 
%This shows that the pooled estimator $\overline{\mu}_N(x_0)$ outperforms the global estimator $\hat{\mu}_N(x_0)$.
The second model considered (right table) illustrates the phenomenon described in Theorem~\ref{superefficient}. We lower bound the supremum risk over $\mathcal{M}_0$ by considering a sequence of alternatives in $\mathcal{M}_0$ (obtained from local perturbations to $\mu(x) = x$ around $x_0 = 0.5$) for which the ratio of the mean squared errors falls dramatically below 1, suggesting that in such a scenario it is better to use the global estimator $\hat{\mu}_N^{-1}(a)$. 
\newline
\newline
It is interesting to note that the super-efficiency phenomenon noted in connection with the pooled estimator in the monotone regression model is also seen with sample-splitting with smoothing based procedures, e.g., kernel based estimation, if the bandwidth used in the divide and conquer method is not appropriately adjusted. We describe the phenomenon in a density estimation setting, since this is the easiest to deal with, in Section \ref{supeffkernel} of the Appendix. Indeed, several authors have criticized such super-efficiency phenomena in nonparametric function estimation; see e.g.,~\cite{BrownEtAl97}, \cite[Section 1.2.4]{Sasha09}, where the authors study super-efficiency in density estimation contexts using kernel methods with a plug-in estimator of the asymptotically ``optimal'' bandwidth. Indeed, it is shown in the second reference that (under the usual twice differentiability assumptions) there exist infinitely many bandwidths that, under any fixed density, produce kernel estimates with asymptotically strictly smaller MSE than the \emph{Epanechnikov oracle} and argued therein that the criterion of assessment of an estimator should therefore be quantified in terms of its maximal risk over an entire class of densities. 

While this is certainly a reasonable perspective --- and indeed, super-efficiency is also encountered with sample-splitting as we have shown above --- we believe that there is also some merit in studying the pointwise behavior of estimators such as in~\eqref{eq:IsoRegLimDist} (as opposed to a uniform measure such as~\eqref{eq:IsoRegUnifBd}). For construction of CIs statisticians usually rely on such pointwise asymptotic results as it is often quite difficult to obtain useful practical procedures that have justification in a uniform sense. Moreover, in the regime of massive datasets, where $N$ is astronomically large, sample-splitting can provide practical gains over the global estimator which might be impossible to compute.

\section{Conclusion}\label{sec: discuss}
We have established rigorous results on sample-splitting in the specific setting of monotone regression and demonstrated both its pros and cons in this problem. The super-efficiency phenomenon demonstrated in this paper is expected to arise more broadly in many of the cube-root $M$-estimation problems as mentioned in the Introduction and developed in \cite{KP90} since the inverse function estimation problem treated in this paper is as an $M$-estimation problem of the type considered in \cite{KP90}. A generic treatment of super-efficiency in these problems should provide an interesting avenue for future research but is outside the scope of this paper. A more general (and harder) question worth considering is a broad characterization of non-standard problems (not necessarily with cube-root convergence rates) where sample-splitting improves the point-wise risk but produces out-of-control uniform risk bounds. 

\section{Proofs of the main results}
\label{proof-section} 
\subsection{Proof of Theorem~\ref{thm:m_n_infty}}\label{proof:thm:m_n_infty}
Since $\{\xi_{n,1}^2\}_{n \ge 1}$ is uniformly integrable and $\xi_{n,1} \stackrel{d}{\to} G$, $\sigma_n^2 := \mbox{Var}(\xi_{n,1})  \rightarrow \sigma^2$  as $n\to\infty$. Set $$Z_n := \sum_{j=1}^{m_n}(\xi_{n,j} - b_n)$$ and let $B_n^2 := \mbox{Var}(Z_n) = m_n \sigma_n^2$.
Now, with $\bar{\xi}_n=m_{n}^{-1}\sum_{j=1}^{m_n}\xi_{n,j} $ we have
\begin{eqnarray*} 
\frac{Z_n}{B_n} & = & \frac{\sum_{j=1}^{m_n}(\xi_{n,j} - b_n)}{\sqrt{m_n}\sigma_n} = \frac{\sqrt{m_n}(\bar{\xi}_n - b_n)}{\sigma_n} \\
& = & 
\frac{\sqrt{m_n}r_n (\bar {\theta}_{m_n} - \theta_0)}{\sigma_n} - \frac{\sqrt{m_n} \,b_n}{\sigma_n} \equiv I_n - II_n.
\end{eqnarray*} 
We show that $Z_n/B_n \stackrel{d}{\rightarrow} N(0,1)$. To this end, we just need to verify the Lindeberg condition: for every $\epsilon > 0$, 
\[ \frac{1}{\sigma_n^2}\,\E [(\xi_{n,1} - b_n)^2\mathbf{1}\{|\xi_{n,1} - b_n| > \epsilon \sqrt{m_n}\,\sigma_n\}] \rightarrow 0.\] 
Since $\sigma_n^2$ converges to $\sigma^2 > 0$ and $m_n \rightarrow \infty$, the above condition is implied by the uniform integrability of $\{(\xi_{n,1} - b_n)^2\}_{n \ge 1}$ which is guaranteed by the uniform integrability of $\{\xi_{n,1}^2\}$ (since the sequence $b_n$ goes to 0 and is therefore bounded). Hence, $Z_n/B_n \stackrel{d}{\to} N(0,1)$.

Now assume that $m_n$ is as in $(i)$. Then, $II_n \rightarrow 0$, which implies that $$I_n = \frac{\sqrt{m_n}r_n (\bar {\theta}_{m_n} - \theta_0)}{\sigma_n} \stackrel{d}{\to} N(0,1),$$  and therefore the conclusion of $(i)$. Next, if $m_n$ is as in $(ii)$, $II_n \rightarrow \tau/\sigma$, and the conclusion of $(ii)$ follows. \qed

\subsection{Proof of Theorem~\ref{theo: currentstatus}}\label{proof:theo: currentstatus}
Let $\mu=-F_T$ and for all $i=1,\dots,n$, let $Y_{i}=-\BBone_{T_{i}\leq X_{i}}$ and
$\eps_{i}=Y_{i}-\mu(X_{i})\in [-1,1].$
Moreover, define $v^2(x) :=\E(\eps_i^2|X_i=x)$ for all $x\in[0,1]$. We then have
$$v^2(x)=\mathrm{Var}(\BBone_{T\leq x})=F_T(x)(1-F_T(x)).$$
With  $\hat \mu_{n}$  defined as in Section \ref{sec: LSE}, we have $\hat F_{Tn}=-\hat \mu_{n}$; see~\cite{groeneboom1992information}. This means that $F_{Tn}^{-1}(a)=\hat \mu_{n}^{-1}(-a)$. Moreover, $F_{T}^{-1}(a)=\mu^{-1}(-a)$. 
Now, under the assumptions of  Theorem \ref{theo: currentstatus}, {\monorandom} and {\density} hold true. The assumption {\sigmamin} holds since $T$ has a density function that is bounded away from zero on $[0,1]$. Moreover, {\momentq} holds true since the $\eps_{i}$'s are bounded. Hence, Theorem \ref{theo: UnifIntRD}  applies to $\hat{\mu}_n$ which translate to conclusions in 1 of this theorem. The  conclusions in 2 (on the orders of the bias of $\hat{F}_{Tn}$ and $\hat{F}_{Tn}^{-1}$) follow by a direct application of Theorems \ref{theorem: BiasInverse} 
and \ref {theo: biasDirect} and the conclusions in 3 follow exactly in the same fashion as for the simple signal plus noise regression model considered above. $\qed$

\subsection{Proof of Theorem~\ref{superefficient}}
\label{supeff-proof} 
For the proof of this theorem we assume that $\mu_0$ is \emph{nondecreasing} --- this is convenient as we borrow several results from other papers stated in the context when $\mu_0$ is nondecreasing. The neighborhood $\mathcal{M}_0$ in the statement of the theorem needs to be similarly modified. (Of course, appropriate changes will lead to the proof of the case when $\mu_0$ is nonincreasing.) 

By conclusion 2 of Corollary~\ref{maximal-risk} (adapted to nondecreasing functions), with $p = 2$ and noting that $\mathcal{M}_0$ is a subset of an appropriate $\mathcal{F}_1$ we conclude that $E < \infty$. Letting 
\[V_1 := \limsup_{N \to \infty}\,\sup_{\mu \in \mathcal{M}_0}\mathrm{Var}_{\mu}[N^{1/3}(\hat{\theta}_N - \mu^{-1}(a))],\]
and 
\[ V_2 := \limsup_{N \to \infty}\,\sup_{\mu \in \mathcal{M}_0}\,N^{2/3}[\E_{\mu}\,\hat{\theta}_N - \mu^{-1}(a)]^2 \,,\]
we have $V_1 \vee V_2 \leq E \leq V_1 + V_2 < \infty.$ Recall that as $\overline{\theta}_{m}$ is the average of the $m$ i.i.d.~random variables $\hat{\mu}_{n,j}^{-1}(a)$, $j = 1,\ldots, m$, $\E_{\mu}(\overline{\theta}_m) = \E_{\mu}(\hat{\mu}_{n,1}^{-1}(a))$. Now, consider 
\begin{eqnarray*}
V_{2,m} &:=&  
\liminf_{N}\,\sup_{\mu \in \mathcal{M}_0}\,N^{2/3}[\E_{\mu}\,\overline{\theta}_m - \mu^{-1}(a)]^2 \\
& = & m^{2/3} \liminf_{n \to \infty}\,\sup_{\mu \in \mathcal{M}_0}\,n^{2/3}[\E_{\mu}\,\hat{\mu}_{n,1}^{-1}(a) - \mu^{-1}(a)]^2 \;\; =: \;\; m^{2/3}\,\widetilde{V}_2.
\end{eqnarray*} 
Note that, $E_m \geq V_{2,m} = m^{2/3}\,\widetilde{V}_2$. We will show below that $\widetilde{V}_2 > 0$; thus $c_0$ in the statement of the theorem can be chosen to be $\widetilde{V}_2$. To this end, consider the monotone regression model under a sequence of  \emph{local alternatives} $\mu_n$ which eventually lie in $\mathcal{M}_0$.  Let $Y = \mu_n(X) + \epsilon$ where everything is as before but $\mu_0$ changes to $\mu_n$ which is defined as  $$\mu_n(x) = \mu_0(x) + n^{-1/3}B\,(n^{1/3}(x- \theta_0))$$ and $B$ is \emph{a non-zero function} continuously differentiable on $\R$, vanishing outside $(-1,1)$, such that $\mu_n$ is monotone for each $n$ and  lies eventually in the class $\mathcal{M}_0$\footnote{There is nothing special about $(-1,1)$ as far as constructing the $B$ is concerned. Any $(-c,c)$, for $c > 0$ can be made to work.}. Note that $\mu_n$ and $\mu_0$ can differ on 
$(\theta_0 - n^{-1/3}, \theta_0 + n^{-1/3})$ only, and that $\mu_n'(x) = \mu_0'(x) + B'(n^{1/3}(x - \theta_0))$ for $x \in [\theta_0 - n^{-1/3}, \theta_0 + n^{-1/3}]$ and $\mu_n'(x) = \mu_0'(x)$ otherwise. It is clear that this can be arranged for infinitely many $B$'s.

The above sequence of local alternatives was considered in \cite{banerjee2005likelihood} in a more general setting, namely that of monotone response models, where (in a somewhat unfortunate collision of notation) $X$ denotes \emph{response} and $Z$ the covariate. We invoke the results of that paper \emph{using the $(Y,X)$ notation of this paper and ask the reader to bear this in mind}.  Using our current notation for the problem in \cite{banerjee2005likelihood}, $X$ follows density $p_X(x) = \BBone_{(0,1)}(x)$ and $Y \mid X = x \sim p(y,\psi(x))$, $\psi$ being a monotone function and $p(y,\theta)$ a regular parametric model.  The monotone regression model with homoscedastic normal errors under current consideration is a special case of this setting with $p(y,\theta)$ being the $N(\theta, v^2)$ density, the $\psi_n$'s in that paper defining the local alternatives are the monotone functions $\mu_n$, $\psi_0 = \mu_0$, $c=1$ and $A_n(x) = B(n^{1/3}(x-x_0))$ for all $n$. 
Invoking Theorems 1 and 2 of \cite{banerjee2005likelihood} with the appropriate changes, we conclude that under $\mu_n$,
\[ X_n(h) := n^{1/3}(\hat{\mu}_n(\theta_0 + h n^{-1/3}) - \mu_0(\theta_0)) \stackrel{d}{\rightarrow}g_{c,d,\mathcal{D}}(h),\]
where $c = v, d = \mu_0'(x_0)/2$, $\mathcal{D}$ is a shift function 
%defined on Page 514 of \cite{banerjee2005likelihood} 
given by\footnote{There is a typo in the drift term as stated on page 514 of \cite{banerjee2005likelihood}: there should be a negative sign before the integral that defines $\mathcal{D}(h)$ for $h < 0$ on page 514.}:
\[ \mathcal{D}(t) = \left(\int_{0}^{t \wedge 1}\,B(u) du \right)\,\BBone_{(0,\infty)}(t) -  \left(\int_{t \vee -1}^0\,B(u) du \right)\,\BBone_{(-\infty,0)}(t),\]
and $g_{c,d,\mathcal{D}}$ is the right-derivative process of the greatest convex minorant (GCM) of $X_{c,d,\mathcal{D}}(t) : = c W(t) + d t^2 + \mathcal{D}(t)$ with $W$ being a two-sided Brownian motion. Now, by essentially the same calculation as on Page 422 of \cite{BW05},  
\[ P(n^{1/3}[\hat{\mu}_n^{-1}(a + \lambda n^{-1/3}) - \mu_0^{-1}(a)] \leq x) = P(n^{1/3}(\hat{\mu}_n(\theta_0 + x n^{-1/3}) - \mu_0(\theta_0)) \geq \lambda) 
\rightarrow P(g_{c,d,\mathcal{D}}(x) \geq \lambda) \,.\]
Setting $\lambda = 0$, we get:
\[ P(n^{1/3}[\hat{\mu}_n^{-1}(a) - \mu_0^{-1}(a)] \leq x) = P(n^{1/3}(\hat{\mu}_n(\theta_0 + x n^{-1/3}) - \mu_0(\theta_0)) \geq 0) 
\rightarrow P(g_{c,d,\mathcal{D}}(x) \geq 0) \,.\]
Next, by the switching relationship\footnote{For the details, see Section \ref{super-efficient-details} of the Appendix.},  
\[ P(g_{c,d,\mathcal{D}}(x) \geq 0) = P(\arg\min_h X_{c,d,\mathcal{D}}(h)  \leq x) \,,\]
and it follows that: 
\[ n^{1/3}(\hat{\mu}_n^{-1}(a) - \mu_0^{-1}(a)) \stackrel{d}{\rightarrow}  \arg\min_h\,X_{c,d,\mathcal{D}}(h) \,.\] 
Choosing $B$ such that $B(0) = 0$, we note that $\mu_n^{-1}(a) = \mu_0^{-1}(a) = \theta_0$, and therefore, under the sequence of local alternatives $\mu_n$,
\begin{equation}
\label{dist-conv-localt}
n^{1/3}(\hat{\mu}_n^{-1}(a) - \mu_n^{-1}(a)) \stackrel{d}{\rightarrow} \arg\min_h X_{c,d,\mathcal{D}}(h).
\end{equation}
Since the $\mu_n$'s eventually fall within the class $\mathcal{M}_0$, by conclusion 2 of (the version of) Corollary~\ref{maximal-risk} (for nondecreasing functions), we conclude that: 
\[ \limsup_{n \to \infty} n^{2/3}\, \E_{\mu_n}\left(|\hat \mu_{n}^{-1}(a) -\mu_n^{-1}(a)|^2\right)\leq K_{2}.\]
Thus the sequence $\{n^{1/3}(\hat{\mu}_n^{-1}(a) - \mu_n^{-1}(a))\}_{n \ge 1}$ is uniformly integrable under the sequence (of probability distributions corresponding to) $\{\mu_n\}_{n \ge 1}$ and in conjunction with (\ref{dist-conv-localt}) it follows that
\[ \lim_{n \to \infty}\,n^{1/3}[\E_{\mu_n}(\hat{\mu}_n^{-1}(a) - \mu_n^{-1}(a))] = \E( \arg\min_h X_{c,d,\mathcal{D}}(h)) \,.\]
[{\bf Claim $\mathbf{C}$}] (proved in Section \ref{super-efficient-details} of the Appendix): For any non-negative function $B$ that satisfies the conditions imposed above, and is additionally symmetric about 0, 
$$ \E(\arg\min_h X_{c,d,\mathcal{D}}(h)) \ne 0\,.$$ 
It follows that for any such $B$, 
\[ [\E(\arg\min_h X_{c,d,\mathcal{D}}(h))]^2 \leq \widetilde{V}_2 ,\] 
and hence $\widetilde{V}_2 > 0$. This delivers the assertions of the theorem for fixed $m$.
\newline
\newline
When $m = m_n \to \infty$, note that
\begin{eqnarray*}
\liminf_{N \to \infty}\,\sup_{\mu \in \mathcal{M}_0}\E_{\mu}\left[ N^{2/3}(\overline{\theta}_{m_n} - \mu^{-1}(a))^2 \right] &\geq & 
\liminf_{N\to \infty}\,\sup_{\mu \in \mathcal{M}_0}\,N^{2/3}[\E_{\mu}\,\overline{\theta}_{m_n} - \mu^{-1}(a)]^2 \\
& \geq & \liminf_{n \rightarrow \infty}\,m_n^{2/3}\sup_{\mu \in \mathcal{M}_0}\,n^{2/3}[\E_{\mu}\,\overline{\theta}_{m_n} - \mu^{-1}(a)]^2 \\
& = & \liminf_{n \rightarrow \infty}\,m_n^{2/3}\,\sup_{\mu \in \mathcal{M}_0}\,n^{2/3}[\E_{\mu}\,\hat{\mu}_{n,1}^{-1}(a) - \mu^{-1}(a)]^2 \,.
\end{eqnarray*}
By our derivations above, 
\[ \sup_{\mu \in \mathcal{M}_0}\,n^{2/3}[\E_{\mu}\,\hat{\mu}_{n,1}^{-1}(a) - \mu^{-1}(a)]^2 \ge \frac{1}{2} 
[\E(\arg\min_h X_{c,d,\mathcal{D}}(h))]^2 > 0 \]
for all sufficiently large $n$, and it follows that the liminf of the maximal normalized risk of $\overline{\mu}_N$ is infinite. 
\qed

\subsection{Some Selected Proofs for Section \ref{sec: intdirect}}\label{sec: randomU}
%%%%%%%%%%%%%%%%%%%%%%%%%%%%%%
We start with a precise exposition of the characterization of the estimators as this is critical to the subsequent analysis. From \eqref{brunk0} we have 
\begin{equation}\label{brunk}
\hat \mu_n(X_{(i)})=\hat\lambda_n(i/n)=\hat\lambda_n \circ F_n(X_{(i)}),\qquad i=1,\dots,n,
\end{equation}
where $F_n$ is the empirical distribution function of $X_1,\dots,X_n$. A convenient way of studying $\hat \mu_{n}$ is to first study $\hat\lambda_{n}$ and then go back to $\hat \mu_{n}$ thanks to \eqref{brunk}. Note that   $\hat\lambda_{n}(i/n)=\hat \mu_n\circ F_{n}^{-1}(i/n)$ for all $i\in\{1,\dots,n\}$, where $ F_{n}^{-1}(a)$ is the smallest $t\in[0,1]$ that satisfies $F_n(t)\geq a$, for all $a\in\R$. Both functions $\hat\lambda_{n}$ and $\hat \mu_n \circ F_{n}^{-1}$ are piecewise constant, so $\hat\lambda_n= \hat \mu_n \circ F_{n}^{-1}$ on $[0,1]$ and
 $\hat\lambda_{n}$ can be viewed as an estimator of the function $\lambda$ defined on $[0,1]$ by
\begin{equation}\label{eq: lambdavsm}
\lambda=\mu\circ F^{-1}.
\end{equation}
Hereafter, we denote by $\mu^{-1}$ and $g$ the respective generalized inverses of $\mu$ and $\lambda$. This means that $\mu^{-1}$ and $g$ extend the usual inverses to the whole real line in such a way that they remain constant on $(-\infty,0]$ and on $[1,\infty)$. Letting $\hat\mu_{n}^{-1}$ and $\hat U_n$ be the respective generalized inverses of $\hat\mu_{n}$ and $\hat \lambda_n$, it follows from \eqref{brunk} that
\begin{equation}\label{invmm}
\hat \mu_n^{-1}=F_n^{-1}\circ \hat U_n,
\end{equation}
and it can be shown that 
\begin{equation}\label{UnRandom}
\hat U_n(a)=\argmax_{u\in[0,1]}\{\Lambda_n(u)-au\}, \qquad \mbox{ for all } a\in\R
\end{equation}
where argmax denotes the greatest location of maximum (which  is achieved on the set $\{i/n,\ i=0,\dots,n\}$ since $\Lambda_{n}$ is 
piecewise-linear).  Part of the proofs below consist in first establish a result for $\hat U_n$ using the above characterization, and then go from $\hat U_{n}$ to $\hat \mu_{n}^{-1}$ using \eqref{invmm}. To this end, we will use a precise bound for the uniform distance between $F^{-1}$ and  $F_{n}^{-1}$, as well as a strong approximation of the empirical quantile function, see Section~\ref{sec: preliminaries} in the Appendix.

In what follows, we will repeatedly use the fact that because $g'=1/\lambda'\circ g$ on $(\lambda(1),\lambda(0))$ where $\lambda'=\mu'\circ F^{-1}/f\circ F^{-1}$ is bounded away from zero (see {\monorandom} and {\densityUI}), we have
\begin{equation}\label{eq : gLips}
|g(u)-g(v)|\leq \frac{1}{\inf_{t\in[0,1]}|\lambda'(t)|}|u-v|
\end{equation}
for all real numbers $u$ and $v$. Furthermore, we recall that from Fubini's theorem, it follows that for all $r\geq 1$ and all random variables $Z$, 
\begin{equation}\label{fubini}
\E\vert Z\vert ^r=\int _{0}^\infty \P(\vert Z\vert^r >x)dx= \int _{0}^\infty \P(\vert Z\vert>t)rt^{r-1}dt.
\end{equation}
We denote by $\P^X$ the conditional probability given $(X_{1},\dots,X_{n})$ and by $\E^X$ the corresponding conditional expectation.

\subsubsection{Preliminaries: Exponential bounds for tail probabilities}\label{sec: expobounds} 
In this subsection, we provide exponential bounds, which are proved in Appendix~\ref{Appendix}, for the tail probabilities of $\hat\mu_{n}^{-1}$ and $\hat U_{n}$. We begin with a generalization to our setting of Theorem 11.3 in \cite{groeneboom2014nonparametric}. Also, the lemma is a stronger version of inequality (11) in \cite{Durot08} where an assumption (A5) was postulated instead of the stronger assumption {\momentexpo}. The lemma will be used in the proof of Theorem \ref{theo: biasDirect}.

\begin{lemma}\label{lem: expoboundinvm}
Assume {\momentexpo}, $X$ has a density function $f$, $\mu$ is nonincreasing and there exist positive numbers $A_{1},\dots,A_{4}$ such that $A_{1}<|\mu'(t)|<A_{2}$ and $A_{3}<f(t)<A_{4}$
for all $t\in[0,1]$. Then, there exist positive numbers $K_{1}$ and $K_{2}$ that depend only on $A_{1},\dots,A_{4},K,\alpha$, where $K$ and $\alpha$ are taken from  {\momentexpo}, such that for all $n$, $a\in\R$ and $x>0$, we have
\begin{equation}\label{eq: expoboundinvm}
\P\left(|\hat \mu_{n}^{-1}(a)-\mu^{-1}(a)|>x\right)\leq K_1\exp(-K_{2}nx^3).
\end{equation}
\end{lemma}

To prove Lemma \ref{lem: expoboundinvm}, we first prove a similar bound for $\hat U_{n}$. The exponential bound for $\hat U_{n}$ is given in the following lemma. It will  be used also in the proof of Theorem \ref{theo: UnifIntRD}.

\begin{lemma}\label{lem: expoboundU}
Assume {\momentexpo}, $X$ has a density function $f$, $\mu$ is nonincreasing and there exist positive numbers $A_{1},\dots,A_{4}$ such that 
$A_{1}<|\mu'(t)|<A_{2}$ and $A_{3}<f(t)<A_{4}$
for all $t\in[0,1]$. Then, there exist positive numbers $K_{1}$ and $K_{2}$ that depend only on $A_{1},\dots,A_{4},K,\alpha$, where $K$ and $\alpha$ are taken from  {\momentexpo}, such that for all $n$, $a\in\R$ and $x>0$, we have
\begin{equation}\label{eq: expoboundU}
\P\left(|\hat U_{n}(a)-g(a)|>x\right)\leq K_1\exp(-K_{2}nx^3).
\end{equation}
\end{lemma}

To prove Theorem \ref{theo: UnifIntRD}, we also need a sharper  inequality for the cases when $a\not\in[\lambda(1),\lambda(0)]$.

\begin{lemma}\label{lem: expoboundU+-}
Assume {\momentexpo}, $X$ has a density function $f$, and $\mu$ is nonincreasing. Then, there exist positive numbers $K_{1}$ and $K_{2}$ that depend only on $K$ and $\alpha$, which are taken from {\momentexpo}, such that
\begin{equation}\label{eq: expoboundU+}
\P^X\left(\hat U_{n}(a)\geq x\right)\leq K_1\exp(-K_{2}(a-\lambda(0))^2nx)
\end{equation}for all $n$, $a>\lambda(0)$ and $x\geq n^{-1}$, and
\begin{equation}\label{eq: expoboundU-}
\P^X\left(1-\hat U_{n}(a)\geq x\right)\leq K_1\exp(-K_{2}(a-\lambda(1))^2nx)
\end{equation}for all $n$, $a<\lambda(1)$ and $x\geq n^{-1}$.
\end{lemma}

\subsubsection{Proof of Theorem \ref{theo: UnifIntRD}}
Integrating the inequality in Lemma \ref{lem: expoboundinvm} according to \eqref{fubini} proves the first assertion.  To prove the second one, we first prove a similar result for $\hat\lambda_{n}$.
\begin{lemma}\label{lem: UnifIntlambda}
Assume {\momentexpo}, $X$ has a density function $f$, $\mu$ is nonincreasing, and there exist positive numbers $A_{1},\dots,A_{4}$ such that $A_{1}<|\mu'(t)|<A_{2}$ and $A_{3}<f(t)<A_{4}$
for all $t\in[0,1]$. Then, for all $p> 0$ and $A>0$,  there exist positive numbers $K_{1}$ and $K_{2}$ that depend only on $A_{1},\dots,A_{4},K,\alpha,p,A$, where $K$ and $\alpha$ are taken from  {\momentexpo}, such that
$$\E\left(n^{1/3}|\hat \lambda_{n}(t)-\lambda(t)|\right)^p\leq K_{p,A}$$
for all $n$ and $t\in[n^{-1/3}A,1-n^{-1/3}A]$.
\end{lemma}

\paragraph{\bf Proof} As is customary, we denote $y_{+}=\max(y,0)$ and $y_{-}=-\min(y,0)$ for all $y\in\mathbb R$. To go from $\hat U_{n}$ to $\hat\lambda_{n}$ we will make use of the following switch relation, that holds for all $t\in(0,1]$ and $a\in\R$:
\begin{equation}\label{eq: switchU}
\hat \lambda_{n}(t)\geq a \Longleftrightarrow t\leq \hat U_{n}(a).
\end{equation}
With $a_{x}=\lambda(t)+x$, it then follows from \eqref{fubini} and the switch relation \eqref{eq: switchU} that 
\begin{eqnarray}\label{eq: EpintU}\notag
\E\left(( \hat\lambda_{n}(t)- \lambda(t))_{+}\right)^p&=&\int_{0}^\infty\P\left(\hat  \lambda_{n}(t)- \lambda(t)\geq x\right)px^{p-1}dx\\
&=&\int_{0}^\infty\P\left(\hat  U_{n}(a_{x})\geq  t\right)px^{p-1}dx\\
&=&I_{1}+I_{2}\end{eqnarray}
where
$$I_{1}=\int_{0}^{\lambda(0)-\lambda(t)}\P\left(\hat  U_{n}(a_{x})\geq  t\right)px^{p-1}dx\quad\mbox{ 
and }\quad
I_{2}=\int_{\lambda(0)-\lambda(t)}^\infty\P\left(\hat  U_{n}(a_{x})\geq  t\right)px^{p-1}dx.$$

Consider $I_{1}$. Since $\lambda=\mu\circ F^{-1}$, it follows from the Taylor expansion that with $c=A_{3}/A_{2}$, we have
$t-\lambda^{-1}(a_{x})> cx$
for all $x\in(0,\lambda(0)-\lambda(t))$. Therefore, \eqref{eq: expoboundU}
 implies that
\begin{equation}\notag
\P\left(\hat  U_{n}(a_{x})\geq  t\right)\leq \P\left(\hat U_{n}(a_x)-\lambda^{-1}(a_{x})>cx\right)\leq K_1\exp(-K_{2}c^3nx^3)
\end{equation}
for all $x\in(0,\lambda(0)-\lambda(t))$. Hence,
\begin{eqnarray*}
I_{1}\; \leq \; K_{1} \int_{0}^{\lambda(0)-\lambda(t)}\exp(-K_{2}c^3nx^3)px^{p-1}dx \; \leq \; K_{1} n^{-p/3} \int_{0}^{\infty}\exp(-K_{2}c^3y^3)py^{p-1}dy,
\end{eqnarray*}
using the change of variable $y=n^{1/3}x$. The integral on the right hand side depends only on $c$ and $p$, and is finite for all $p> 0$. Hence, with $C_{p}/K_{1}$ greater than this integral we obtain
\begin{equation}\label{eq: majI1U}
I_{1}\leq C_{p}n^{-p/3}.
\end{equation}

Now consider $I_{2}$. We have $a_{x}>\lambda(0)$ for all $x>\lambda(0)-\lambda(t)$ so it follows from \eqref{eq: expoboundU+} together with \eqref{eq: expoboundU} (where $g(a_{x})=0$) that
\begin{eqnarray*}
I_{2}&\leq& K_{1} \int_{\lambda(0)-\lambda(t)}^{2(\lambda(0)-\lambda(t))} \exp(-K_{2}nt^3)px^{p-1}dx+K_{1} \int_{2(\lambda(0)-\lambda(t))}^\infty \exp(-K_{2}(a_{x}-\lambda(0))^2nt)px^{p-1}dx\\
&\leq& K_{1} \exp(-K_{2}nt^3)2^p(\lambda(0)-\lambda(t))^p+K_{1} \int_{2(\lambda(0)-\lambda(t))}^\infty \exp(-K_{2}x^2nt/4)px^{p-1}dx,
\end{eqnarray*}
since $a_{x}-\lambda(0)\geq x/2$ for all $x\geq 2(\lambda(0)-\lambda(t))$. Since the sup-norm of $\lambda'$ is smaller than or equal to $ A_{2}/A_{3}$ we then have 
\begin{eqnarray*}
I_{2}&\leq& K_{1}2^p (A_{2}/A_{3})^p \exp(-K_{2}nt^3)t^p+K_{1} (nt)^{-p/2} \int_{0}^\infty \exp(-K_{2}y^2/4)py^{p-1}dy
\end{eqnarray*}
using the change of variable $y=x\sqrt{nt}$. The function $t\mapsto \exp(-K_{2}nt^3)t^p$ achieves its maximum on $[0,\infty)$ at the point $(3K_{2}n/p)^{-1/3}$. This means that for all $t\geq 0$ we have
$$\exp(-K_{2}nt^3)t^p\leq \exp(-p/3)\left(\frac{3K_{2}n}p\right)^{-p/3}.$$
On the other hand, we have $(nt)^{-p/2}\leq A^{-p/2}n^{-p/3}$ for all $t\geq n^{-1/3}A$, where $A>0$ is fixed. Combining this with the two preceding displays, we arrive at
\begin{eqnarray*}
I_{2}&\leq& K_{1}2^p (A_{2}/A_{3})^p \exp(-p/3)\left(\frac{3K_{2}n}p\right)^{-p/3}+K_{1} A^{-p/2}n^{-p/3}\int_{0}^\infty \exp(-K_{2}y^2/4)py^{p-1}dy
\end{eqnarray*}
 for all $t\geq n^{-1/3}A$, where the integral on the right hand side is finite. This means that  there exists $K_{p,A}>0$ such that
 $I_{2}\leq K_{p,A}n^{-p/3}/2$
 for all $t\geq n^{-1/3}A$.
Combining this with \eqref{eq: EpintU} and \eqref{eq: majI1U} and possibly enlarging  $K_{p,A}>0$, we obtain
$$\E\left((\hat  \lambda_{n}(t)- \lambda(t))_{+}\right)^p\leq K_{p,A}n^{-p/3}$$
 for all $t\geq n^{-1/3}A$. It can be proved with similar arguments  that the above inequality remains valid with $(\cdot)_+$ replaced by $(\cdot)_-$, and Lemma \ref{lem: UnifIntlambda} follows. \cqfd

It is known that Grenander type estimators are inconsistent at the boundaries. However,  the following lemma shows that  such estimators remain bounded in the $\L_{p}$-sense. The lemma, which is proved in Appendix~\ref{Appendix}, will be useful to go from Lemma \ref{lem: UnifIntlambda} to Theorem \ref{theo: UnifIntRD}.

\begin{lemma}\label{lem: boundaries}
%Assume {\monorandom}, {\densityUI} and {\momentexpo}. For all $p>0$ we then have $\E|\hat\lambda_{n}(0)|^p=O(1)$ and $\E|\hat\lambda_{n}(1)|^p=O(1)$. 
Assume  {\momentexpo} and $\mu$ is nonincreasing with
$|\mu(t)|\leq A_{5}$
for some $A_{5}>0$ and all $t\in[0,1]$. Then, for all $p>0$, there exists  $K_{1}>0$  that depends only on $p,A_{5,}K$ and $\alpha$, where $K$ and $\alpha$ are taken from {\momentexpo}, such that $\E|\hat\lambda_{n}(0)|^p\leq K_{1}$ and $\E|\hat\lambda_{n}(1)|^p\leq K_{1}$, $\forall n$.
\end{lemma}

We are now in a position to prove the second assertion in Theorem \ref{theo: UnifIntRD}.
Since  $\hat \mu_{n}$ is constant on all intervals $(X_{(i)},X_{(i+1)}]$ for $i\in\{1,\dots,n-1\}$ and also on the interval $[0,X_{(1)}]$, and $F_{n}$ is constant on all intervals $[X_{(i)},X_{(i+1)})$ for $i\in\{1,\dots,n-1\}$ and also on the interval $[0,X_{(1)})$, it follows from \eqref{brunk} that for all $t\not\in\{X_{(1)},\dots,X_{(n)}\}$ we have
$\hat \mu_{n}(t)=\hat\lambda_{n}(F_{n}(t)+n^{-1}).$
But the $X$ has a continuous distribution so for a fixed $t$, we indeed have $t\not\in\{X_{(1)},\dots,X_{(n)}\}$ with probability one.
Hence, for all $p\geq 1$ we have
\begin{eqnarray*}
\E\left((\hat \mu_{n}(t)-\mu(t))_{+}\right)^p&=&\E\left(\left(\hat\lambda_{n}(F_{n}(t)+n^{-1})-\lambda(F(t))\right)_{+}\right)^p.
\end{eqnarray*}
Using monotonicity of $\hat \lambda_{n}$, this means that
\begin{eqnarray}\label{eq: UIm1}\notag
\E\left((\hat \mu_{n}(t)-\mu(t))_{+}\right)^p&\leq &\E\left(\left(\hat\lambda_{n}(F(t)-n^{-1/2}\log n)-\lambda(F(t))\right)_{+}\right)^p\\
& &\qquad+\E\left(\left(\hat\lambda_{n}(0)-\lambda(1)\right)_{+}^p\BBone_{F_{n}(t)+n^{-1}\leq F(t)-n^{-1/2}\log n}\right).
\end{eqnarray}
It follows from H\"older's inequality that
\begin{equation}\notag\begin{split}
&\E\left(\left(\hat\lambda_{n}(0)-\lambda(1)\right)_{+}^p\BBone_{F_{n}(t)+n^{-1}\leq F(t)-n^{-1/2}\log n}\right)\\
&\qquad\leq
\E^{1/2}\left(\left(\hat\lambda_{n}(0)-\lambda(1)\right)^{2p}\right)\P^{1/2}\left(F_{n}(t)+n^{-1}\leq F(t)-n^{-1/2}\log n\right)\\
&\qquad\leq
\E^{1/2}\left(\left(\hat\lambda_{n}(0)-\lambda(1)\right)^{2p}\right)\P^{1/2}\left(\sup_{t\in[0,1]}|F_{n}(t)-F(t)|> n^{-1/2}\log n\right).
\end{split}
\end{equation}
Combining this with Lemma \ref{lem: boundaries} together with  Corollary 1 in \cite{massart1990tight} yields
\begin{equation}\notag\begin{split}
\E\left(\left(\hat\lambda_{n}(0)-\lambda(1)\right)_{+}^p\BBone_{F_{n}(t)+n^{-1}\leq F(t)-n^{-1/2}\log n}\right)\leq
O(1)\left( 2\exp(-2(\log n)^2)\right)^{1/2}
\end{split}
\end{equation}
where the big-$O$ term is uniform for all functions $\mu$ satisfying the assumptions of the lemma.
This means that there exists $C_{p}>0$ such that
\begin{equation}\label{eq: UIm2}
\E\left(\left(\hat\lambda_{n}(0)-\lambda(1)\right)_{+}^p\BBone_{F_{n}(t)+n^{-1}\leq F(t)-n^{-1/2}\log n}\right)\leq C_{p}n^{-p/3}
\end{equation}
for all $t\in[0,1]$. Now, consider the first term on the right hand side of \eqref{eq: UIm1}. 
It follows from the convexity of the function $x\mapsto x^p$ that $(x+y)^p\leq 2^{p-1}(x^p+y^p)$ for all positive numbers $x$ and $y$. Therefore, with $t\geq n^{-1/3}$ and $x_{n}=F(t)-n^{-1/2}\log n$  we have
\begin{eqnarray*}
\E\left(\left(\hat\lambda_{n}(x_{n})-\lambda(F(t))\right)_{+}\right)^p&\leq&2^{p-1}\E\left(|\hat\lambda_{n}(x_{n})-\lambda(x_{n}))|^p\right)+2^{p-1}|\lambda(x_{n})-\lambda(F(t)))|^p\\
&\leq&2^{p-1}\E\left(|\hat\lambda_{n}(x_{n})-\lambda(x_{n}))|^p\right)+2^{p-1} (A_{2}/A_{3})^p n^{-p/2}(\log n)^p
\end{eqnarray*}
since  the sup-norm of $\lambda'$ is less than or equal to $A_{2}/A_{3}$. Let 
$A \leq A_{3}/2$. 
For $n$ sufficiently large, we have $x_{n}\in[n^{-1/3}A,1-n^{-1/3}A]$ for all $t\in[n^{-1/3},1-n^{-1/3}]$. This means that the previous display combined with Lemma \ref{lem: UnifIntlambda} ensures that there exists $C_{p}>0$ such that
\begin{eqnarray*}
\E\left(\left(\hat\lambda_{n}(x_{n})-\lambda(F(t))\right)_{+}\right)^p\leq C_{p}n^{-p/3}
\end{eqnarray*}
for all $t\in[n^{-1/3},1-n^{-1/3}]$ and $n$ sufficiently large.  Together with \eqref{eq: UIm2} and \eqref{eq: UIm1}, this yields
\begin{eqnarray}\notag
\E\left((\hat \mu_{n}(t)-\mu(t))_{+}\right)^p \leq 2C_{p}n^{-p/3}
\end{eqnarray}
for all $t\in[n^{-1/3},1-n^{-1/3}]$ and $n$ sufficiently large. Possibly enlarging $C_{p}$, the previous inequality remains true for all $n$. To see this, suppose that the above display holds for all $n \geq n_{min}$. Now, 
\[ \E((\hat \mu_n(t)-\mu(t))_+)^p  \leq  2^{p-1}E( |\hat\mu_n(0)|^p \vee |\hat\mu_n(1)|^p) + 2^{p-1} |\mu(0)|^p\vee|\mu(1)|^p \,.\] 
by monotonicity of both $\mu$ and $\hat\mu_n$, and using convexity of the function $x\mapsto x^p$. Hence, for $n < n_{min}$, 
$$n^{p/3}\E((\hat \mu_n(t)-\mu(t))_+)^p  \leq  (2^{p}K_1 + 2^{p} A_5) n_{min}^{p/3}\,,$$
where $K_1$ and $A_5$ are taken from Lemma \ref{lem: boundaries}. It can be proved likewise that there exists $C_{p}>0$ such that
$\E\left((\hat \mu_{n}(t)-\mu(t))_{-}\right)^p \leq 2C_{p}n^{-p/3}$ for all $t\in[n^{-1/3},1-n^{-1/3}]$ and all $n$. This completes the proof of Theorem \ref{theo: UnifIntRD}. \cqfd

\subsubsection{Proof of Theorem \ref{theorem: BiasInverse}}
Theorem \ref{theorem: BiasInverse}  immediately follows from Lemma \ref{lemma: ConnectingBias2} combined to Theorem \ref{theo: randomU}  below by noticing that $\mu(1)=\lambda(1)$ and $\mu(0)=\lambda(0)$.
Theorem \ref{theo: randomU} provides  a precise bound for the bias of $\hat U_{n}$  whereas Lemma \ref{lemma: ConnectingBias2} makes the connection between the biases of $\hat\mu_{n}^{-1}$ and $\hat U_{n}$. The lemma is proved in Section \ref{sec: proof ConnectingBias2} in the Appendix, using that $\mu^{-1}=F^{-1}\circ g$  and $\hat \mu_n^{-1}=F_n^{-1}\circ \hat U_n,$ where $F_n^{-1}$ estimates $F^{-1}$.

\begin{lemma}\label{lemma: ConnectingBias2}
Assume {\monorandom}, {\density}  and {\momentq}. Denote by $\mu^{-1}$ and $g$ the respective generalized inverses  of $\mu$ and $\lambda$.
We then have
$$\E\left(\hat \mu_{n}^{-1}(a)-\mu^{-1}(a)\right)=\frac{1}{f\circ F^{-1}(g(a))}\E\left(\hat U_{n}(a)-g(a)\right)+o(n^{-1/2})$$
where the small-$o$ term is uniform in $a\in\R$. 
\end{lemma}

\begin{theorem}\label{theo: randomU}
Assume {\monorandom}, {\density}, {\sigmamin}, {\momentq},  $v^2$ has a bounded second derivative on $[0,1]$ and $\mu$ satisfies \eqref{eq: holder} for some $C>0$ and $s>3/4$.  For an arbitrary constant $K>0$ we then have
$$\E(\hat U_n(a))-g(a)=o(n^{-1/2})$$
where the small-$o$ term is uniform in $a\in {\cal J}_n:=[\lambda(1)+Kn^{-1/6}\log n,\lambda(0)-Kn^{-1/6}\log n].$ 
\end{theorem}

\paragraph{\bf Proof} 
We first localize.   For a given $a$ we define 
\begin{equation}\label{eq: Unhathat}
\hat{\hat U}_n(a)=\argmax_{|u-g(b)|\leq T_nn^{-1/3},\ u\in[0,1]}
\left\{\Lambda_n(u)-au\right\}
\end{equation}
with $T_n=n^{\eps}$  and $b$ a random variable such that $b=a+O_{p}(n^{-1/2})$. Here, $\eps>0$ is arbitrarily small. The variable $b$ will be chosen in a convenient way later.  Note that $\hat{\hat U}_{n}(a)$ is defined in a similar way as $\hat U_{n}(a)$, see \eqref{UnRandom}, but with the location of the maximum taken on a shrinking neighborhood of $g(b)$ instead of being taken over the whole interval $[0,1]$.  Although it may seem more natural to consider $b=a$, we will see  that this choice is not the better one to derive precise bounds on the bias of $\hat U_{n}(a)$. For notational convenience, we do not make it explicit in the notation that $\hat{\hat U}_n(a)$ depends on $b$. The following lemma makes the connection between the bias of $\hat U_{n}(a)$ and that of the localized version; it is proved in Appendix~\ref{Appendix}.

\begin{lemma}\label{lemma: local}
Assume {\monorandom}, {\densityUI}  and {\momentq}. Let $a\in\R$ and let $b$ be a random variable such that
\begin{equation}\label{eq: atob}
\P(|a-b|>x)\leq K_{1}\exp(-K_{2}nx^2)
\end{equation}
for all $x>0$ where $K_{1}$ and $K_{2}$ depend only on $f$, $\mu$ and $\sigma$. We then have
$\E|\hat U_n(a)-\hat{\hat U}_n(a)|=o(n^{-1/2})$
uniformly in $a\in\R$. 
\end{lemma}

In the sequel, we use the notation 
\begin{equation}\label{eq: L}
L(t)=\int_0^tv^2\circ F^{-1}(u)\d u\mbox{ for }t\in[0,1].
\end{equation}
We recall moreover that the notation ${\cal J}_{n}$ has been defined in Theorem \ref{theo: randomU}. We use $L$ to normalize $\hat{\hat U}_n(a)$. This is done in the following lemma, which proof is given in Appendix~\ref{Appendix}. Thanks to the normalization with $L$, $\hat{\hat U}_n(a)$ can be approached by a drifted Brownian motion, see \eqref{eq: argmaxRn}.

\begin{lemma}\label{lem: changeL}
Assume {\monorandom}, {\density}, {\sigmamin} and {\momentq}.  
Let  $a\in{\cal J}_{n}$ and let $b$ be such that \eqref{eq: atob} holds for all $x>0$, where $K_{1}$ and $K_{2}$ depend only on $f$, $\mu$ and $v$. Assume, furthermore, that $\E(b)=a+o(n^{-1/2})$ and  that  $v^2$ and $\mu$ have a continuous first derivative on $[0,1]$. We then have
$$\E(\hat{\hat U}_n(a)-g(a))=\E\left(\frac{L(\hat{\hat U}_n(a))-L(g(b))}{L'(g(a))}\right)+o(n^{-1/2})$$
where the small-$o$ term is uniform in $a\in{\cal J}_{n}$. 
\end{lemma}

Let 
\begin{equation}\label{eq: phin}
\phi_n(t)=\frac{L''(t)}{\sqrt nL'(t)}B_n(t)
\end{equation}
where $B_{n}$ and $L$  are taken from (\ref{kmtinv}) and \eqref{eq: L} respectively.
Moreover, let $A_{n}$ be the event that all inequalities in \eqref{An1} and \eqref{An2} below hold true : 
\begin{equation}\label{An1}\sup_{u\in[0,1]} |B_n(u)|\leq\log n,\quad\sup_{|u-v|\leq T_nn^{-1/3}
\sqrt{\log n}} |B_n(u)-B_n(v)|\leq \sqrt T_n n^{-1/6}\log n,
\end{equation}
\begin{equation}\label{An2}
\quad\sup_{u\in[0,1]}\left|F^{-1}_{n}(u)-F^{-1}(u)-\frac{1}
{\sqrt n f(F^{-1}(u))}B_n(u)\right|\leq n^{\delta-1},
\end{equation}
where we recall that $T_n=n^\eps$ for an arbitrarily small $\eps>0$, and where $\delta\in(0,1/3)$ can be chosen as small as we wish.  We will prove below that
$\P(A_{n})\to 1\mbox{ as }n\to\infty$, see \eqref{eq : P(An)}.
The following lemma is proved in Appendix~\ref{Appendix}.

\begin{lemma}\label{lem: argmax}
Assume {\monorandom}, {\density}, {\sigmamin} and {\momentq}.  
Assume, furthermore, that $v^2$ has a bounded second derivative on $[0,1]$ and $\mu'$ satisfies \eqref{eq: holder} for some $C>0$ and $s>1/2$. Let  $a\in{\cal J}_{n}$ and 
\begin{equation}\label{eq: b}
b=a-\frac{B_{n}(g(a))}{\sqrt n}\lambda'(g(a)).
\end{equation}
Let $q>0$. Then on $A_{n}$, conditionally on $(X_1,\dots,X_n)$, the variable
\begin{equation}\label{eq: norm}
n^{1/3}(L(\hat{\hat U}_n(a))-L(g(b)))
\end{equation}
has the same distribution as 
\begin{equation}\label{eq: argmaxRn}
\argmax_{u\in I_n(b)}
\{D_n(b,u)+W_{g(b)}(u)+R_n(a,b,u)\},
\end{equation}
where for all $t\in[0,1]$,
\begin{equation}\label{eq: Wt}
W_t(u)=\frac{n^{1/6}}{\sqrt{1+\phi_n(t)}}\left[W_n\left(L_n(t)+
n^{-1/3}u(1+\phi_n(t)\right)-W_n(L_n(t))\right],\ u\in\R,
\end{equation}
with $W_n$ being a standard Brownian motion under $\P^X$, 
$$I_{n}(b)=\left[n^{1/3}\left(L(g(b)-n^{-1/3}T_{n})-L(g(b))\right)\,,\, n^{1/3}\left(L(g(b)+n^{-1/3}T_{n})-L(g(b))\right)\right],$$
\begin{eqnarray*}
D_n(b,u)=n^{2/3}\left(\Lambda\circ L^{-1}(L(g(b))+n^{-1/3}u)-\Lambda(g(b))-bL^{-1}(L(g(b))+n^{-1/3}u)+bg(b)\right),
\end{eqnarray*}
and with $T_{n}=n^\eps$ for some sufficiently small $\eps>0$,
\begin{equation}\label{eq: Rn}
\P^X\left(\sup_{u\in I_{n}(b)}|R_{n}(a,b,u)|>x\right)\leq K_{q}x^{-q}n^{1-q/3}
\end{equation}
for all $x>0$, where $K_{q}>0$ does not depend on $n$.
\end{lemma}
It follows from Lemma \ref{lem: argmax} that conditionally on $(X_{1},\dots,X_{n})$, on $A_{n}$ the variable in \eqref{eq: norm} has the same expectation as the variable defined in \eqref{eq: argmaxRn}. The following lemma, which is proved in Appendix~\ref{Appendix},  shows that $R_{n}$ is negligible in \eqref{eq: argmaxRn} in the sense that this expectation, up to a negligible remainder term, is equal to the expectation of the variable
$$V_n(b)=\argmax_{|u|\leq (L'(g(b)))^{4/3}\log n }\{D_n(b,u)+W_{g(b)}(u)\}.$$
\begin{lemma}\label{lem: argmax2}
Assume {\monorandom}, {\density}, {\sigmamin} and {\momentq}.  Assume, furthermore, that $v^2$ has a bounded second derivative on $[0,1]$ and $\mu'$ satisfies \eqref{eq: holder} for some $C>0$ and $s>1/2$. 
Let  $a\in{\cal J}_{n}$ and let $b$ be given by \eqref{eq: b}. With $T_{n}=n^\eps$ for some sufficiently small $\eps>0$, there exists $K>0$ such that on $A_{n}$, we have 
\begin{equation}\notag
\left\vert \E^X\left(n^{1/3}(L(\hat{\hat U}_n(a))-L(g(b)))\right)-\E^X(V_n(b))\right\vert \leq Kn^{-1/6}L'(g(b))(\log n)^{-1}.
\end{equation}
\end{lemma}

Next, we give a precise bound for the conditional expectation of $V_{n}(b)$ (see Appendix~\ref{Appendix} for a proof). For this, we assume that $s>3/4$.

\begin{lemma}\label{lem: argmaxasympt}
Assume {\monorandom}, {\density} and {\sigmamin}.  Assume, furthermore, that $v^2$ has a bounded second derivative on $[0,1]$ and $\mu'$ satisfies \eqref{eq: holder} for some $C>0$ and $s>3/4$. 
Let  $a\in{\cal J}_{n}$ and let $b$ be given by \eqref{eq: b}.  With $T_{n}=n^\eps$ for some sufficiently small $\eps>0$, there exists $K>0$ such that on $A_{n}$, we have 
\begin{equation}\notag
\left\vert \E^X(V_n(b))\right\vert \leq Kn^{-1/6}L'(g(b))(\log n)^{-1}.
\end{equation}
\end{lemma}

We are now in a position to prove Theorem \ref{theo: randomU}. Let  $a\in{\cal J}_{n}$ and let $\hat{\hat U}_n(a))$ be defined by \eqref{eq: Unhathat} where $b$ is taken from \eqref{eq: b}. Since $\lambda'$ is bounded, there exists $K>0$ such that 
\begin{equation}\notag
\P(|a-b|>x)\leq \P\left(\sup_{u\in[0,1]}|B_{n}(u)|>Kx\sqrt n\right)\, \mbox{ for all }x>0.
\end{equation}
Then, with the representation $B_n(u)=W(u)-uW(1)$ in distribution of processes, where $W$ is a standard Brownian motion, we conclude from the triangle inequality that
\begin{eqnarray}\notag
\P(|a-b|>x)\leq \P\left(\sup_{u\in[0,1]}|W(u)|>Kx\sqrt n/2\right) = 2\P\left(\sup_{u\in[0,1]}W(u)>Kx\sqrt n/2\right).
\end{eqnarray}
For the last equality, we used symmetry of $W$.
Then, it follows from the exponentiel inequality for the Brownian motion (see {\it e.g.} Proposition 1.8 in \cite{revuz2013continuous})  that \eqref{eq: atob} holds for all $x>0$, where $K_{1}=2$ and $K_{2}$ depends only on $\lambda$. By lemma \ref{lemma: local}, we then have
\begin{equation}\notag 
\E(\hat U_n(a)-g(a))=\E(\hat{\hat U}_n(a)-g(a))+o(n^{-1/2})
\end{equation}
where the small-$o$ term is uniform in $a\in{\cal J}_{n}$.  Since $B_{n}$ is a centered process, we have $\E(b)=a$, so Lemma \ref{lem: changeL} combined with the preceding display ensures that 
\begin{equation}\label{eq: lemlocal}
\E(\hat U_n(a)-g(a))=\E\left(\frac{L(\hat{\hat U}_n(a))-L(g(b))}{L'(g(a))}\right)+o(n^{-1/2})
\end{equation}
uniformly in $a\in{\cal J}_{n}$. Now, conditionally on $(X_{1},\dots,X_{n})$, on $A_{n}$ we have
\begin{equation}\notag
\left\vert \E^X\left(n^{1/3}(L(\hat{\hat U}_n(a))-L(g(b)))\right)-\E^X(V_n(b))\right\vert \leq K_{3}n^{-1/6}L'(g(b))(\log n)^{-1}
\end{equation}
and
\begin{equation}\notag
\left\vert \E^X(V_n(b))\right\vert \leq K_{3}n^{-1/6}L'(g(b))(\log n)^{-1}.
\end{equation}
Here, we use Lemma \ref{lem: argmax2} and Lemma \ref{lem: argmaxasympt} with $A_{n}$ being the event that all inequalities in \eqref{An1} and \eqref{An2}  hold true. It then follows from the triangle inequality that
$$\E\left(\left\vert \E^X\left(n^{1/3}(L(\hat{\hat U}_n(a))-L(g(b)))\right)\right\vert \BBone_{A_{n}}\right)\leq 2K_{3}n^{-1/6}\E\left(L'(g(b))\right)(\log n)^{-1}.$$
But $L'\circ g$ is a Lipschitz function, so we have
\begin{eqnarray*}
\E\left|L'(g(b))-L'(g(a))\right|\leq K_{4}\E|b-a|\leq K_{5}n^{-1/2},
\end{eqnarray*}
using \eqref{eq: momenta-b} together with the Jensen inequality for the last inequality. Using \eqref{eq: minL'} and the two previous displays yields
\begin{equation}\label{eq: onAn}
\E\left(\left\vert \E^X\left(n^{1/3}(L(\hat{\hat U}_n(a))-L(g(b)))\right)\right\vert \BBone_{A_{n}}\right)\leq 3K_{3}n^{-1/6}L'(g(a))(\log n)^{-1}
\end{equation}
for $n$ sufficiently large.
On the other hand, denoting by $\bar{A_{n}}$ the complementary of $A_{n}$, it follows from H\"{o}lder's inequality together with the Jensen inequality that
\begin{equation}\label{eq: onAnbar}
\begin{split}
&\E\left(\left\vert \E^X\left(n^{1/3}(L(\hat{\hat U}_n(a))-L(g(b)))\right)\right\vert \BBone_{\bar{A_{n}}}\right)\\
&\qquad \leq \E^{1/2}\left(n^{1/3}(L(\hat{\hat U}_n(a))-L(g(b)))\right)^2\P^{1/2}(\bar{A_{n}}).
\end{split}
\end{equation}
Then, we derive from  \eqref{eq: momentU} and  \eqref{eq: momenta-b} that the expectation on the right-hand side is finite. Now, consider $\P(\bar{A_{n}})$ on the right-hand side. It follows from the Markov inequality together with Lemma \ref{lem: kmt} that for all $r\geq 1$ we have
{\small \begin{equation*}
\P\left(\sup_{u\in[0,1]}\left|F_{n}^{-1}(u)-F^{-1}(u)-\frac{B_n(u)}
{\sqrt n f(F^{-1}(u))}\right|> n^{\delta-1}\right) \leq K_{6}\left(\log n\right)^r n^{-r\delta} \leq \;  K_{6}\left(n^{-1/6}L'(g(a))(\log n)^{-1}\right)^2
\end{equation*}}
for $n$ sufficiently large, provided that $r>2/(3\delta)$. Using {\it e.g.} Lemma 5.1 in \cite{durot2014kiefer}, since the Brownian motion satisfies the assumption (A2) of that paper with $\tau=1$  (see the proof of Corollary 3.1 in that paper),  we conclude that 
\begin{equation}\label{eq : P(An)}
\P^{1/2}(\bar{A_{n}})\leq K_{7}n^{-1/6}L'(g(a))(\log n)^{-1}
\end{equation}
for $n$ sufficiently large. Hence, \eqref{eq: onAnbar} yields
$$\E\left(\left\vert \E^X\left(n^{1/3}(L(\hat{\hat U}_n(a))-L(g(b)))\right)\right\vert \BBone_{\bar{A_{n}}}\right)\leq K_{8}n^{-1/6}L'(g(a))(\log n)^{-1}.$$
Together with \eqref{eq: onAn}, this yields 
$$\E\left(\left\vert \E^X\left(n^{1/3}(L(\hat{\hat U}_n(a))-L(g(b)))\right)\right\vert \right)\leq (3K_{3}+K_{8})n^{-1/6}L'(g(a))(\log n)^{-1}.$$
Hence, with the Jensen inequality we arrive at
$$\left\vert \E\left(\E^X\left(n^{1/3}(L(\hat{\hat U}_n(a))-L(g(b)))\right) \right)\right\vert\leq (3K_{3}+K_{8})n^{-1/6}L'(g(a))(\log n)^{-1}.$$
This means that
$$\E\left(\frac{n^{1/3}(L(\hat{\hat U}_n(a))-L(g(b)))}{L'(g(a))}\right)=o(n^{-1/6}).$$
Combining this with \eqref{eq: lemlocal} completes the proof of  Theorem \ref{theo: randomU}. \cqfd

\subsubsection{Proof of Theorem \ref{theo: biasDirect}}
%%%%%%%%%%%%%%%%%%%%%%%%%%%%%%
We begin with a lemma whose proof is available in Appendix~\ref{Appendix}.  

\begin{lemma}\label{lem: biasStep1}
 Assume {\monorandom}, {\densityUI} and {\momentexpo}. 
With $K>0$ an arbitrary positive constant, there exists $K_{1},K_{2}$ both positive such that   
\begin{equation}\label{eq: expoBoundm}
\P\left(|\hat \mu_{n}(t)-\mu(t)|> n^{-1/3}\log n\right) \leq K_{1}\exp(-K_{2}(\log n)^3)
\end{equation}
for all $t\in [Kn^{-1/6}\log n,1-Kn^{-1/6}\log n]$, and
$$\E\left(\hat \mu_{n}(t)-\mu(t)\right)=\E\left[\left(\hat \mu_{n}(t)-\mu(t)\right)\BBone_{|\hat \mu_{n}(t)-\mu(t)|\leq n^{-1/3}\log n}\right]+o(n^{-1/2})$$
where the small-$o$ term is uniform in $t\in [Kn^{-1/6}\log n,1-Kn^{-1/6}\log n]$.
\end{lemma}

Now we turn to the proof of Theorem \ref{theo: biasDirect}.
Distinguishing the positive and negative parts of $\hat \mu_{n}(t)-\mu(t)$, we derive from \eqref{fubini} together with Lemma \ref{lem: biasStep1} that $\E\left(\hat \mu_{n}(t)-\mu(t)\right)=I_1-I_{2}+o(n^{-1/2})$ where
$$I_1= \int_{0}^{n^{-1/3}\log n}\P\left(\hat \mu_{n}(t)-\mu(t)\geq x\right)dx\mbox{ and }
I_{2}= \int_{0}^{n^{-1/3}\log n}\P\left(\mu(t)-\hat \mu_{n}(t)>x\right)dx.$$
%\begin{eqnarray*}I_1&=& \int_{0}^{n^{-1/3}\log n}\P\left(\hat \mu_{n}(t)-\mu(t)\geq x\right)dx\\
%I_{2}&=& \int_{0}^{n^{-1/3}\log n}\P\left(\mu(t)-\hat \mu_{n}(t)>x\right)dx.
%\end{eqnarray*}
Consider $I_{1}$. Since $\hat \mu_{n}^{-1}=F_{n}^{-1}\circ \hat U_{n}$, it follows from the switch relation and \eqref{eq: DKWinv} that
\begin{eqnarray*}
I_{1}&=&\int_{0}^{n^{-1/3}\log n}\P\left(\hat \mu_{n}^{-1}(x+\mu(t))\geq t\right)dx\\
&=&\int_{0}^{n^{-1/3}\log n}\P\left(F^{-1}\circ \hat U_{n}(x+\mu(t))\geq t -O\left(n^{-1/2}\log n\right)\right)dx +o(n^{-1/2})\\
%&=&\int_{0}^{n^{-1/3}\log n}\P\left(\hat U_{n}(x+\mu(t))\geq F\left(t -O\left(n^{-1/2}\log n\right)\right)\right)dx +o(n^{-1/2})\\
&=&\int_{0}^{n^{-1/3}\log n}\P\left(\hat U_{n}(x+\mu(t))\geq F(t)-O(n^{-1/2}\log n)\right)dx+o(n^{-1/2}),
\end{eqnarray*}
where the small $o$-term is uniform in $t\in [c_{1},c_{2}]$. We have $g\circ \mu=F$ and $g'\circ \mu=(\lambda'\circ F)^{-1}$ so it follows from the Taylor expansion  that 
$$g(x+\mu(t))=F(t)-\frac{x}{|\lambda'\circ F(t)|}+O(x^{1+s})$$
for all $t\in [c_{1},c_{2}]$ and $x\in[0,n^{-1/3}\log n]$, where $s$ is taken from  \eqref{eq: holder} and $c_1,c_2$ are as in the statement of the theorem. Since $x^{1+s}\leq n^{-1/2}\log n$ for all $x\leq n^{-1/3}\log n$ provided that $n$ is sufficiently large, we conclude that
\begin{eqnarray*}
I_{1}&=&\int_{0}^{n^{-1/3}\log n}\P\left(\hat U_{n}(a_{x})-g(a_{x})>\frac{x}{|\lambda'\circ F(t)|}-O(n^{-1/2}\log n)\right)dx+o(n^{-1/2}),
\end{eqnarray*}
uniformly, where we set $a_{x}=\mu(t)+x$. But it follows from \eqref{UnRandom} together with  \eqref{eq: Unhathat}  that 
\begin{eqnarray}\label{eq: fromHatToHatHat}
 \P\left(\hat{\hat U}_n(a_x)\neq\hat U_n(a_x)\right)&\leq &
\P\left(|{\hat U}_n(a_x)-g(b_x)|> T_nn^{-1/3}
\right)
\end{eqnarray}
for all $x>0$,  where we recall that $T_{n}=n^{\eps}$ for some arbitrarily small $\eps>0$, and $b_x$ satisfies \eqref{eq: atob} with $a$ replaced by $a_x$. Together with
Lemma \ref{lem: expoboundU}, this yields
\begin{eqnarray*}
I_{1}&=&\int_{0}^{n^{-1/3}\log n}\P\left(\hat {\hat U}_{n}(a_{x})-g(a_{x})>\frac{x}{|\lambda'\circ F(t)|}-O(n^{-1/2}\log n)\right)dx+o(n^{-1/2}),
\end{eqnarray*}
uniformly in $t$. Using again \eqref{eq: fromHatToHatHat} and Lemma \ref{lem: expoboundU}, we then derive from \eqref{eq: toL} in Appendix~\ref{Appendix} that
\begin{eqnarray*}
I_{1}&=&\int_{0}^{n^{-1/3}\log n}\P\left(\frac{L(\hat {\hat U}_{n}(a_{x}))-L(g(b_{x}))}{L'(g(a_{x}))}>\frac{x}{|\lambda'\circ F(t)|}-O(n^{-1/2}\log n)\right)dx+o(n^{-1/2}),
\end{eqnarray*}
where $b_{x}$ is given by \eqref{eq: b} with $a$ replaced by $a_{x}$ and $B_{n}$ being taken from Lemma \ref{lem: kmt}. Since $L'\circ g=v^2\circ \mu^{-1}$, we have 
$$\P(L'(g(b_{x}))\leq c_{0}\gamma)\leq \P(\mu^{-1}(b_{x})\leq \gamma)+\P(1-\mu^{-1}(b_{x})\leq \gamma)$$
for all $\gamma>0$ and $x\in(0,n^{-1/3}\log n)$, where $c_{0}$ is taken from {\sigmamin}. Consider the first probability on the right-hand side. Assume that $\gamma>0$ is chosen small enough so that $c_{1}>\gamma$.  By monotonicity of $\mu$ and the definition of $b_{x}$, there exists a positive constant $K_{1}$ such that for $x\in(0,n^{-1/3}\log n]$ we have 
\begin{eqnarray*}
\P(\mu^{-1}(b_{x})\leq\gamma)&\leq&\P\left(\mu(t)+x-\frac{B_{n}(g(x_{a}))}{\sqrt n}\lambda'(g(a_{x})\geq \mu(\gamma)\right)\\
&\leq&\P\left(|B_{n}(g(x_{a}))|\geq K_{1}\sqrt n(c_{1}-\gamma)\right) \; \leq \;  4\exp(-K_{1}^2n(c_{1}-\gamma)^2/2).
\end{eqnarray*}
It can be proved likewise that $\P(1-\mu^{-1}(b_{x})\leq\gamma)\leq 4\exp(-K_{1}^2n(1-c_{2}-\gamma)^2/2)$ provided $\gamma>0$ is chosen sufficiently small so that $c_{2}+\gamma<1$. Hence, we can restrict attention to the event $\{L'(g(b_{x}))>c_{0}\gamma\}$, which mean that $L'(g(b_{x})$ cannot go to zero. Then, using \eqref{th4D02bis} with $\delta=n^{1/3}\gamma_{n}$ for some $\gamma_{n}\in(n^{-1/2}\log n,n^{-1/3}\log n)$ to be chosen later, we have 
\begin{eqnarray*}
I_{1}&=&\int_0^{n^{-1/3}\log n}\E\P^X\left(\frac{n^{-1/3}V_{n}(b_{x})}{L'(g(a_{x}))}>\frac{x}{|\lambda'\circ F(t)|}-O(\gamma_{n})\right)dx+o(n^{-1/2})\\
&&\qquad+\;O\left(n^{-1/3}(\log n)^2n^{(3-q)/(3(q+1))}(n^{1/3}\gamma_{n})^{-3q/(2(q+1))}\right)
\end{eqnarray*}
where $q$ can be chosen arbitrarily large. For arbitrary $\phi>0$ we can choose $q$ large enough so that
\begin{eqnarray*}
I_{1}&=&\int_0^{n^{-1/3}\log n}\E\P^X\left(\frac{n^{-1/3}V_{n}(b_{x})}{L'(g(a_{x}))}>\frac{x}{|\lambda'\circ F(t)|}-O(\gamma_{n})\right)dx+o(n^{-1/2})\\
&&\qquad+\;O(n^{-7/6+\phi}\gamma_{n}^{-3/2-\phi})\\
&=&\int_0^{n^{-1/3}\log n}\E\P^X\left(n^{-1/3}V_{n}(b_{x})>\frac{xv^2(t)}{|\lambda'\circ F(t)|}-O(\gamma_{n})\right)dx+o(n^{-1/2})\\
&&\qquad+\;O(n^{-7/6+\phi}\gamma_{n}^{-3/2-\phi}).
\end{eqnarray*}
Now, using \eqref{eq: VnToV} in Appendix~\ref{Appendix} with $s=1$ and again $\delta=n^{1/3}\gamma_{n}$ we arrive at
\begin{eqnarray*}
I_{1}
&=&\int_0^{n^{-1/3}\log n}\E\P^X\left(n^{-1/3}V(b_{x})>\frac{xv^2(t)}{|\lambda'\circ F(t)|}-O(\gamma_{n})\right)dx+o(n^{-1/2})\\
&&\qquad+\;O(n^{-7/6+\phi}\gamma_{n}^{-3/2-\phi}).
\end{eqnarray*}
Recall that $g\circ \mu=F$ and define 
$Z(t)=\argmax_{u\in\R }\{-d(F(t))u^2+W(u)\},$
where $d=|\lambda'|/(2(L')^2)$ and $W$ is a standard Brownian motion. Then, $Z(t)$ has the same distribution as 
$$ \argmax_{u\in\R }\{-d(F(t))u^2+W_{g(b_{x})}(u)\}$$
under $\P^X$. On the event $\{\sup_{t\in[0,1]}|B_{n}(t)|\leq\log n\}$ we have
$$V(b_{x})= \argmax_{|u|\leq(L'(g(b_{x})))^{4/3}\log n }\{-d(F(t))u^2+W_{g(b_{x})}(u)+R_{n}(u,x,t)\}$$
where
$$\sup_{|u|\leq(L'(g(b_{x})))^{4/3}\log n}|R_{n}(u,x,t)|=O(n^{-s/3}(\log n)^{2+s})$$
uniformly in $t\in[c_{1},c_{2}]$ and $x\in(0,n^{-1/3}\log n)$. It then  follows from Proposition 1 in \cite{durot2002sharp}  (see also the comments just
 above this proposition) that there are versions of $Z(t)$ and $V(b_{x})$,  and  constants  $K_{1},K_{2},K_{3}>0$, such that on the event $\{\sup_{t\in[0,1]}|B_{n}(t)|\leq\log n\}$ and for $n$ sufficiently large, we have
\begin{equation}\notag
\begin{split}
&\P^X\left(|V(b_{x})-Z(t)|>n^{1/3}\gamma_{n}\right)\leq 
\P^X\left(2\sup_{|u|\leq(L'(g(b_{x})))^{4/3}\log n}|R_{n}(u,x,t)|>x(n^{1/3}\gamma_{n})^{3/2}\right)\\
&\qquad+ K_{1}x\log n+2\P^X\left(|Z(t)|>K_{2}\log n\right)
\end{split}
\end{equation}
where $x=K_{3} (n^{1/3}\gamma_{n})^{-3/2} n^{-s/3}(\log n)^{2+s}.$ We can chose $K_{3}$ large enough so that the probability on the right hand side is equal to zero. Hence, there exists $K_4>0$ such that on the event $\{\sup_{t\in[0,1]}|B_{n}(t)|\leq\log n\}$ we have
\begin{eqnarray*}
\P^X\left(|V(b_{x})-Z(t)|>n^{1/3}\gamma_{n}\right)&\leq &K_{4}(n^{1/3}\gamma_{n})^{-3/2}n^{-s/3}(\log n)^{3+s}+2\P^X\left(|Z(t)|>K_{2}\log n\right)\\
&\leq& K_{4}(n^{1/3}\gamma_{n})^{-3/2}n^{-s/3}(\log n)^{3+s}+4\exp(-K_{5}(\log n)^3)
\end{eqnarray*}
for some $K_5>0$. For the last inequality, we used Theorem 4 in \cite{durot2002sharp}. The second term on the right hand side is negligible as compared to the first one, so we conclude that there exists $K_{6}>0$ such that
\begin{eqnarray*}
\P^X\left(|V(b_{x})-Z(t)|>n^{1/3}\gamma_{n}\right)&\leq &K_{6}(n^{1/3}\gamma_{n})^{-3/2}n^{-s/3}(\log n)^{3+s}.
\end{eqnarray*}
Since $s=1$, we obtain 
\begin{eqnarray}\label{eq: I1}\notag
I_{1}=\int_0^{n^{-1/3}\log n}\P\left(n^{-1/3}Z(t)>\frac{xv^2(t)}{|\lambda'\circ F(t)|}-O(\gamma_{n})\right)dx+o(n^{-1/2})+O(n^{-7/6+\phi}\gamma_{n}^{-3/2-\phi}).
\end{eqnarray}
Consider the integral on the right-hand side. There exists $K>0$ such that 
\begin{equation}\notag
\begin{split}
&\int_0^{n^{-1/3}\log n}\P\left(n^{-1/3}Z(t)>\frac{xv^2(t)}{|\lambda'\circ F(t)|}-O(\gamma_{n})\right)dx\\
&\qquad \leq \int_0^{n^{-1/3}\log n}\P\left(n^{-1/3}Z(t)>\frac{(x-K\gamma_{n})v^2(t)}{|\lambda'\circ F(t)|}\right)dx\\
&\qquad \leq \int_0^{n^{-1/3}\log n}\P\left(n^{-1/3}Z(t)>\frac{yv^2(t)}{|\lambda'\circ F(t)|}\right)dy+O(\gamma_{n})
\end{split}
\end{equation}
using the change of variable $y=x-K\gamma_{n}$. Similarly,
\begin{equation}\notag
\begin{split}
&\int_0^{n^{-1/3}\log n}\P\left(n^{-1/3}Z(t)>\frac{xv^2(t)}{|\lambda'\circ F(t)|}-O(\gamma_{n})\right)dx\\
&\qquad \geq  \int_0^{n^{-1/3}\log n}\P\left(n^{-1/3}Z(t)>\frac{yv^2(t)}{|\lambda'\circ F(t)|}\right)dy+O(\gamma_{n})
\end{split}
\end{equation}
and therefore,
\begin{eqnarray}\label{eq: I1}\notag
I_{1}
&=&\int_0^{n^{-1/3}\log n}\P\left(n^{-1/3}Z(t)>\frac{xv^2(t)}{|\lambda'\circ F(t)|}\right)dx+O(\gamma_{n})+O(n^{-7/6+\phi}\gamma_{n}^{-3/2-\phi}).
\end{eqnarray}
We choose $\gamma_{n}$ that approximately realize the best trade-of between the two big-$O$-terms. This means that we choose $\gamma_{n}$ such that $\gamma_{n}=n^{-7/6}\gamma_{n}^{-3/2}$, that is $\gamma_{n}=n^{-7/15}$, we conclude that for arbitrarily small $\phi>0$,
\begin{eqnarray*}
I_{1}&=&\int_0^{n^{-1/3}\log n}\P\left(n^{-1/3}Z(t)>\frac{xv^2(t)}{|\lambda'\circ F(t)|}\right)dx+O(n^{-14/30+\phi}).\end{eqnarray*}
With similar arguments, we obtain that for arbitrarily small $\phi>0$,
\begin{eqnarray*}
I_{2}&=&\int_0^{n^{-1/3}\log n}\P\left(n^{-1/3}Z(t)<-\frac{xv^2(t)}{|\lambda'\circ F(t)|}\right)dx+O(n^{-7/15+\phi}).\end{eqnarray*}
But $Z(t)$ has the same distribution as $-Z(t)$ for all $t$ so the two preceding displays yield that $I_{1}-I_{2}=O(n^{-7/15+\phi}).$ This completes the proof of Theorem \ref{theo: biasDirect}. \cqfd

\section*{Acknowledgements} The authors would like to thank Arijit Chakrabarty, Guang Cheng, Subhashis Ghoshal and Jon Wellner for some illuminating discussions and correspondences. 

\appendix
\section{Appendix}\label{Appendix}

\subsection{Approximation of the quantile empirical function}\label{sec: preliminaries} 

To go from $\hat U_{n}$ to $\hat \mu_{n}^{-1}$ using \eqref{invmm}, we need  a precise bound for the uniform distance between the quantile function $F^{-1}$ and the corresponding empirical quantile function $F_{n}^{-1}$. The bound we use is given in Lemma \ref{lem: DKWinv} below. It compares to the well known Dvoretzky-Kiefer-Wolfowitz \cite{dvoretzky1956asymptotic} inequality for the empirical distribution function. 
\begin{lemma}\label{lem: DKWinv} 
Let $F$ be a distribution function on $\R$ with a density $f$ supported on $[0,1]$ and bounded away from zero on $[0,1]$. Let $F_{n}$ be the empirical distribution function associated with a $n$-sample from  $F$ and let $F_{n}^{-1}$ be the corresponding empirical quantile function. We then have
\begin{equation}\label{eq: DKWinv}
\P\left(\sup_{t\in[0,1]}|F_{n}^{-1}(t)-F^{-1}(t)|>x\right)\leq 4\exp(-2nc^2x^2)
\end{equation}
for all $n$ and $x>0$, where $c$ is a lower bound for $f$. Moreover, for all $p>0$ there exists $K_p>0$ that depends on $c$ and $p$ only, such that for all $n$,
\begin{equation}\label{eq: DKWinvE}
\E\left(\sup_{t\in[0,1]}|F_{n}^{-1}(t)-F^{-1}(t)|^p\right)\leq K_pn^{-p/2}.
\end{equation}
\end{lemma}

\paragraph{Proof}
Since $f$ is supported on $[0,1]$, both $F_{n}^{-1}$ and  $F^{-1}$ take values in $[0,1]$ so the sup-distance between those functions is less than or equal to one. This means that the probability on the left hand side of~\eqref{eq: DKWinv} is equal to zero for all $x\geq 1$. Hence,  it suffices to prove \eqref{eq: DKWinv} for $x\in(0,1)$. As is customary, we use the notation $y_{+}=\max(y,0)$ and $y_{-}=-\min(y,0)$ for all real numbers $y$. This means that $|y|=\max(y_{-},y_{+})$. Recall the switching relation for the empirical distribution and empirical quantile functions: for arbitrary $a\in[0,1]$ and $t\in[0,1]$, we have \begin{equation}\label{eq: switchEmpir}
F_{n}(a)\geq t\Longleftrightarrow a\geq F_{n}^{-1}(t).
\end{equation} For all $x\in(0,1)$ we then have
\begin{eqnarray*}
\P\left(\sup_{t\in[0,1]}(F_{n}^{-1}(t)-F^{-1}(t))_{+}>x\right)& = &\P\left(\exists t\in[0,1]:\ F_{n}^{-1}(t)>x+F^{-1}(t)\right)\\
& = &\P\left(\exists t\in[0,1]:\ t>F_{n}(x+F^{-1}(t))\right).
\end{eqnarray*}
Using $t=F(F^{-1}(t))$ together with the change of variable $u=x+F^{-1}(t)$ we obtain
\begin{eqnarray*}
\P\left(\sup_{t\in[0,1]}(F_{n}^{-1}(t)-F^{-1}(t))_{+}>x\right)
& \leq &\P\left(\exists u>0:\ F(u-x)>F_{n}(u)\right)\\
& = &\P\left(\exists u\in(x,1):\ F(u-x)>F_{n}(u)\right).
\end{eqnarray*}
For the last equality, we use the fact that $F(u-x)\leq 1=F_{n}(u)$ for all $u\geq 1$, and $F(u-x)=0\leq F_{n}(u)$ for all $u\leq x$. With $c$ a lower bound for $f$ we have 
 $F(u-x)< F(u)-cx$ for all $x\in(0,1)$ and $u\in(x,1)$. Combining this to the  previous display yields
 \begin{eqnarray}\notag \label{eq: DKWinv+}
\P\left(\sup_{t\in[0,1]}(F_{n}^{-1}(t)-F^{-1}(t))_{+}>x\right)& \leq &\P\left(\exists u\in(x,1):\ F(u)-F_{n}(u)>cx\right)\\
\notag & \leq &\P\left(\sup_{u\in\R}|F(u)-F_{n}(u)|>cx\right)\\
&\leq&2\exp(-2nc^2x^2).
\end{eqnarray}
For the last inequality, we used Corollary 1 in \cite{massart1990tight}. On the other hand, for all $x\in(0,1)$  we have
\begin{eqnarray*}
\P\left(\sup_{t\in[0,1]}(F_{n}^{-1}(t)-F^{-1}(t))_{-}>x\right)& \leq &\P\left(\exists t\in[0,1]:\ F_{n}^{-1}(t)< F^{-1}(t)-x\right)\\
& \leq &\P\left(\exists u\in(x,1):\ F_{n}^{-1}(F(u))\leq  u-x\right),
\end{eqnarray*}
using the change of variable $u=F^{-1}(t)$. Hence, with the switching relation we obtain
\begin{eqnarray*}
\P\left(\sup_{t\in[0,1]}(F_{n}^{-1}(t)-F^{-1}(t))_{-}>x\right)
&\leq&\P\left(\exists u\in(x,1):\ F(u)\leq F_{n}(u-x)\right)\\
& \leq &\P\left(\exists u\in(x,1):\ F(u-x)+cx < F_{n}(u-x)\right),
\end{eqnarray*}
using that $F(u-x)< F(u)-cx$
 for all $x\in(0,1)$ and $u\in(x,1)$.  Using again Corollary 1 in \cite{massart1990tight} together with the change of variable $v=u-x$, we arrive at
\begin{eqnarray*}
\P\left(\sup_{t\in[0,1]}(F_{n}^{-1}(t)-F^{-1}(t))_{-}>x\right)
& \leq &\P\left(\sup_{v\in\R}|F(v)-F_{n}(v)|>cx\right)\\
&\leq&2\exp(-2nc^2x^2).
\end{eqnarray*}
Combining the previous display  with \eqref{eq: DKWinv+} completes the proof of \eqref{eq: DKWinv} since $|y|\leq y_{-}+y_{+}$ for all $y\in\R$. Then, \eqref{eq: DKWinvE} follows from  \eqref{fubini} combined to \eqref{eq: DKWinv}. \cqfd

The following lemma, which is a consequence of the strong approximation of the uniform quantile process by Brownian Bridges proved in~\cite{csorgo1978strong}, will also be useful.

\begin{lemma}\label{lem: kmt}
Assume  {\density}. Then, there exist versions of $F_{n}$ and the Brownian bridge $B_{n}$ such that for all  $n$ and $r\geq 1$,
\begin{equation}\label{kmtinv}
\E^{1/r}\left[\sup_{y\in[0,1]}\left|F^{-1}_{n}(y)-F^{-1}(y)-\frac{1}
{\sqrt n f(F^{-1}(y))}B_n(y)\right|^r\right]=O
\left(\frac{\log n}{n}\right).
\end{equation}
\end{lemma}

\paragraph{Proof}
With probability one, the empirical distribution function corresponding to $F(X_{1}),$ \dots$,F(X_{n})$ is $F_{n}\circ F^{-1}$ so the corresponding quantile function is $Q_{n}=F\circ F_{n}^{-1}$. Since the random variables $F(X_{1}),\dots,F(X_{n})$ are i.i.d. and uniformly distributed on $[0,1]$, it follows from Theorem 1 in \cite{csorgo1978strong} that there exist versions of $Q_{n}$ and the Brownian bridge $B_{n}$ such that 
\begin{equation}\label{kmt}\P\left(\sup_{y\in[0,1]}\left\vert Q_{n}(y)-y-\frac{1}
{\sqrt n}B_n(y)\right\vert>\frac{A\log n+z}{n}\right)\leq B\exp(-C z)\end{equation}
for all $z$, where $A$, $B$ and $C$ are positive absolute constants. Thanks to \eqref{fubini}, 
integrating the inequality in \eqref{kmt} where we recall that $Q_{n}=F\circ F_{n}^{-1}$, we obtain that
for all $r\geq 1$,
\begin{equation}\label{kmtE}
\E^{1/r}\left[\sup_{y\in[0,1]}\left|F\circ F_{n}^{-1}(y)-y-\frac{1}
{\sqrt n}B_n(y)\right|^r\right]=O\left(\frac{\log n}{n}\right).
\end{equation}

Now,  $F$ is strictly monotone on $[0,1]$ and  has a bounded second derivative, so it follows from the Taylor expansion that for all $y\in[0,1]$,
\begin{eqnarray*}
F\circ F_{n}^{-1}(y)-y&=&F\circ F_{n}^{-1}(y)-F\circ F^{-1}(y)\\
&=&\left(F_{n}^{-1}(y)-F^{-1}(y)\right)f\circ F^{-1}(y)+\frac12\left(F^{-1}_{n}(y)-F^{-1}(y)\right)^2f'(\theta_{y})
\end{eqnarray*}
for some $\theta_{y}$ lying between $F^{-1}(y)$ and $F_{n}^{-1}(y)$. Combining this with \eqref{kmtE} together with the triangle inequality we get
\begin{equation}\notag
\begin{split}
&\E^{1/r}\left[\sup_{y\in[0,1]}\left\vert \left(F^{-1}_{n}(y)-F^{-1}(y)\right)f\circ F^{-1}(y)-\frac{1}
{\sqrt n}B_n(y)\right\vert^r\right]\\
&\qquad \leq O\left(\frac{\log n}{n}\right)+\frac{1}{2}\sup_{t}\vert f'(t)\vert\E^{1/r}\left[\sup_{y\in[0,1]}\left(F^{-1}_{n}(y)-F^{-1}(y)\right)^{2r}\right]
\end{split}
\end{equation}
for all $r\geq 1$.
With \eqref{eq: DKWinvE} we conclude that
\begin{equation}\notag
\E^{1/r}\left[\sup_{y\in[0,1]}\left\vert \left(F^{-1}_{n}(y)-F^{-1}(y)\right)f\circ F^{-1}(y)-\frac{1}
{\sqrt n}B_n(y)\right\vert^r\right]= O\left(\frac{\log n}{n}\right)\end{equation}
for all $r\geq 1$. The lemma then follows, using that $f$ is bounded away from zero. \cqfd

\subsection{Proof of Lemma \ref{lem: expoboundU}}
By definition, both $\hat U_{n}$ and $g$ take values in $[0,1]$, so
$|\hat{U}_{n}(a)-g(a)|\leq 1.$
This means that the probability on the left hand side of \eqref{eq: expoboundU} is equal to zero for all  $x\geq 1$. Moreover, the right-hand side in \eqref{eq: expoboundU} is greater than one for appropriate $K_{1}$ and $K_{2}$ for all $x\leq n^{-1/3}$. Hence,  it remains to prove  \eqref{eq: expoboundU}  for $x\in(n^{-1/3},1)$.

Let $\Lambda$ be defined on $[0,1]$ by 
\begin{equation}\label{eq: Lambda}
\Lambda(t)=\int_{0}^t\lambda(u) du
\end{equation}
where $\lambda=\mu \circ F^{-1}$ on $[0,1]$.  Let $M_{n}=\Lambda_{n}-\Lambda$ where by definition, $\Lambda_{n}$ is linear on $[(i-1)/n,i/n]$ for all $i\in\{1,\dots,n\}$ and with $\eps_{(j)}=Y_{(j)}-m(X_{(j)})$, satisfies for all $i$
\begin{eqnarray}\notag\label{eq: Lambdan}
\Lambda_n\left(\frac{i}{n}\right)
&=&\frac{1}{n}\sum_{j\leq i}\eps_{(j)}+\frac{1}{n}\sum_{j\leq i}\mu \circ F_{n}^{-1}(j/n)\\
&=&\frac{1}{n}\sum_{j\leq i}\eps_{(j)}+\int_{0}^{i/n}\mu\circ F_{n}^{-1}(u)du.
\end{eqnarray}
For the latter equality, we used the fact that $F_{n}^{-1}$ is piecewise constant. 
Let  $d=A_{1}/A_{4}$ so that
$$\sup_{t\in[0,1]}\lambda'(t)<-d.$$
It follows from the Taylor expansion that
\begin{equation}\label{eq: taylorLambda}
\Lambda(u)-\Lambda(g(a))\leq (u-g(a))a-\frac d2(u-g(a))^2
\end{equation}
for all $u\in[0,1]$ and $a\in[\lambda(1),\lambda(0)]$. For the case $a>\lambda(0)$, we have  $g(a)=0$ and therefore, it follows from the Taylor expansion that
\begin{eqnarray*}
\Lambda(u)-\Lambda(g(a))
&\leq& u\lambda(0)-\frac d2u^2,%\\
%&\leq& (u-g(a))a-\frac d2(u-g(a))^2
\end{eqnarray*}
whence the inequality in \eqref{eq: taylorLambda} also holds for all $a>\lambda(0)$. The case $a<\lambda(1)$ can be handled similarly so we conclude that the inequality in \eqref{eq: taylorLambda} holds for all $a\in\R$. Combining this with \eqref{UnRandom}, where (because $\Lambda_{n}$ is piecewise-linear) the maximum is achieved
on the set $\{i/n,\ i=0,\dots,n\}$, we conclude  that for all $a\in\R$ and $x\in(n^{-1/3},1)$,
\begin{eqnarray*}
\P\left(|\hat U_{n}(a)-g(a)|>x\right)&\leq& \P\left(\sup_{i:\ |g(a)-i/n|>x}\{\Lambda_{n}(i/n)-ai/n\}\geq\Lambda_{n}(g(a))-ag(a)\right)\\
&\leq& \P\left(\sup_{i:\ |g(a)-i/n|>x}\{M_{n}(i/n)- M_{n}(g(a))-\frac d2(in^{-1}-g(a))^2\}\geq 0\right).
\end{eqnarray*}
Define
\begin{equation}\notag%\label{def: calEn}
{\cal E}_n(u)=M_n(u)-\int_{0}^u (\mu\circ F_{n}^{-1}(t)-\mu\circ F^{-1}(t))dt
\end{equation}
for all $u\in[0,1]$. Since $F_{n}^{-1}$ is piecewise-constant , this means that ${\cal E}_{n}(0)=0$ and
\begin{equation}\label{eq: calEn}
{\cal E}_{n}\left(\frac in\right)=\frac 1n\sum_{j\leq i}\eps_{(j)}
\end{equation}
for all $i=1,\dots,n$ with linear interpolation between those points. We then have
\begin{eqnarray}\label{eq: majP1P2}
\P\left(|\hat U_{n}(a)-g(a)|>x\right)&\leq& \P_1+\P_{2}
\end{eqnarray}
where
$$\P_{1}= \P\left(\sup_{u:\ |g(a)-u|>x}\{\int_{g(a)}^u (\mu\circ F_{n}^{-1}(t)-\mu\circ F^{-1}(t))dt-\frac d4(u-g(a))^2\}\geq 0\right)$$
and
$$\P_{2}= \P\left(\sup_{i:\ |g(a)-i/n|>x}\{{\cal E}_{n}(i/n)-{\cal E}_{n}(g(a))-\frac d4(in^{-1}-g(a))^2\}\geq 0\right).$$

We first deal with $\P_{1}$. %Denote by $\|\mu'\|_{\infty}$ the sup-norm of $\mu'$, which is finite thanks to {\monorandom}. We then have
Recall that $A_{2}$ is an upper bound for the sup-norm of $\mu'$. Hence, we have, 
$$\int_{g(a)}^u (\mu\circ F_{n}^{-1}(t)-\mu\circ F^{-1}(t))dt\leq A_{2} |u-g(a)|\sup_{t\in[0,1]}|F_{n}^{-1}(t)-F^{-1}(t)|.$$ 
Combining this with Lemma \ref{lem: DKWinv}, we conclude that for all $x\in(n^{-1/3},1)$, we have
\begin{eqnarray*}
\P_{1}&\leq& \P\left(\sup_{t\in[0,1]}|F_{n}^{-1}(t)-F^{-1}(t)|\geq \frac{dx}{4 A_2}\right)\\
&\leq& 4\exp\left(-\frac{n A_3^2d^2x^2}{8 A_2^2}\right).
\end{eqnarray*}
Since $x^2\geq x^3$, this means that
\begin{equation}\label{eq: majP1}
\P_{1}\leq 4\exp(-Knx^3)
\end{equation}
for all $K\leq A_3^2d^2 A_2^{-2}/8$. 

Next consider $\P_{2}$. For $x\in(n^{-1/3},1)$ we have
\begin{eqnarray*}
\P_{2}&\leq& \sum_{k\geq 1}\P\left(\sup_{i:\ |g(a)-i/n|\in(kx,(k+1)x]}\{{\cal E}_{n}(i/n)-{\cal E}_{n}(g(a))-\frac d4(in^{-1}-g(a))^2\}\geq 0\right)\\
&\leq& \sum_{k\geq 1}\P\left(\sup_{i:\ |g(a)-i/n|\leq(k+1)x}\{{\cal E}_{n}(i/n)-{\cal E}_{n}(g(a))\}\geq \frac d4k^2x^2\right).
\end{eqnarray*}
Using that ${\cal E}_{n}$ is piecewise linear and satisfies \eqref{eq: calEn}, we get
$${\cal E}_{n}(g(a))={\cal E}_{n}\left(\frac{\lfloor ng(a) \rfloor}n\right)+\left(g(a)-\frac{\lfloor ng(a) \rfloor}n\right)\eps_{(\lfloor ng(a) \rfloor+1)}$$
where $\lfloor ng(a) \rfloor$ denotes the integer part of $ng(a)$. Combining the two previous displays yields

\begin{equation}\label{eq: majS1S2}
\P_2\leq S_{1}+S_{2}
\end{equation}
where
$$S_{1}=\sum_{k\geq 1}\P\left(\left(\frac{\lfloor ng(a) \rfloor}n-g(a)\right)\eps_{(\lfloor ng(a) \rfloor+1)}\geq \frac d8k^2x^2\right)$$
and
$$S_{2}=\sum_{k\geq 1}\P\left(\sup_{i:\ |g(a)-i/n|\leq (k+1)x}\left\{\sum_{j\leq i}\eps_{(j)}-\sum_{j\leq ng(a)}\eps_{(j)}\right\}\geq \frac {nd}8k^2x^2\right).$$
We will argue conditionally on $(X_{1},\dots,X_{n})$ to deal with $S_{1}$ and $S_{2}$. It follows from the Markov inequality that for all $\theta>0$, $k\geq 1$, $a\in\R$ and $x\in(n^{-1/3},1)$,
\begin{equation}\notag
\begin{split}
&   \P\left(\left(\frac{\lfloor ng(a) \rfloor}n-g(a)\right)\eps_{(\lfloor ng(a) \rfloor+1)}\geq \frac d8k^2x^2\right)\\
&\qquad\qquad\leq \exp\left(- \frac { \theta d}8k^2x^2\right) \E\left(\exp\left(\theta\left(\frac{\lfloor ng(a) \rfloor}n-g(a)\right)\eps_{(\lfloor ng(a) \rfloor+1)}\right)\right)\\
&\qquad\qquad = \exp\left(- \frac{ \theta d}8k^2x^2\right) \E\left[\E^X\left(\exp\left(\theta\left(\frac{\lfloor ng(a) \rfloor}n-g(a)\right)\eps_{(\lfloor ng(a) \rfloor+1)}\right)\right)\right].
\end{split}
\end{equation}
Recall that $X_{(1)}<\dots<X_{(n)}$ is the order statistics corresponding to $X_{1},\dots,X_{n}$ and that $\eps_{(j)}=\eps_{i}$ if $X_{(j)}=X_{i}$. Therefore, it follows from {\momentexpo} that
\begin{eqnarray*}
\E^X\left(\exp\left(\theta\left(\frac{\lfloor ng(a) \rfloor}n-g(a)\right)\eps_{(\lfloor ng(a) \rfloor+1)}\right)\right)&\leq& K\exp\left(\theta^2\left(\frac{\lfloor ng(a) \rfloor}n-g(a)\right)^2\alpha\right)\\
&\leq& K\exp\left(\frac{\theta^2\alpha}{n^2}\right).
\end{eqnarray*}
Combining the two preceding displays yields that for all $\theta>0$,
\begin{equation}\notag
\P\left(\left(\frac{\lfloor ng(a) \rfloor}n-g(a)\right)\eps_{(\lfloor ng(a) \rfloor+1)}\geq \frac d8k^2x^2\right) \leq K\exp\left(- \frac{ \theta d}8k^2x^2+\frac{\theta^2\alpha}{n^2}\right) .
\end{equation}
Choosing $\theta=dk^2x^2n^2/(16\alpha)$ we arrive at
\begin{equation}\notag
\P\left(\left(\frac{\lfloor ng(a) \rfloor}n-g(a)\right)\eps_{(\lfloor ng(a) \rfloor+1)}\geq \frac d8k^2x^2\right) \leq K\exp\left(- \frac{ d^2k^4x^4n^2}{16^2\alpha}\right) .
\end{equation}
Putting this in the definition of $S_{1}$ and using that $k^4\geq k$ for all $k\geq 1$ and $nx\geq 1$ for all $x\in(n^{-1/3},1)$  we conclude that for all $a\in\R$ and $x\in(n^{-1/3},1)$
\begin{eqnarray}\notag \label{eq: S1}
S_{1}&\leq & K\sum_{k\geq 1}\exp\left(- \frac{ d^2kx^3n}{16^2\alpha}\right) \\ \notag
&\leq & K\exp\left(- \frac{ d^2x^3n}{16^2\alpha}\right)  \sum_{k\geq 0}\exp\left(- \frac{ d^2kx^3n}{16^2\alpha}\right)\\
&\leq& K'\exp(-K_{2}nx^3)
\end{eqnarray}
with any finite $K'$ that satisfies
$K'\geq K\sum_{k\geq 0}\exp\left(- \frac{ d^2k}{16^2\alpha}\right)$
and $K_{2}\leq d^2/(16^2\alpha)$. 

Next, consider $S_{2}$. For this task, recall  that conditionally on $(X_{1},\dots,X_{n})$, the variables $\eps_{(1)},\dots,\eps_{(n)}$ are mutually independent. This means that we can use the Doob's inequality: for all $\theta>0$ we have
\begin{equation}\notag
\begin{split}
&\P^X\left(\sup_{i:\ |g(a)-i/n|\leq (k+1)x}\left\{\sum_{j\leq i}\eps_{(j)}-\sum_{j\leq ng(a)}\eps_{(j)}\right\}\geq \frac {nd}8k^2x^2\right)\\
&\qquad\qquad\leq 2 \exp\left(-\frac {\theta nd}8k^2x^2\right)\sup_{i:\ |g(a)-i/n|\leq (k+1)x}\E^X\left[\exp\left(\theta \left(\sum_{j\leq i}\eps_{(j)}-\sum_{j\leq ng(a)}\eps_{(j)}\right)\right)\right]\\
&\qquad\qquad\leq 2K \exp\left(-\frac {\theta nd}8k^2x^2\right)\sup_{i:\ |g(a)-i/n|\leq (k+1)x}\exp\left(\theta^2\alpha|i-\lfloor ng(a)\rfloor|\right),
\end{split}
\end{equation}
using {\momentexpo} for the last inequality. We have 
$$|i-\lfloor ng(a)\rfloor\leq |i-ng(a)|+1\leq |i-ng(a)|+nx$$
for all $x\in(n^{-1/3},1)$ and therefore,
\begin{equation}\notag
\begin{split}
&\P^X\left(\sup_{i:\ |g(a)-i/n|\leq (k+1)x}\left\{\sum_{j\leq i}\eps_{(j)}-\sum_{j\leq ng(a)}\eps_{(j)}\right\}\geq \frac {nd}8k^2x^2\right)\\
&\qquad\qquad\leq 2K \exp\left(-\frac {\theta nd}8k^2x^2+\theta^2\alpha(k+2)nx\right)
\end{split}
\end{equation}
for all $\theta>0$. Choosing $\theta=dk^2x/(16\alpha(k+2))$ and taking the expectation on both sides we arrive at
\begin{equation}\notag
\begin{split}
&\P\left(\sup_{i:\ |g(a)-i/n|\leq (k+2)x}\left\{\sum_{j\leq i}\eps_{(j)}-\sum_{j\leq ng(a)}\eps_{(j)}\right\}\geq \frac {nd}8k^2x^2\right)\\
&\qquad\qquad\leq 2K \exp\left(-\frac {d^2k^4nx^3}{16^2\alpha(k+2)}\right)\\
&\qquad\qquad\leq 2 K\exp\left(-\frac {d^2knx^3}{3\times 16^2\alpha}\right),
\end{split}
\end{equation}
since $3k^3\geq k+2$ for all $k\geq 1$. By definition of $S_{2}$ we then have
\begin{eqnarray*}
S_2&\leq& 2K \sum_{k\geq 1}\exp\left(-\frac {d^2knx^3}{3\times 16^2\alpha}\right)\\
&\leq& 2K \exp\left(-\frac {d^2nx^3}{3\times 16^2\alpha}\right) \sum_{k\geq 0}\exp\left(-\frac {d^2k}{3\times 16^2\alpha}\right)
\end{eqnarray*}
for all $x\in(n^{-1/3},1)$. For all finite $K'$ such that $K'/2K$ is greater than the sum in the previous display, and $K_{2}\leq d^2/(3\times 16^2\alpha)$, we arrive at
\begin{equation}\notag
S_{2}\leq K'\exp(-K_{2}nx^3).
\end{equation}
Combining this with \eqref{eq: S1} and \eqref{eq: majS1S2} yields that $\P_{2}\leq 2K'\exp(-K_{2}nx^3)$ for appropriate $K'$ and $K_{2}$. Combining this with \eqref{eq: majP1} and \eqref{eq: majP1P2} completes the proof of Lemma \ref{lem: expoboundU}.
\cqfd

\subsection{Proof of Lemma \ref{lem: expoboundinvm}}
Similar to the proof of Lemma \ref{lem: expoboundU}, it suffices to prove the inequality for $x\in(n^{-1/3},1)$.
 Since $\mu^{-1}=F^{-1}\circ g$, it follows from  \eqref{invmm} combined to the triangle inequality that for all $a\in\R$,
\begin{eqnarray*}
|\hat \mu_{n}^{-1}(a)-\mu^{-1}(a)|&=&|F_{n}^{-1}(\hat U_{n}(a))-F^{-1}(g(a))|\\
%&\leq &\sup_{t\in[0,1]}|F_{n}^{-1}(t)-F^{-1}(t)|+|F^{-1}(\hat U_{n}(a))-F^{-1}(g(a))|.
&\leq &\sup_{t\in[0,1]}|F_{n}^{-1}(t)-F^{-1}(t)|+A_{3}^{-1}|\hat U_{n}(a)-g(a)|
\end{eqnarray*}
using that  the first derivative of $F^{-1}$ is bounded by $A_{3}^{-1}$. This means that for all $x\in(n^{-1/3},1)$, we have
\begin{equation}
\label{eq: Utominv}
\P\left(|\hat \mu_{n}^{-1}(a)-\mu^{-1}(a)|>x\right)\leq \P\left(\sup_{t\in[0,1]}|F_{n}^{-1}(t)-F^{-1}(t)|>\frac x2\right)+\P\left(|\hat U_{n}(a)-g(a)|>\frac{x A_{3}}{2}\right).
\end{equation}
Combining this with Lemma \ref{lem: expoboundU} together with \eqref{eq: DKWinv}, we arrive at
$$\P\left(|\hat \mu_{n}^{-1}(a)-\mu^{-1}(a)|>x\right)\leq 4\exp\left(-\frac{n A_{3}^2x^2}2\right)+K_{1}\exp\left(-\frac{K_{2} A_{3}^3 nx^3}{8}\right),$$
for some $K_1,K_{2}>0$. Since $x^2\geq x^3$ for all $x\in(n^{-1/3},1)$, this completes the proof of Lemma \ref{lem: expoboundinvm}. \cqfd

\subsection{Proof of Lemma \ref{lem: expoboundU+-}}

We begin with the proof of \eqref{eq: expoboundU+}.
Similar to the proof of Lemma \ref{lem: expoboundU}, it suffices to prove the inequality for $x\in[n^{-1},1]$.
Let $\Lambda$ and $\lambda$ be taken from \eqref{eq: Lambda} and \eqref{eq: lambdavsm} respectively. Let $\Lambda_{n}$ be defined by \eqref{eq: Lambdan} with linear interpolation between the points $0,1/n, 2/n,\dots,n/n$. It follows from \eqref{UnRandom} (where because $\Lambda_{n}$ is piecewise-linear, the maximum is achieved
on the set $\{i/n,\ i=0,\dots,n\}$) together with $\Lambda_{n}(0)=0$,   that for all $a>\lambda(0)$ and $x\in[n^{-1},1]$ we have
\begin{eqnarray*}\notag
\P^X\left(\hat U_{n}(a)\geq x\right)&\leq&  \P^X\left(\sup_{i\geq nx}\{\Lambda_{n}(i/n)-ai/n\}\geq 0\right)\\
&\leq& \P^X\left(\sup_{i\geq nx}\left\{\frac{1}{n}\sum_{j\leq i}\eps_{(j)}-(a-\mu(0))i/n\right\}\geq 0\right)
\end{eqnarray*}
using the monotonicity of $\mu$. Hence,
\begin{eqnarray*}\notag
\P^X\left(\hat U_{n}(a)\geq x\right)&\leq&\sum_{k\geq 1} \P^X\left(\sup_{i\in [knx,(k+1)nx)}\left\{\sum_{j\leq i}\eps_{(j)}-(a-\mu(0))i\right\}\geq 0\right)\\
&\leq&\sum_{k\geq 1} \P^X\left(\sup_{i\leq  (k+1)nx}\left\{\sum_{j\leq i}\eps_{(j)}\right\}\geq  (a-\mu(0))knx\right).
\end{eqnarray*}
Conditionally on $(X_{1},\dots,X_{n})$, the variables $\eps_{(1)},\dots,\eps_{(n)}$ are mutually independent. This means that we can use the Doob's inequality: for all $\theta>0$ we have
\begin{equation}\notag
\begin{split}
&\P^X\left(\sup_{i\leq (k+1)nx}\left\{\sum_{j\leq i}\eps_{(j)}\right\}\geq (a-\mu(0))knx\right)\\
&\qquad\qquad\leq  \exp\left(-\theta (a-\mu(0))knx\right)\sup_{i\leq (k+1)nx}\E^X\left[\exp\left(\theta\sum_{j\leq i}\eps_{(j)}\right)\right]\\
&\qquad\qquad\leq K \exp\left(-\theta (a-\mu(0))knx\right)\exp\left(\theta^2\alpha (k+1)nx\right),
\end{split}
\end{equation}
using {\momentexpo} for the last inequality.  Choosing $\theta=(a-\mu(0))k/(2\alpha (k+1))$  we arrive at
\begin{equation}\notag
\begin{split}
\P^X\left(\sup_{i\leq (k+1)nx}\left\{\sum_{j\leq i}\eps_{(j)}\right\}\geq (a-\mu(0))knx\right)
&\leq K \exp\left(-\frac {(a-\mu(0))^2k^2nx}{4\alpha (k+1)}\right)\\
&\leq K \exp\left(-\frac {(a-\mu(0))^2knx}{8\alpha}\right),
\end{split}
\end{equation}
since $2k\geq k+1$ for all $k\geq 1$. This means that
\begin{eqnarray*}
\P^X\left(\hat U_{n}(a)\geq x\right)&\leq& K \sum_{k\geq 1}\exp\left(-\frac {(a-\mu(0))^2knx}{8\alpha}\right)\\
&\leq&K \exp\left(-\frac {(a-\mu(0))^2nx}{8\alpha}\right) \sum_{k\geq 0}\exp\left(-\frac {(a-\mu(0))^2k}{3\times 16\alpha}\right) \end{eqnarray*}
for all $x\in[ n^{-1},1]$. For all finite $K'$ such that $K'/K$ is greater than the sum in the previous display, and $K_{2}\leq (8\alpha)^{-1}$, we arrive at
\begin{equation}\notag
\P^X\left(\hat U_{n}(a)\geq x\right)\leq K'\exp(-K_{2}(a-\mu(0))^2nx)
\end{equation}
for all $x\in[ n^{-1},1]$. 
This completes the proof of \eqref{eq: expoboundU+} since $\mu(0)=\lambda(F(0))=\lambda(0)$. The inequality in \eqref{eq: expoboundU-} can be proved in a similar way.
\cqfd

\subsection{Proof of Lemma \ref{lem: boundaries}}
Recall that $\hat\lambda_{n}(0)$ is the right-hand slope at point 0 of the least concave majorant of $\Lambda_{n}$, where $\Lambda_{n}(0)=0$ and $\Lambda_{n}$ is piecewise linear and changes its slope only at points in $\{1/n,2/n,\dots,(n-1)/n\}$. This means that for all $x\geq 0$ we have
\begin{eqnarray*}
\P(\hat\lambda_{n}(0)\geq x)&\leq&\P(\exists i\in\{1,\dots,n\}:\ \Lambda_{n}(i/n)\geq xi/n)\\
&\leq&\P\left(\exists i\in\{1,\dots,n\}:\ \frac{1}{n}\sum_{j\leq i}\eps_{(j)} \geq (x-\mu(0))i/n\right),
\end{eqnarray*}
by monotonicity of $\mu$,
where we recall that $\mu(0)=\lambda(0)$.
With similar arguments as for the proof of Lemma \ref{lem: expoboundU+-} we conclude that there exists $K_{1}>0$ and $K_{2}>0$ such that 
\begin{eqnarray*}
\P(\hat\lambda_{n}(0)\geq x)\leq K_{1}\exp\left(-K_{2}(x-\lambda (0))^2\right)
\end{eqnarray*}
for all $x>\lambda (0)$. Here again, we use the notation $y_{+}=\max(y,0)$ and $y_{-}=-\min(y,0)$ for all real numbers $y$.  Combining the preceding display together with \eqref{fubini} and the fact that a probability is less than or equal to one yields
\begin{eqnarray*}
\E(\hat\lambda_{n}(0)_{+})^p&=&\int_{0}^\infty\P(\hat\lambda_{n}(0)\geq x)px^{p-1}dx\\
&\leq &\lambda^p(0)+ K_{1}\int_{\lambda(0)}^\infty\exp\left(-K_{2}(x-\lambda(0))^2\right)px^{p-1}dx\\
&\leq&A_{5}^p(0)+ K_{1}\int_{0}^\infty\exp\left(-K_{2}x^2\right)p(x+ A_5)^{p-1}dx.
\end{eqnarray*}
The integral on the right hand side is finite so we conclude that $\E(\hat\lambda_{n}(0)_{+})^p\leq K_{3}$ for some $K_{3}>0$. It can be proved likewise that $\E(\hat\lambda_{n}(1)_{-})^p \leq K_3$. Then by monotonicity,
\begin{eqnarray*}\E|\hat\lambda_{n}(0)|^p&\leq& \E(\hat\lambda_{n}(0)_{+})^p+ \E(\hat\lambda_{n}(0)_{-})^p\\
&\leq& \E(\hat\lambda_{n}(0)_{+})^p+ \E(\hat\lambda_{n}(1)_{-})^p
%\\&\leq &2K_{3}.
\end{eqnarray*}
which is at most $2K_3$. Likewise, $\E|\hat\lambda_{n}(1)|^p \leq 2K_3$, which completes the proof. \cqfd

\subsection{Proof of Lemma \ref{lemma: ConnectingBias2}}\label{sec: proof ConnectingBias2}
It follows from \eqref{eq: lambdavsm} that  with generalized inverses, $\mu^{-1}=F^{-1}\circ g$  on $ \R.$
Combined with \eqref{invmm}, this yields
\begin{equation}\label{eq: biasminv}
\begin{split}
&\E\left(\hat \mu_{n}^{-1}(a)-\mu^{-1}(a)\right)=\E\left(F_n^{-1}(\hat U_{n}(a))-F^{-1}(g(a))\right)\\
&\qquad=\E\left(F^{-1}(\hat U_{n}(a))-F^{-1}(g(a))\right)+\E\left(F_{n}^{-1}( \hat U_{n}(a))-F^{-1}(\hat U_{n}(a))\right).
\end{split}
\end{equation}
Consider the first term on the right hand side. Since $F$ has a density function $f$ that is bounded away from zero with a bounded first derivative, see {\density}, it follows from the Taylor expansion that there exists $\theta_{n}$ lying between $\hat U_{n}(a)$ and $g(a)$ such that
$$F^{-1}( \hat U_{n}(a))-F^{-1}( g(a))=\frac{\hat U_{n}(a)-g(a)}{f(F^{-1}(g(a)))}-\frac{f'(F^{-1}(\theta_{n}))}{2(f(F^{-1}(\theta_{n})))^3}(\hat U_{n}(a)-g(a))^2.$$
Hence,
$$\E\left\vert F^{-1}(\hat U_{n}(a))-F^{-1}(g(a))-\frac{\hat U_{n}(a)-g(a)}{f(F^{-1}(g(a)))}\right\vert\leq \frac{\sup_{t}\vert f'(t)\vert}{2(\inf_{t }f(t))^3}\; \E(\hat U_{n}(a)-g(a))^2.$$
It follows from Lemma \ref{lem: expoboundU} combined to \eqref{fubini} that the right-hand side if of maximal order $n^{-2/3}$ uniformly in $a$, whence
\begin{equation}\label{eq: biasminv1}
\E\left( F^{-1}(\hat U_{n}(a))-F^{-1}(g(a))\right)=\E\left(\frac{\hat U_{n}(a)-g(a)}{f(F^{-1}(g(a)))}\right)+O(n^{-2/3})
\end{equation}
uniformly in $a\in\R$. 
Next, consider the second term on the right hand side of \eqref{eq: biasminv}. By Lemma \ref{lem: kmt}, there are versions of $F_{n}$ and the Brownian bridge $B_{n}$ such that
\begin{equation}\label{eq: lem5}
\E\left(F_{n}^{-1}(\hat U_{n}(a))-F^{-1}(\hat U_{n}(a))\right) =O\left(\frac{\log n}n\right)+\E\left(\frac{B_n(\hat U_{n}(a))}
{\sqrt n f(F^{-1}(\hat U_{n}(a)))}\right)
\end{equation}
where the big-$O$ term is uniform in $a\in\R$. 
Now, it follows from Taylor expansion and H\"{o}lder's inequality that
\[\begin{split}
&
\left\vert\E\left(\frac{B_n(\hat U_{n}(a))}{\sqrt n f(F^{-1}(\hat U_{n}(a)))}\right)-\E\left(\frac{B_n(\hat U_{n}(a))}{\sqrt n f(F^{-1}(g(a)))}\right)\right\vert\\
&\qquad\leq\frac{\sup_{t}\vert f'(t)\vert}{\sqrt n(\inf_{t}f(t))^3}\E\left\vert(\hat U_{n}(a)-g(a))B_n(\hat U_{n}(a))\right\vert\\
&\qquad\leq \frac{\sup_{t}\vert f'(t)\vert}{\sqrt n(\inf_{t}f(t))^3}\E^{1/2}\left(\hat U_{n}(a)-g(a)\right)^2\E^{1/2}\left(\sup_{t\in[0,1]}\vert B_n(t)\vert \right)^2
%\\&\qquad=O(n^{-5/6})
\end{split}
\]
which is of order $O(n^{-5/6})$ uniformly in $a\in\R$. Together with \eqref{eq: lem5}, this implies that
$$\E\left(F_{n}^{-1}(\hat U_{n}(a))-F^{-1}(\hat U_{n}(a))\right) =O(n^{-5/6})+\E\left(\frac{B_n(\hat U_{n}(a))}
{\sqrt n f(F^{-1}(g(a)))}\right).$$
Next, we have
$$\left|\E\left(B_n(\hat{U}_n(a)) - B_n(g(a))\right)\right|\leq \E\left(\sup_{|u-g(a)|\leq n^{-1/3}\log n}|B_n(u) - B_n(g(a))|\right)$$
$$+2\E\left(\sup_{u\in[0,1]}|B_n(u)|1(|\hat{U}_n(a)-g(a)|>n^{-1/3}\log n)\right).$$
The first expectation on the right hand side tends to zero by rescaling the Brownian motion $W_n$ and the representation $B_n(t)=W_n(t)-tW_n(1)$ in distribution. For the second expectation, use H\"older's inequality together with Lemma \ref{lem: expoboundU} to conclude that it tends to zero as well as $n\to\infty$. We conclude that 
\begin{eqnarray}\notag\label{eq: biasminv2}
\E\left(F_{n}^{-1}(\hat U_{n}(a))-F^{-1}(\hat U_{n}(a))\right)& =&o(n^{-1/2})+\frac{\E\left(B_n(g(a))\right)}
{\sqrt n f(F^{-1}(g(a)))}\\
&=& o(n^{-1/2})
\end{eqnarray}
uniformly in $a\in\R$. 
For the last equality, we simply used the fact that $B_n$ is a centered process. 
Combining together \eqref{eq: biasminv}, \eqref{eq: biasminv1} and \eqref{eq: biasminv2} completes the proof. \cqfd

\subsection{Proof of Lemma \ref{lemma: local}}
By Fubini's theorem we have
\begin{eqnarray*}
\E\vert\hat{\hat U}_n(a)-\hat U_n(a)\vert&=&\int_{0}^\infty\P\left(\vert\hat{\hat U}_n(a)-\hat U_n(a)\vert>x\right)dx.
\end{eqnarray*}
But it follows from \eqref{UnRandom} together with  \eqref{eq: Unhathat}  that for all $x>0$, 
\begin{eqnarray*}%\label{eq: fromHatToHatHat}
\P\left(\vert\hat{\hat U}_n(a)-\hat U_n(a)\vert>x\right)\leq \P\left(\hat{\hat U}_n(a)\neq\hat U_n(a)\right)&\leq &
\P\left(|{\hat U}_n(a)-g(b)|> T_nn^{-1/3}
\right).
\end{eqnarray*}
Hence, $|\hat{\hat U}_n(a)-g(b)|\leq |{\hat U}_n(a)-g(b)|$ and we obtain
\begin{eqnarray}\label{eq: localE}\notag
&&\E\vert\hat{\hat U}_n(a)-\hat U_n(a)\vert\\
&&\qquad \leq\int_{0}^\infty\min\left\{\P\left(2\vert{\hat U}_n(a)-g(b)\vert>x\right)\,;\, \P\left(|{\hat U}_n(a)-g(b)|> T_nn^{-1/3}
\right)\right\}dx\\ \notag
&&\qquad\leq2T_{n}n^{-1/3}\P\left(|{\hat U}_n(a)-g(b)|> T_nn^{-1/3}
\right)+2\int_{T_{n}n^{-1/3}}^\infty\P\left(\vert{\hat U}_n(a)-g(b)\vert>x\right)dx.
\end{eqnarray}
For all $x>0$ we have
\begin{eqnarray*}
\P\left(|{\hat U}_n(a)-g(b)|> x
\right)&\leq &\P\left(|{\hat U}_n(a)-g(a)|> \frac{x}{2}
\right)+\P\left(K|a-b|>\frac x2\right)
\end{eqnarray*}
for some $K>0$, using that  $g$ is Lipshitz on $\R$.
Using  \eqref{eq: atob} and Lemma \ref{lem: expoboundU}, we conclude that there exist positive constants $K_{1}$ and $K_{2}$ such that 
\begin{eqnarray*}
\P\left(|{\hat U}_n(a)-g(b)|> x
\right)&\leq &K_{1}\exp\left(-K_2nx^3\right)+K_{1}\exp\left(-K_2nx^2\right).
\end{eqnarray*}
Hence, for all $x\leq 1$ we have
\begin{eqnarray}\label{eq: boudUnab}
\P\left(|{\hat U}_n(a)-g(b)|> x
\right)&\leq &2K_{1}\exp\left(-K_2nx^3\right).
\end{eqnarray}
The previous  inequality is trivially true for $x>1$ since in that case, the probability on the left-hand side is equal to zero. Hence, the inequality holds for all $x>0$. Hence, it follows from \eqref{eq: localE} that
\begin{eqnarray*}\E\vert\hat{\hat U}_n(a)-\hat U_n(a)\vert&\leq&4K_{1}T_{n}n^{-1/3}\exp\left(-K_2T_{n}^3\right)+4K_{1}\int_{T_{n}n^{-1/3}}^\infty \exp\left(-K_2nx^3\right) dx\\
&=&o(n^{-1/2})\end{eqnarray*}
by definition of $T_n$, and Lemma \ref{lemma: local} follows. \cqfd

\subsection{Proof of Lemma \ref{lem: changeL}}
It follows from the Taylor expansion that  for all $a\in\R$, there exists $\theta_{a}\in[0,1]$ such that 
$$\frac{L(\hat{\hat U}_n(a))-L(g(a))}{L'(g(a))}=\hat{\hat U}_n(a)-g(a)+\frac12(\hat{\hat U}_n(a)-g(a))^2\frac{L''(\theta_{a})}{L'(g(a))}.$$
Since $F^{-1}\circ g=\mu^{-1}$ we have
\begin{equation}\label{eq: minL'R3}
L'(g(a))=v^2\circ \mu^{-1}(a)\geq c_{0}(\mu^{-1}(a)\wedge(1-\mu^{-1}(a)))
\end{equation}
where $c_{0}$ is taken from {\sigmamin}. On the interval $(\mu(1),\mu(0))$, the function $\mu^{-1}$ has a negative first derivative that is bounded away from zero. Denoting by $c>0$ a lower bound for the absolute value of the derivative, we have
$$\mu^{-1}(a)-\mu^{-1}(\lambda(0))=\int_a^{\lambda(0)}|(\mu^{-1})'(u)|du\geq c(\lambda(0)-a)\geq cKn^{-1/6}\log n$$
for all $a\in{\cal J}_{n}.$ Since $\mu^{-1}(\lambda(0))=0$, we arrive at
$\mu^{-1}(a)\geq cKn^{-1/6}\log n.$
Likewise,
$1-\mu^{-1}(a)\geq cKn^{-1/6}\log n$
for all $a\in{\cal J}_{n}.$ Using \eqref{eq: minL'R3}, this means that
\begin{equation}\label{eq: minL'}
L'(g(a))\geq c_{0}cKn^{-1/6}\log n\mbox{ for all }a\in{\cal J}_{n}.
\end{equation} 
Since, furthermore,  $L''$ is bounded, we conclude that there exists $K>0$ such that
\begin{equation}\label{eq: toL}
\left\vert\frac{L(\hat{\hat U}_n(a))-L(g(a))}{L'(g(a))}-(\hat{\hat U}_n(a)-g(a))\right\vert\leq \frac{Kn^{1/6}}{\log n}(\hat{\hat U}_n(a)-g(a))^2.
\end{equation}
Repeating the same arguments as in the proof of Lemma \ref{lemma: local}, it can be seen that for all $p\geq 1$,
\begin{eqnarray*}\E\vert\hat{\hat U}_n(a)-\hat U_n(a)\vert^p=o(n^{-p/3})\end{eqnarray*}
uniformly in $a$, so it follows from Lemma \ref{lem: expoboundU} combined with the triangle inequality that  for all $p\geq 1$, there exists $K_{p}>0$ such that
\begin{equation}\label{eq: momentU}
\E|\hat{\hat U}_n(a)-g(a)|^p\leq K_{p}n^{-p/3}
\end{equation}
for all $a\in\R$.  With $p=2$, we conclude from \eqref{eq: toL} that
\begin{equation}\label{eq: changeL1}
\E(\hat{\hat U}_n(a)-g(a))=\E\left(\frac{L(\hat{\hat U}_n(a))-L(g(a))}{L'(g(a))}\right)+o(n^{-1/2})
\end{equation}
uniformly for $a\in{\cal J}_{n}$.
Now, it follows from \eqref{fubini} that
\begin{eqnarray}\label{eq: momenta-b}\notag
\E(a-b)^2 &=&\int_{0}^\infty\P(|b-a|>\sqrt x)dx\\
&\leq&\int_{0}^\infty K_{1}\exp(-K_{2}nx)dx=\frac{K_{1}}{nK_{2}}
\end{eqnarray}
where $K_{1}$ and $K_{2}$ are taken from \eqref{eq: atob}. Moreover, from the Taylor expansion it follows that  for all $a\in\R$, there exists $\theta_{a}$ lying between $g(b)$ and $g(a)$ and $\eta_{a}$ lying between $ b$ and $a$ such that 
\begin{eqnarray*}
\frac{L(g(b))-L(g(a))}{L'(g(a))}
&=&(b-a)g'(\eta_{a})+\frac12((b-a)g'(\eta_{a}))^2\frac{L''(\theta_{a})}{L'(g(a))}.
\end{eqnarray*}
The function $L$ has a bounded second derivative on $[0,1]$ so using \eqref{eq : gLips}, \eqref{eq: momenta-b} and \eqref{eq: minL'}, we then arrive at
\begin{eqnarray}\notag\label{eq: devL}
\E\left(\frac{L(g(b))-L(g(a))}{L'(g(a))}\right)&=&\E((b-a)g'(\eta_{a}))+o(n^{-1/2})\\ \notag
&=&\E\left((b-a)(g'(\eta_{a})-g'(a))\right)+o(n^{-1/2}).
\end{eqnarray}
For the last equality, we used the assumption that $\E(b)=a+o(n^{-1/2})$. Consider the expectation on the right hand side. It follows from H\"older's inequality together with~\eqref{eq: momenta-b} that
\begin{eqnarray*}
\E\left|(b-a)(g'(\eta_{a})-g'(a))\right|
&\leq&n^{-1/2}\sqrt{\frac{K_{1}}{K_{2}}}\E^{1/2}\left(g'(\eta_{a})-g'(a)\right)^{1/2}.
\end{eqnarray*}
On the other hand, it follows from the Borel-Cantelli Lemma together with \eqref{eq: atob} that $b$ converges to $a$ as $n\to\infty$ with probability one. Since $g'$ is continuous on $(\lambda(1),\lambda(0))$, this implies that $g'(\eta_{a})$ converges to $g'(a)$ as $n\to\infty$ with probability one. Since $g'$ is bounded, it then follows from the dominated convergence theorem that $\E\left(g'(\eta_{a})-g'(a)\right)^{1/2}$ tends to zero as $n\to\infty$. Hence, it follows from the preceding display that
\begin{eqnarray*}
\E\left|(b-a)(g'(\eta_{a})-g'(a))\right|=o(n^{-1/2}).
\end{eqnarray*}
Combining this with \eqref{eq: devL} and \eqref{eq: changeL1} completes the proof of Lemma \ref{lem: changeL}. \cqfd

\subsection{Proof of Lemma \ref{lem: argmax}}
The location of the maximum of a process is invariant under addition of constants or multiplication by $n^{2/3}$ so it follows from  \eqref{eq: Unhathat} that 
\begin{equation}\notag
\begin{split}
n^{1/3}(L(\hat{\hat U}_n(a))-L(g(b)))=\argmax_{u\in I_n(b)}
\left\{P_{n}(a,b,u)\right\}
\end{split}
\end{equation}
where for all $a,b,u$, 
{\small \begin{equation}\notag
P_{n}(a,b,u) = n^{2/3}\left\{\Lambda_{n}\circ L^{-1}\left(L(g(b))+n^{-1/3}u\right)-\Lambda_n(g(b))-aL^{-1}\left(L(g(b))+n^{-1/3}u\right)-ag(b)\right\}.
\end{equation}}
Recall \eqref{eq: Lambdan} and $\Lambda_{n}$ linearly interpolates between the points $i/n,$ $i=0,\dots,n$.  The $\eps_i$'s are independent under $\P^X$ and we have
\begin{eqnarray*}
\mathrm{Var}^X\left(\frac{1}{n}\sum_{j\leq i}\eps_{(j)}\right)=\frac{1}{n^2}\sum_{j\leq i}v^2\circ F_{n}^{-1}(j/n)
=\frac{1}{n}\int_{0}^{i/n}v^2\circ F_{n}^{-1}(u)du.
\end{eqnarray*}
With {\momentq}, the function $x\mapsto\E(|\eps_i|^q|X_i=x)$
is bounded on $[0,1]$ with an arbitrary $q>0$. 
It then follows from Theorem A in \cite{sakhanenko2006estimates}  that  there exist a positive constant $C_{q}$, and
versions of $\Lambda_{n}$ and the Brownian motion $W_n$ under $\P^X$,  such that for all $x>0$,
\begin{equation}\label{sak}
\P^X\left[\sup_{t\in[0,1]}\left|\Lambda_n(t)-\int_0^t\mu\circ F_n^{-1}(u)\d u-
\frac{1}{\sqrt n}W_n\left(L_n(t)\right)\right|>x \right]\leq C_{q}n^{1-q}x^{-q},
\end{equation}
where 
\begin{equation}\label{eq: Ln}
L_n(t)=\int_0^tv^2\circ F_n^{-1}(u)\d u.
\end{equation}
For these versions of $\Lambda_{n}$ and $W_{n}$ we have
\begin{equation}\notag
\begin{split}
P_{n}(a,b,u)=& n^{2/3}\left\{\int_{g(b)}^{L^{-1}\left(L(g(b))+n^{-1/3}u\right)}\mu\circ F_n^{-1}(t)\d t-a\left(L^{-1}\left(L(g(b))+n^{-1/3}u\right)-g(b)\right)\right\}\\
&\quad+n^{1/6}\left\{W_{n}\circ L_{n}\circ L^{-1}\left(L(g(b))+n^{-1/3}u\right)-W_{n}\circ L_{n} (g(b))\right\}+R_{n1}(b,u),
\end{split}
\end{equation}
where
\begin{equation}\label{eq: R1}
|R_{n1}(b,u)|\leq 2n^{2/3}\sup_{t\in[0,1]}\left|\Lambda_n(t)-\int_0^t\mu\circ F_n^{-1}(u)\d u-
\frac{1}{\sqrt n}W_n\left(L_n(t)\right)\right|.
\end{equation}
We then have
\begin{equation}\notag
\begin{split}
&P_{n}(a,b,u)=D_{n}(b,u)+W_{g(b)}(u)+R_{n1}(b,u)+R_{n2}(a,b,u)+R_{n3}(b,u),
\end{split}
\end{equation}
where by definition of $\Lambda$,
\begin{equation}\begin{split}\notag
R_{n2}(a,b,u)&=n^{2/3}\int_{g(b)}^{L^{-1}\left(L(g(b))+n^{-1/3}u\right)}(\mu\circ F_n^{-1}(t)-\mu\circ F^{-1}(t))\d t\\
&\qquad -n^{2/3}(a-b)\left(L^{-1}\left(L(g(b))+n^{-1/3}u\right)-g(b)\right)
\end{split}
\end{equation}
and
\begin{equation}\begin{split}\notag
R_{n3}(b,u)&=n^{1/6}\left\{W_{n}\circ L_{n}\circ L^{-1}\left(L(g(b))+n^{-1/3}u\right)-W_{n}\circ L_{n} (g(b))\right\}-W_{g(b)}(u).
\end{split}
\end{equation}
To complete the proof of Lemma \ref{lem: argmax}, it remains to prove that $R_{n}$ satisfies \eqref{eq: Rn} for all $x>0$, where
$$R_{n}(a,b,u)=R_{n1}(b,u)+R_{n2}(a,b,u)+R_{n3}(b,u).$$

To do this, note that  from \eqref{eq: R1} and \eqref{sak}, it follows  that 
$$
\P^X\left(\sup_{u\in I_{n}(b)}|R_{n1}(b,u)|>\frac{x}2\right)\leq C_{q}4^qx^{-q}n^{1-q/3}
$$
for all $x>0$.
Therefore, it remains to prove that there exists $K_{q}>0$ such that
\begin{eqnarray}\label{eq: Rnbis}
\P^X\left(\sup_{u\in I_{n}(b)}|R_{n2}(a,b,u)+R_{n3}(b,u)|>\frac x2\right)&\leq& K_{q}x^{-q}n^{1-q/3}
\end{eqnarray}
for all $x>0$. By choosing $K_{q}\geq 1$, the inequality clearly holds for all $x< n^{-1/3+1/q}$ since for such $x$'s, the bound on the right-hand side is greater than one. Therefore, it remains to prove \eqref{eq: Rnbis} for all $x\geq n^{-1/3+1/q}$.

Consider $R_{n2}$. It follows from the Taylor expansion that
\begin{equation}\begin{split}\notag
R_{n2}(a,b,u)&=n^{2/3}\int_{g(b)}^{L^{-1}\left(L(g(b))+n^{-1/3}u\right)}(F_{n}^{-1}(t)-F^{-1}(t))\mu'\circ F^{-1}(t)\d t\\
&\qquad -n^{2/3}(a-b)\left(L^{-1}\left(L(g(b))+n^{-1/3}u\right)-g(b)\right)+R_{n4}(b,u)
\end{split}
\end{equation}
where
$$R_{n4}(b,u)= n^{2/3}\int_{g(b)}^{L^{-1}\left(L(g(b))+n^{-1/3}u\right)}(F_{n}^{-1}(t)-F^{-1}(t))(\mu'(\theta_{nt})-\mu'\circ F^{-1}(t))\d t$$
for some $\theta_{nt}$ lying between $F^{-1}(t)$ and $F_{n}^{-1}(t)$. But it follows from the definition of $I_n(b)$ together with the monotonicity of $L$ that
\begin{equation}\label{eq: accrL}
\left|L^{-1}\left(L(g(b))+n^{-1/3}u\right)-g(b)\right|\leq n^{-1/3}T_{n}
\end{equation}
for all $u\in I_{n}(b)$, so thanks to  the triangle inequality and \eqref{eq: holder}, we obtain
\begin{eqnarray*}
|R_{n4}(b,u)|
&\leq& Cn^{1/3}T_{n} \sup_{t\in[0,1]}|F_{n}^{-1}(t)-F^{-1}(t))|^{1+s}
\end{eqnarray*}
for all $u\in I_{n}(b)$.
On $A_{n}$, the inequalities in \eqref{An1}, \eqref{An2} hold and therefore, with $\delta<1/2$  in \eqref{An2}, we obtain that there exists $K_{1}>0$ such that 
$$\sup_{u\in I_{n}(b)}|R_{n4}(b,u)|\leq K_{1}n^{1/3}T_{n}(n^{-1/2}\log n)^{1+s}.$$ 
Using again \eqref{An1} and \eqref{An2}  together with the fact that 
$\mu'\circ F^{-1}/f\circ F^{-1}=\lambda'$ where $\mu$ satisfies \eqref{eq: holder},
we arrive at
\begin{equation}\notag
R_{n2}(a,b,u)=n^{2/3}\left(L^{-1}\left(L(g(b))+n^{-1/3}u\right)-g(b)\right)\left(\frac{B_{n}(g(a))}{\sqrt n}\lambda'(g(a))-(a-b)\right)+R_{n5}(b,u)
\end{equation}
where 
\begin{eqnarray*}
\sup_{u\in I_{n}(b)}|R_{n5}(b,u)|&\leq &K_{1}n^{1/3}T_{n}(n^{-1/2}\log n)^{1+s}\\
&&+K_{2}n^{1/3}T_{n}\left(n^{\delta-1}+\frac{\log n}{\sqrt n}(n^{-1/3}T_{n})^s+n^{-1/2}\sqrt {T_{n}}n^{-1/6}\log n\right)
\end{eqnarray*}
for some $K_{2}>0$ that does not depend on $n$. Using  \eqref{eq: b} and the assumption that $\delta<1/3$, we conclude from the two preceding displays that
\begin{eqnarray*}
\sup_{u\in I_{n}(b)}|R_{n2}(a,b,u)|&\leq& \frac{n^{-1/3+1/q}}4
\end{eqnarray*}
for $n$ sufficiently large and $T_{n}=n^\eps$ for a sufficiently small $\eps>0$. This means that for all $x\geq n^{-1/3+1/q}$, 
\begin{eqnarray}\label{eq: Rnter}
\P^X\left(\sup_{u\in I_{n}(b)}|R_{n2}(a,b,u)+R_{n3}(b,u)|>\frac x2\right) \leq  \P^X\left(\sup_{u\in I_{n}(b)}|R_{n3}(b,u)|>\frac x4\right).
\end{eqnarray}
Now, consider $R_{n3}$. By definition of $L_{n}$, on $A_{n}$ we have
\begin{equation}\notag
\begin{split}
&L_n\circ L^{-1}\left(L(g(b))+n^{-1/3}u\right)-L_{n}(g(b))\\
&\qquad =\int_{g(b)}^{L^{-1}\left(L(g(b))+n^{-1/3}u\right)}v^2\circ F_{n}^{-1}(t)\d t\\
&\qquad =\int_{g(b)}^{L^{-1}\left(L(g(b))+n^{-1/3}u\right)}v^2\left(F^{-1}(t)+\frac{B_{n}(g(b))}{\sqrt n f(F^{-1}(g(b)))}\right)\d t +O\left(\frac{T_{n}^{3/2}}{n}\log n\right)\\
&\qquad =\int_{g(b)}^{L^{-1}\left(L(g(b))+n^{-1/3}u\right)}v^2\circ F^{-1}(t)+\frac{B_{n}(g(b))}{\sqrt n f(F^{-1}(g(b)))}(v^2)'\circ F^{-1}(g(b))\d t +O\left(\frac{T_{n}^{3/2}}{n}\log n\right)
\end{split}
\end{equation}
uniformly in $u\in I_{n}(b)$. Here, we used the assumption that the function $v^2$ has a bounded second derivative, together with \eqref{An1} and \eqref{An2} with $\delta<1/3$. By definition of $L$, the second derivative of $L$ is given by  $L''=(v^2)'\circ F^{-1}/f\circ F^{-1}$ and therefore,
\begin{equation}\notag
\begin{split}
&L_n\circ L^{-1}\left(L(g(b))+n^{-1/3}u\right)-L_{n}(g(b))\\
&\qquad =n^{-1/3}u+\left(L^{-1}\left(L(g(b))+n^{-1/3}u\right)-g(b)\right)\frac{B_{n}(g(b))}{\sqrt n}L''(g(b)) +O\left(\frac{T_{n}^{3/2}}{n}\log n\right)\\
&\qquad =n^{-1/3}u\left(1+\frac{B_{n}(g(b))}{\sqrt nL'(g(b))}L''(g(b)) \right)+O\left(\frac{T_{n}^{2}}{n}\right)\\
&\qquad =n^{-1/3}u\left(1+\phi_{n}(g(b)) \right)+O\left(\frac{T_{n}^{2}}{n}\right)
\end{split}
\end{equation}
where the big-$O$ term is uniform in $a$ and $u\in I_{n}(b)$. Here, we used that $u=O(T_{n}L'(g(b)))$ uniformly on $I_{n}(b)$. This means that
$$\vert L_n\circ L^{-1}\left(L(g(b))+n^{-1/3}u\right)-\left(L_{n}(g(b))+n^{-1/3}u\left(1+\phi_{n}(g(b))\right) \right)\vert \leq \frac{T_{n}^{2}}{n}(\log n)^2$$
provided that $n$ is sufficiently large.  By definition of $R_{n3}$, we then get
\begin{equation}\begin{split}\notag
R_{n3}(b,u)&=R_{n6}(b,u)+R_{n7}(b,u)
\end{split}
\end{equation}
where for all $u\in I_{n}(b)$,
\begin{eqnarray*}
|R_{n6}(b,u)|&=&n^{1/6}\left\vert W_{n}\left(L_{n}(g(b))+n^{-1/3}u\left(1+\phi_{n}(g(b)) \right)\right)-W_{n}\circ L_{n}\circ L^{-1}\left(L(g(b))+n^{-1/3}u\right)\right\vert\\
&\leq &n^{1/6}\sup_{u\in[0,\log n],\, |u-v|\leq T_{n}^{2}n^{-1}(\log n)^2}\vert W_{n}(v)-W_{n}(u) \vert
\end{eqnarray*}
and
\begin{eqnarray*}
|R_{n7}(b,u)|&=&\left\vert n^{1/6}\left\{W_{n}\left(L_{n}(g(b))+n^{-1/3}u\left(1+\phi_{n}(g(b)) \right)\right)-W_{n}\circ L_{n} (g(b))\right\}-W_{g(b)}(u)\right\vert\\
&\leq&\left |W_{g(b)}(u)\right |\left(1-\sqrt{1-|\phi_{n}(g(b))|}\right)\\
&\leq & |W_{g(b)}(u)|\times\frac{|\phi_{n}(g(b))|}{\sqrt 2},
\end{eqnarray*}
using that $1-\sqrt{1-x}\leq x/\sqrt 2$ for all $x\in(0,1/2]$ together with the fact that on $A_{n}$, 
$$|\phi_{n}(g(b))|\leq\frac{(\log n)^2}{L'(g(b))\sqrt n}\leq \frac 12$$
for $n$ sufficiently large. Combining the two previous displays yields
\begin{equation}\notag
|R_{n7}(b,u)|\leq  \sup_{u\in [0,1]}|W_{g(b)}(u)|\times\frac{(\log n)^2}{L'(g(b))\sqrt {2n}}
\end{equation}
for all $u\in I_{n}(b)$. Therefore, for all $x>0$ we have
\begin{equation}\label{eq: R3}
\begin{split}
&\P^X\left(\sup_{u\in I_{n}(b)}|R_{n3}(b,u)|>\frac x4\right)\leq \P^X\left( \sup_{u\in [0,1]}|W_{g(b)}(u)|\frac{(\log n)^2}{L'(g(b))\sqrt {2n}}>\frac x8\right)\\
&\qquad +\P^X\left(n^{1/6}\sup_{u\in[0,\log n],\, |u-v|\leq T_{n}^{2}n^{-1}(\log n)^2}\vert W_{n}(v)-W_{n}(u) \vert>\frac x8\right).
\end{split}
\end{equation}
It follows from \eqref{eq: b} together with \eqref{An1} that on $A_{n}$,
$$|b-a|\leq n^{-1/2}\log n\sup_{u\in[0,1]}|\lambda'(u)|.$$
Similar to \eqref{eq: minL'}, we then have
\begin{equation}\label{eq: minL'b}
L'(g(b))\geq cn^{-1/6}\log n.
\end{equation} 
for some $c>0$ that does not depend on $a,b,n$. Therefore, the first probability on the right-hand side of \eqref{eq: R3} satisfies
\begin{eqnarray*}
 \P^X\left( \sup_{u\in [0,1]}|W_{g(b)}(u)|\frac{(\log n)^2}{L'(g(b))\sqrt {2n}}>\frac x8\right)&\leq&\P^X\left( \sup_{u\in [0,1]}|W_{g(b)}(u)|>\frac {cxn^{-1/6}\sqrt{2n}}{8\log n}\right)\\
 &\leq& 2\exp\left(-\frac {c^2x^2n^{2/3}}{(8\log n)^2}\right)\\
 &\leq& K_{3}x^{-q}n^{1-q/3}
 \end{eqnarray*}
 for some $K_{3}>0$. Here, we use the fact that $W_{g(b)}$ is distributed as a standard Brownian motion under $\P^X$. Since $W_{n}$ is distributed as a standard Brownian motion under $\P^X$, there exist $K_{4}>0$, $K_{5}>0$ and $K_{6}>0$ such that
 \begin{equation}\notag
 \begin{split}
& \P^X\left(n^{1/6}\sup_{u\in[0,\log n],\, |u-v|\leq T_{n}^{2}n^{-1}(\log n)^2}\vert W_{n}(v)-W_{n}(u) \vert>\frac x8\right)\\
 &\qquad\qquad \leq K_{4}T_{n}^{-2}n(\log n)^{-1}\exp\left(-\frac{K_{5}}{(\log n)^2}x^2T_{n}^{-2}n^{2/3}\right)\\
 &\qquad\qquad \leq K_{6}x^{-q}n^{1-q/3}
 \end{split}
 \end{equation}
  for all $x>0$. Combining the two preceding displays with \eqref{eq: R3}, we conclude that 
 $$\P^X\left(\sup_{u\in I_{n}(b)}|R_{n3}(b,u)|>\frac x4\right)\leq (K_{3}+K_{6})x^{-q}n^{1-q/3}$$
 for all $x>0$. Together with \eqref{eq: Rnter}, this proves that  \eqref{eq: Rnbis} holds for all $x\geq n^{-1/3+1/q}$. This completes the proof of Lemma \ref{lem: argmax}. \cqfd

\subsection{Proof of Lemma \ref{lem: argmax2}}
We use Lemma \ref{lem: argmax} with some $q>18$. We assume without loss of generality that the variables in \eqref{eq: norm} and \eqref{eq: argmaxRn} are equal and defined on the same probability space as $V_{n}(b)$. Define
$$\tilde V_{n}(b)=\argmax_{u\in I_n(b)}\{D_n(b,u)+W_{g(b)}(u)\}.$$ 
It  follows from Proposition 1 in \cite{durot2002sharp}  (see also the comments just
 above this proposition) that there exists $K_{1}>0$ such that  for $n$ sufficiently large, we have
\begin{equation}\notag
\begin{split}
&\P^X\left(|n^{1/3}(L(\hat{\hat U}_n(a))-L(g(b)))-\tilde V_n(b)|>\delta\right)\leq 
\P^X\left(2\sup_{u\in I_{n}(b)}|R_{n}(a,b,u)|>x\delta^{3/2}\right)\\
&\qquad+ K_{1}x(L'(g(b)))^{4/3}\log n+\P^X\left(|\tilde V_n(b)|>(L'(g(b)))^{4/3}\log n\right).
\end{split}
\end{equation}
for every pair $(x,\delta)$ that satisfies
\begin{equation}\label{eq: condD02}
\delta\in(0,(L'(g(b)))^{4/3}\log n],\ x>0,\ \frac{(\log n)^3}{ (L'(g(b)))^{4/3}}\leq-\frac{1}{\delta\log(2x\delta)}.\end{equation}
Here, we use the fact that with $T=(L'(g(b)))^{4/3}\log n$, there exists $K_{2}>0$ such that
$$\sup_{|t|\leq T}\left(\frac{\partial}{\partial t}D_{n}(b,t)\right)^2\leq \left(K_{2}\frac{T}{(L'(g(b)))^2}\right)^2\leq \frac{(\log n)^3}{(L'(g(b)))^{4/3}}$$
for $n$ sufficiently large. By definition, $\tilde V_n(b)$ can differ from
$V_n(b)$ only if its absolute value exceeds $(L'(g(b)))^{4/3}\log n$. Hence we get
\begin{equation}\notag
\begin{split}
&\P^X\left(|n^{1/3}(L(\hat{\hat U}_n(a))-L(g(b)))-V_n(b)|>\delta\right)\leq 
\P^X\left(2\sup_{u\in I_{n}(b)}|R_{n}(a,b,u)|>x\delta^{3/2}\right)\\
&\qquad+ K_{1}x(L'(g(b)))^{4/3}\log n+2\P^X\left(|\tilde V_n(b)|>(L'(g(b)))^{4/3}\log n\right).
\end{split}
\end{equation}
for  every pair $(x,\delta)$ that satisfies \eqref{eq: condD02}. Using \eqref{eq: Rn} with $x$ replaced by $x\delta^{3/2}/2$ proves that the first probability on the right-hand side is less than or equal to $K_{3}(x\delta^{3/2})^{-q}n^{1-q/3}$ for some $K_{3}>0$.
Moreover,  $\tilde V_n(b)(L'(g(b)))^{-4/3}$ is distributed as the location of the maximum
of
$$\frac{D_n(b,(L'(g(b)))^{4/3}u)}{(L'(g(b)))^{2/3}}+W(u),$$
 where $W$ is a standard Brownian motion, and
$$\frac{D_n(b,(L'(g(b)))^{4/3}u)}{(L'(g(b)))^{2/3}}\leq -K_{4}u^2$$
for some $K_{4}>0$ that only depends on $\lambda$ and $v^2$. By Theorem 4 in \cite{durot2002sharp}, we then have
\begin{equation}\notag
\P^X\left(\left|\tilde V_n(b)(L'(g(b)))^{-4/3}\right|>\log n\right)\leq 2\exp(-K_{4}^2(\log n)^3/2).
\end{equation}
Therefore,
\begin{equation}\label{th4D02}
\begin{split}
&\P^X\left(|n^{1/3}(L(\hat{\hat U}_n(a))-L(g(b)))-V_n(b)|>\delta\right)\leq 
K_{3}(x\delta^{3/2})^{-q}n^{1-q/3}\\
&\qquad+ K_{1}x(L'(g(b)))^{4/3}\log n+4\exp(-K_{4}^2(\log n)^3/2)
\end{split}
\end{equation}
for every pair $(x,\delta)$ that satisfies \eqref{eq: condD02}. 
For every $\delta>0$, let 
$$x_\delta=n^{(3-q)/(3(q+1))}\delta^{-3q/(2(q+1))}(L'(g(b)))^{-4/(3(q+1))}.$$
Note that $x_\delta$ is defined in such a way that
$$(x\delta^{3/2})^{-q}n^{1-q/3}=x(L'(g(b)))^{4/3}.$$
This means that the first two terms on the right-hand side of \eqref{th4D02} are of the same order of magnitude, up to a $\log n$ factor.
Therefore, using \eqref{eq: minL'b} we get
\begin{equation}\label{th4D02bis}
\P^X\left(|n^{1/3}(L(\hat{\hat U}_n(a))-L(g(b)))-V_n(b)|>\delta\right)\leq 
K_{5}x_\delta(L'(g(b)))^{4/3}\log n 
\end{equation}
for all $\delta\in(\delta_n,n^{-\eps}(L'(g(b)))^{4/3}]$, where $\eps>0$ and $\delta_n=n^{-1/6}L'(g(b))(\log n)^{-1}$. 
Now by definition of the intervals on which the location of the maximum is taken, we have
$$|V_n(b)|\leq (L'(g(b)))^{4/3}\log n$$
and 
\begin{equation}\label{eq: boundUn}
|n^{1/3}(L(\hat{\hat U}_n(a))-L(g(b)))|\leq \left(\sup_{|\theta-g(b)|\leq n^{-1/3}T_{n}}L'(\theta)\right)T_n.
\end{equation}
Hence, by the triangle inequality, there exists $K_6>0$ such that
$$|n^{1/3}(L(\hat{\hat U}_n(a))-L(g(b)))-V_n(b)|\leq K_{6}L'(g(b))T_n.$$
It now follows from Fubini's theorem that
\begin{equation}\notag
\begin{split}
&\E^X|n^{1/3}(L(\hat{\hat U}_n(a))-L(g(b)))-V_n(b)|\\
&\qquad=\int_0^{K_{6}L'(g(b))T_n}\P^X\left(|n^{1/3}(L(\hat{\hat U}_n(a))-L(g(b)))-V_n(b)|>\delta\right)\d
\delta\\
&\qquad \leq\delta_n+ K_{5}\int_{\delta_n}^{n^{-\eps}(L'(g(b)))^{4/3}}
x_\delta(L'(g(b)))^{4/3}\log n\d \delta\\
&\qquad\qquad+K_{5}K_{6}L'(g(b))T_{n}
x_{n^{-\eps}(L'(g(b)))^{4/3}}(L'(g(b)))^{4/3}\log n.
\end{split}
\end{equation}
For the last inequality, we used \eqref{th4D02bis} for $\delta\in(\delta_n,n^{-\eps}(L'(g(b)))^{4/3}]$ together with the fact that for all $\delta>n^{-\eps}(L'(g(b)))^{4/3}$, we have
\begin{equation}\notag
\begin{split}
&\P^X\left(|n^{1/3}(L(\hat{\hat U}_n(a))-L(g(b)))-V_n(b)|>\delta\right)\\
&\qquad \leq \P^X\left(|n^{1/3}(L(\hat{\hat U}_n(a))-L(g(b)))-V_n(b)|>n^{-\eps}(L'(g(b)))^{4/3}\right)\\
&\qquad \leq K_{5}x_{n^{-\eps}(L'(g(b)))^{4/3}}(L'(g(b)))^{4/3}\log n.
\end{split}
\end{equation}
Recall that $L'$ may go to zero at the boundaries of the interval $[0,1]$ but thanks to the assumption {\sigmamin}, we have \eqref{eq: minL'b}. Combining this with the assumption that $q>18$, we conclude that
\begin{equation}\notag
\E^X|n^{1/3}(L(\hat{\hat U}_n(a))-L(g(b)))-V_n(b)|\leq K_{7}n^{-1/6}L'(g(b))(\log n)^{-1}
\end{equation}
for some $K_{7}>0$, uniformly in $a\in{\cal J}_{n}$. This completes the proof of Lemma \ref{lem: argmax2}. \cqfd

\subsection{Proof of Lemma \ref{lem: argmaxasympt}}
First, we approximate $V_{n}(b)$ with the location of the maximum of a Brownian motion with parabolic drift. Define
$$V(b)=\argmax_{|u|\leq (L'(g(b)))^{4/3}\log n }\{-d(g(b))u^2+W_{g(b)}(u)\},$$
where $d=|\lambda'|/(2(L')^2).$ Here, $W_{g(b)}$, $\lambda$ and $L$  are taken from \eqref{eq: Wt}, \eqref{eq: lambdavsm} and \eqref{eq: L} respectively.
Recall that conditionally on $(X_{1},\dots,X_{n})$, $V_n(b)(L'(g(b)))^{-4/3}$ is distributed as the location of the maximum
of
$$\frac{D_n(b,(L'(g(b)))^{4/3}u)}{(L'(g(b)))^{2/3}}+W(u)$$
over $[-\log n,\log n]$, where $W$ is a standard Brownian motion. Hence, we can assume without loss of generality that
$$ V_n(b)(L'(g(b)))^{-4/3}=\argmax_{|u|\leq \log n}\left\{\frac{D_n(b,(L'(g(b)))^{4/3}u)}{(L'(g(b)))^{2/3}}+W(u)\right\}.$$
Likewise, we assume that
\begin{eqnarray*}
V(b)(L'(g(b)))^{-4/3}&=&\argmax_{|u|\leq \log n}\left\{\frac{-d(g(b))}{(L'(g(b)))^{2/3}}\left((L'(g(b)))^{4/3}u\right)^2+W(u)\right\}\\
&=&\argmax_{|u|\leq \log n}\left\{-\frac{|\lambda'(g(b))|}{2}u^2+W(u)\right\}.
\end{eqnarray*}
Note that
$$\left|L^{-1}\left(L(g(b))+n^{-1/3}(L'(g(b)))^{4/3}u\right)-g(b)\right|\leq n^{-1/3}(\log n)^2(L'(g(b)))^{1/3}$$ for all $u$ with $|u|\leq \log n$, provided that $n$ is sufficiently large. It follows from the Tayor expansion that 
\begin{equation}\notag
\begin{split}
\frac{D_n(b,(L'(g(b)))^{4/3}u)}{(L'(g(b)))^{2/3}}&=\frac{n^{2/3}}{2(L'(g(b)))^{2/3}}\left(L^{-1}\left(L(g(b))+n^{-1/3}(L'(g(b)))^{4/3}u\right)-g(b)\right)^2\lambda'(\theta_{nb})
\end{split}
\end{equation}
for some $\theta_{nb}$ satisfying
$$|\theta_{nb}-g(b)|\leq \left|L^{-1}\left(L(g(b))+n^{-1/3}(L'(g(b)))^{4/3}u\right)-g(b) \right|.$$
The second derivative of the function $L^{-1}$ is given by
\begin{eqnarray*}
(L^{-1})''=-\frac{L''\circ L^{-1}}{\left(L'\circ L^{-1}\right)^3}=-\frac{(v^2)'\circ F^{-1}\circ L^{-1}}{(v^2\circ F^{-1}\circ L^{-1})^3f\circ F^{-1}\circ L^{-1}}
\end{eqnarray*}
where $L''$ is a bounded function. Using once more the Taylor expansion together with the fact that $\lambda=\mu\circ F^{-1}$ where $\mu'$ satisfies \eqref{eq: holder},  we conclude that
\begin{equation}\label{eq: VntoV}
\begin{split}
\frac{D_n(b,(L'(g(b)))^{4/3}u)}{(L'(g(b)))^{2/3}}=-\frac{|\lambda'(g(b))|}{2}u^2+O((\log n)^{2+s}n^{-s/3}(L'(g(b)))^{s/3})
\end{split}
\end{equation}
where the big $O$-term is uniform in $|u|\leq \log n$. Then, with similar arguments as for the proof of Lemma \ref{lem: argmax2}, we obtain from \eqref{eq: VntoV} that there exist positive $K_{1}$ and $K_{2}$ such that
\begin{equation}\notag
\begin{split}
&\P^X\left(|V_n(b)-V(b)|(L'(g(b)))^{-4/3}>\delta\right)\leq 
\P^X\left(K_{1}(\log n)^{2+s}n^{-s/3}(L'(g(b)))^{s/3}>x\delta^{3/2}\right)\\
&\qquad+ K_{2}x\log n+\P^X\left(|V(b)|>\log n\right).
\end{split}
\end{equation}
for every pair $(x,\delta)$ that satisfies
$\delta\in(0,\log n],\ x>0$ and $(\log n)^3\leq-(\delta\log(2x\delta))^{-1}.$
But $|V(b)|\leq \log n$ by definition, so setting
$x=2K_{1}(\log n)^{2+s}n^{-s/3}(L'(g(b)))^{s/3},$
the probabilities on the right hand side of the previous display are both equal to zero. This means that
\begin{equation}\label{eq: VnToV}
\begin{split}
&\P^X\left(|V_n(b)-V(b)|(L'(g(b)))^{-4/3}>\delta\right)\leq 
 2K_{2}K_{1}(\log n)^{2+s}n^{-s/3}(L'(g(b)))^{s/3}.
\end{split}
\end{equation}
for all $\delta\in(\delta_n,(\log n)^{-5}]$, where $\delta_n=n^{-1/6}(L'(g(b)))^{-1/3}(\log n)^{-1}$. 
Since by definition,
$$|V_n(b)-V(b)|(L'(g(b)))^{-4/3}\leq |V_n(b)|(L'(g(b)))^{-4/3}+|V(b)|(L'(g(b)))^{-4/3}\leq 2\log n,$$
it follows from Fubini's theorem that
\begin{equation}\notag
\begin{split}
&\E^X|V_n(b)-V(b)|(L'(g(b)))^{-4/3}\\
&\qquad\qquad=\int_0^{2\log n}\P^X\left(||V_n(b)-V(b)|(L'(g(b)))^{-4/3}|>\delta\right)\d
\delta\\
&\qquad\qquad \leq\delta_n+ 
2K_{2}K_{1}(\log n)^{2+s}n^{-s/3}(L'(g(b)))^{s/3}\int_{\delta_n}^{(\log n)^{-5}}\delta^{-3/2}\d \delta\\
&\qquad\qquad\qquad+4K_{2}K_{1}(\log n)^{2+s}n^{-s/3}(L'(g(b)))^{s/3}(\log n)^{15/2}\\
&\qquad\qquad \leq K_3 n^{-1/6}(L'(g(b))^{-1/3}(\log n)^{-1},
\end{split}
\end{equation}
since $s>3/4$.
This means that
$$\E^X\left|V_n(b)-V(b)\right|\leq Kn^{-1/6}L'(g(b))(\log n)^{-1}.$$
Hence, 
$$\left\vert \E^X(V(b))-\E^X(V_n(b))\right\vert \leq Kn^{-1/6}L'(g(b))(\log n)^{-1}.$$
Under $\P^X$, the process $u\mapsto W_{g(b)}(-u)$ is a standard Brownian motion on the real line. Hence, 
$$- V(b)=\argmax_{|u|\leq (L'(g(b)))^{4/3}\log n}\{-d(g(b))u^2+W_{g(b)}(-u)\}$$
has the same distribution as $ V(b)$.  Since this variable has a finite expectation, his means that $\E^X( V(b))=0$, and Lemma \ref{lem: argmaxasympt} follows. \cqfd

\subsection{Proof of Lemma \ref{lem: biasStep1}}
For all $t\in [Kn^{-1/6}\log n,1-Kn^{-1/6}\log n]$ we have
\begin{eqnarray*}
&&\P\left(|\hat \mu_{n}(t)-\mu(t)|> n^{-1/3}\log n\right)\\
&&\qquad\leq\P\left(\hat \mu_{n}(t)>\mu(t)+n^{-1/3}\log n\right)+\P\left(\hat  \mu_{n}(t)<\mu(t)-n^{-1/3}\log n\right)\\
&&\qquad\leq \P\left(\hat \mu_{n}^{-1}(\mu(t)+n^{-1/3}\log n)\geq t \right)+\P\left(\hat \mu_{n}^{-1}(\mu(t)-n^{-1/3}\log n)\leq t \right),
\end{eqnarray*}
using the switch relation for $\hat \mu_{n}$ similar to \eqref{eq: switchU}. Consider the first probability on the right hand side. Note that $\mu(t)+n^{-1/3}\log n\in(\mu(1),\mu(0))$ for sufficiently large $n$. By {\monorandom},  there exists $c>0$ such that
 $$\mu^{-1}(\mu(t)+n^{-1/3}\log n)< t-cn^{-1/3}\log n.$$
Using Lemma \ref{lem: expoboundinvm}, this yields
\begin{eqnarray*}
&&\P\left(\hat \mu_{n}^{-1}(\mu(t)+n^{-1/3}\log n)\geq t \right)\\
&&\qquad \leq \P\left(\hat \mu_{n}^{-1}(\mu(t)+n^{-1/3}\log n)-\mu^{-1}(\mu(t)+n^{-1/3}\log n)> cn^{-1/3}\log n \right)\\
&&\qquad \leq K_{1}\exp(-K_{2}(\log n)^3).
\end{eqnarray*}
Repeating similar arguments for the second probability we conclude that there exist positive constants $K_{1}$ and $K_{2}$ such that 
\eqref{eq: expoBoundm} holds for all $t$. Combining this with H\"older's inequality together with Theorem \ref{theo: UnifIntRD} we obtain that there exists $K_3>0$ such that
\begin{eqnarray*}
&&\E\left[\left(\hat \mu_{n}(t)-\mu(t)\right)\BBone_{|\hat \mu_{n}(t)-\mu(t)|> n^{-1/3}\log n}\right]\\
&&\qquad \leq 
\E^{1/2}\left(\hat \mu_{n}(t)-\mu(t)\right)^2\P^{1/2}\left(|\hat \mu_{n}(t)-\mu(t)|> n^{-1/3}\log n\right)\\
&&\qquad \leq K_{3}n^{-1/3}\exp(-K_{2}(\log n)^3)=o(n^{-1/2})
\end{eqnarray*}
for all $t\in(0,1]$.  This completes the proof of Lemma \ref{lem: biasStep1}. \cqfd

\subsection{Super-Efficiency of the pooled estimator in kernel density estimation}
\label{supeffkernel} 
We describe the phenomenon in a density estimation setting since this is the easiest to deal with. Since isotonic regression is optimal under a once differentiable assumption on the underlying (monotone) function, for a meaningful comparison, we will assume that the smooth density of interest, $f_0$, belongs to the class $\mathcal{F}$ of \emph{all continuously differentiable densities supported on $[0,1]$} that are bounded above by some fixed $M$ and whose derivatives are also uniformly bounded by some $L$.  This is a special case of the generic class of densities considered in Section 1.2 of \cite{Sasha09} with $\beta = \ell = 1$. 

Given i.i.d.~data $X_1, X_2, \ldots, X_N$ from a density $f$, consider the kernel density estimator (KDE) $\hat{f}_n(t_0)$ at some interior point $t_0$ using bandwidth $h_N$ and a continuously differentiable symmetric unimodal (at 0) kernel $K$ supported on $[-1,1]$: 
$$\hat f_N(x) := \frac{1}{n h_N} \sum_{i=1}^N K \left(\frac{x - X_i}{h_N} \right). $$ 
An optimal choice of $h_N$ is $N^{-1/3}$, under a once differentiability assumption on $f$. Suppose that we divide our sample as usual, $N = m  \times n$, and compute $m$ KDEs $\hat{f}_{n,1}, \ldots, \hat{f}_{n,m}$ using bandwidth $h_n = n^{-1/3}$ as above and let $\bar f_N$ denote the pooled estimator, i.e, $$\bar{f}_N(t_0) = \frac{1}{m } \sum_{j=1}^{m } \hat{f}_{j,n}(t_0).$$ Then we have the following result; see Section~\ref{kernel-pool} of the Appendix for a proof.

\begin{lemma}
\label{kernel-pool-lemma}
With $h_N = N^{-1/3}$, we have 
\begin{equation}\label{eq:f_N-hatKDE}
N^{1/3} (\hat{f}_{N}(t_0) - f (t_0)) \stackrel{d}{\to} N(0, f (t_0) R(K)),\qquad \mbox{as } n \to \infty.
\end{equation}
where $R(K) :=  \int K^2(u) du$; and with $m$ fixed and $\overline{f}_N$ as defined above, 
$$N^{1/3} (\bar{f}_N(t_0) - f (t_0)) \stackrel{d}{\to} N(0, m^{-1/3} f(t_0) R(K)),\qquad \mbox{as } n \to \infty.$$ 
Further, let \[\tilde{E}_m: =  \liminf_{N \to \infty}\,\sup_{f \in \mathcal{F}} \E_{f}\left[ N^{2/3}(\bar{f}_N(t_0) - f(t_0))^2 \right],\]
where the subscript $m$ indicates that the maximal risk of the $m$-fold pooled estimator ($m$ fixed) is being considered. Then $\tilde{E}_m \geq m^{2/3}\, {c}_0$ for some ${c}_0 > 0$. If $m = m_n$ is allowed to increase with $n$, then
\[  \liminf_{N \to \infty}\,\sup_{f \in \mathcal{F}} \E_{f}\left[ N^{2/3}(\bar{f}_N(t_0) - f(t_0))^2 \right] = \infty \,.\] 
\end{lemma}
Thus, similar to isotonic regression, the pooled estimator $\bar{f}_N(t_0)$ has lower asymptotic variance by a factor $m^{1/3}$. However, as in isotonic regression, the maximal risk of the pooled estimator suffers.
\subsubsection{Fixing super-efficiency via undersmoothing in KDE} Suppose now that we \emph{change the bandwidth} for each of the subsample based estimators $\hat{f}_{n,j}$'s to ${h}_n = N^{-1/3} \equiv n^{-1/3} m^{-1/3}$, i.e., we \emph{slightly undersmooth} relative to what we were doing above, using the factor $m^{-1/3} < 1$. Then, it is easy to see that:
\[ n^{1/3} (\hat{f}_{n,1}(t_0) - f (t_0))  \stackrel{d}{\to}  N(0, m^{1/3} f (t_0) R(K)) ,\qquad \mbox{as }n \to \infty, \] 
which translates to 
\[ N^{1/3} (\bar{f}_N(t_0) - f (t_0))  \stackrel{d}{\to}  N(0, f (t_0) R(K)) ,\qquad \mbox{as }n \to \infty,\] for $m$ fixed, showing that this new pooled estimator has comparable asymptotic performance to that of the global KDE $\hat{f}_N$; cf.~\eqref{eq:f_N-hatKDE}. Further, we have the following lemma; see Section~\ref{lem:KDE-UnderSmooth} for a proof of the result. 
\begin{lemma}\label{lem:SupEffKDE}
For $h_n = N^{-1/3}$, there exists a constant $C>0$ such that for all $n,m$, \[\sup_{f \in \mathcal{F}} \E_{f}\left[ N^{2/3}(\bar{f}_N(t_0) - f(t_0))^2 \right] \le C.\] Furthermore, even with $N = m_n \times n$ where $m_n$ increases to infinity with $n$, 
\[N^{1/3} (\bar{f}_N(t_0) - f (t_0))  \stackrel{d}{\to}  N(0, f (t_0) R(K)) \,.\] 
\end{lemma}

Thus, by adjusting the bandwidth, we \emph{obviate} the super-efficiency phenomenon. Further, as shown above, super-efficiency can also be fixed with $m_n$ increasing with $n$ by the same choice of bandwidth as in the finite $m$ case. 

\begin{remark}
The phenomenon of matching the pooled estimator's performance to the global estimator in KDE was noted by~\cite{LiEtAl13} in a twice differentiable setting under a regime where $m$ was allowed to increase with $n$, and also observed in a different (but related) context by~\cite{Zhang13}. Note that this is in direct contrast to the isotonic LSE discussed earlier: the bandwidth in the isotonic regression problem is not user-specified but chosen adaptively by the least squares procedure, and therefore does not permit the kind of adjustment that kernel based estimation does where the flexibility of choosing the bandwidth appropriately allows the global KDE match the performance of the pooled estimator while also preventing the super-efficiency phenomenon.  
\end{remark}
\subsubsection{Proof of Lemma~\ref{kernel-pool-lemma}}\label{kernel-pool}
By considerations similar to those in the proof of Theorem~\ref{superefficient}, 
\[ E_m \geq m^{2/3} \liminf_{n} \sup_{f \in \mathcal{F}}\,n^{2/3}(\E_f(\hat{f}_{n,1}(t_0)) - f(t_0))^2 \]
and it suffices to show that the right-side is larger than some positive number. Fix an $f_0$ in the class and consider a sequence of 
densities defined by 
\[ f_n(t) = f_0(t) + n^{-1/3}\,B(n^{1/3}(t-t_0)) \,,\] 
where $B$ is continuously differentiable and vanishes outside of $[-1,1]$. Note that for the $f_n$'s to be densities, $\int_{-1}^1 B(u)du = 0$. To ensure that the $f_n$'s fall within the class $\mathcal{F}$ we need to ensure that $B$ and its derivative $B'$ are uniformly bounded in absolute value by a sufficiently small number. We now consider the sequence $n^{1/3}\,b_{f_n}(t_0)$ under the sequence 
$f_n$ at the point $t_0$, where $b_{f_n}(t_0) = E_{f_n}(\hat{f}_{n,1}(t_0)) - f(t_0)$. Now,
\begin{eqnarray*}
n^{1/3}\,b_{f_n}(t_0) &=& n^{1/3} \int_{-1}^1[f_n(t_0 + u\,n^{-1/3}) - f_n(t_0)] K(u) du \\
& = & n^{1/3} \int_{-1}^1[(f_n - f_0)(t_0 + u\,n^{-1/3}) - (f_n - f_0)(t_0)] K(u) du  \\
& & \qquad \qquad \qquad + n^{1/3} \int_{-1}^1[f(t_0 + u\,n^{-1/3}) - f(t_0)] K(u) du \\
& = & n^{1/3}\,\int_{-1}^1[n^{-1/3}\,B(u) - n^{-1/3}\,B(0)] K(u) du + o(1) \\
& = & \int_{-1}^1 (B(u) - B(0))K(u)\,du + o(1) \\
& \rightarrow &   \int_{-1}^1 B(u)\,K(u)\,du \ne 0\,,
\end{eqnarray*}
\emph{provided} $B(0) = 0$ and $\int_{-1}^1 B(u) K(u)\,du \ne 0$. We can define $B$ as:
\[ B(u) = -C\,[(1/16) - (u + 3/4)^2]^2\,1(-1 \leq u \leq -3/4) \;+ \;C\,[(1/16) - (u + 1/4)^2]^2\,1(-1/2 \leq u \leq 0) \]
\[ \qquad \qquad  + C\,[(1/16) - (u - 1/4)^2]^2\,1(0 \leq u \leq 1/2) \;- \;C\,[(1/16) - (u - 3/4)^2]^2\,1(1/2 \leq u \leq 1) \,.\]
Then $B$ is continuously differentiable, and by manipulating $C$ to depend just on $f_0$, its derivative can be uniformly bounded by as small as number as we like, $B(0) = 0$ and $\int_{-1}^1 B(u) K(u) du > 0$ for any symmetric unimodal (at 0) kernel on $[-1,1]$. Note that:
\[  m^{2/3} \liminf_{n}\sup_{f \in \mathcal{F}} n^{2/3}(\E_f(\hat{f}_{n,1}(t_0)) - f(t_0))^2 \geq m^{2/3} \liminf_{n \rightarrow \infty} n^{2/3}\,b_{f_n}^2(t_0) 
> 0 \,. \;\;\Box \] 

\subsubsection{Proof of Lemma~\ref{lem:SupEffKDE}}\label{lem:KDE-UnderSmooth}
Note that $h_n = N^{-1/3}$. The maximal risk of $\bar{f}_N(t_0)$ over the class $\mathcal{F}$ is bounded in $n,m$ as shown below. For any $f \in \F$, 
\begin{eqnarray*}
	\E_{f}\left[(\bar{f}_N(t_0) - f(t_0))^2 \right] & = &  \mathrm{Var}(\bar{f}_N(t_0)) + [\E_{f}(\overline{f}_N(t_0)) - f(t_0)]^2 \\
	& = & \frac{1}{m} \mathrm{Var}(\hat{f}_{n,1}(t_0)) + [\E_{f}(\hat{f}_{n,1}(t_0)) - f(t_0)]^2 \\
	& \le & \frac{C_1}{m n h_n} + C_2 {h}_n^2  \;\; \le \;\; (C_1 + C_2) N^{-2/3},
\end{eqnarray*}
where $C_1,C_2 >0$ are constants (see e.g., Propositions 1.1 and 1.2 in~\cite{Sasha09}) and we have used the fact that ${h}_n = N^{-1/3} \equiv n^{-1/3} m^{-1/3}$. Therefore, 
\[ \sup_{f \in \mathcal{F}} \E_f [N^{2/3} (\bar{f}_N(t_0) - f(t_0))^2] \leq C_1 + C_2,\] which yields the first part of the result.

Let $\{X_i^{(j)}\}_{i=1}^n$ denotes the $j$'th split-sample, for $j = 1, 2, \ldots, m_n$ (where $N = m_n \times n$), and let $\tilde{h}_n = N^{-1/3}$. Define the KDE for each split-sample and the pooled estimator as
$$\hat f_{j,n}(t_0) := \frac{1}{n\tilde{h}_n} \sum_{i=1}^n K \left(\frac{t_0 - X_i^{(j)}}{\tilde{h}_n} \right), \qquad \mbox{and} \qquad \bar{f}_N(t_0) = \frac{1}{m_n}\,\sum_{j=1}^{m_n}\,\hat{f}_{j,n}(t_0).$$
Consider the distribution of $N^{1/3}(\bar{f}_N(t_0) - f(t_0))$. This quantity can be written as:
\begin{equation}\label{eq: Z}
 N^{1/3}[\E(\bar{f}_N(t_0)) - f(t_0)] + \frac{N^{1/3}}{m_n}\sum_{j=1}^{m_n}\left[\hat{f}_{j,n}(t_0) - \E (\hat{f}_{j,n}(t_0)) \right].
\end{equation}
The first term is simply $N^{1/3}[\E(\hat{f}_{1,n}(t_0)) - f(t_0)]$ and converges to 0, and it remains to find the distribution of the second term which can be written as 
$$ S_N :=\sum_{i=1}^n\sum_{j=1}^{m_n}\,[Z_{n,i}^{(j)} - \E (Z_{n,i}^{(j)})]\qquad \mbox{ with }\qquad  Z_{n,i}^{(j)} =  \tilde h_n  K \left(\frac{t_0 - X_i^{(j)}}{\tilde{h}_n} \right) \,,$$
using that $\tilde h_n=N^{-1/3}$ together with the definition of $\hat f_{j,n}$. Letting $B_n^2 := \mbox{Var}(S_N)$, we can conclude that $S_N/B_n$ converges to $N(0,1)$ provided the Lindeberg condition can be verified. By a straightforward calculation,
\[ B_n^2 = N\tilde h_n^2\mbox{Var}\left(K\left(\frac{t_0-X_1^{(1)}}{\tilde h_n}\right)\right) \rightarrow f(t_0)R(K) \,\]
using again that $\tilde h_n=N^{-1/3}$, and
where $R(K) :=   \int K^2(u) du$.
Thus, subject to the Lindeberg condition being satisfied, $N^{1/3}(\bar{f}_N(t_0) - f(t_0)) \stackrel{d}{\rightarrow} N(0, f(t_0) R(K))$, matching the performance of the global estimator. Since $B_n$ converges to a non-zero limit and the $Z_{n,i}^{(j)}$'s are i.i.d.~it is easy to see that the Lindeberg condition reduces to checking that for any $\eta > 0$, 
\[ N\,\E[Z_{n}^2\,1(|Z_{n}| > \eta)] \rightarrow 0 \]
where $Z_n$ has the same distribution as $Z_{n,i}^{(j)}$ for arbitrary $i$ and $j$. For arbitrary $\eta>0$ we have
\[ N\,\E[Z_{n}^2\,1(|Z_{n}| > \eta)] =\int K^2(u)f_0(t_0-u\tilde h_n)1(|\tilde h_nK(u)|>\eta)du\to0\]
by the dominated convergence theorem, since $\tilde h_n\to0$. Hence, the Lindeberg condition is satisfied. \qed

\subsection{Additional details of the Proof of Theorem ~\ref{superefficient}}
\label{super-efficient-details} 
\subsubsection{Proof of Claim [C]}
We prove the claim for $c=d=1$; the proof for the general case is similar and involves no new ideas. In fact, we show that 
$\E[\arg\min_{t}\,W(t) + t^2 + \mathcal{D}(t)] < 0$. Note that for any integrable  random variable $X$, 
\[ \E(X) = \int_{0}^{\infty}\,\left(\P(X^{+} > x) - \P(X^{-} > x)\right)dx \,.\]
 %= \int_{0}^{\infty}\,\left(\P(X > x) - \P(X < -x)\right)\,dt. \]
We use the following result (see Theorem 1 of \cite{ArijitResult16}) on the minimizer of drifted Brownian motion. 
\begin{theorem}
\label{arijit-theorem} 
Let $\{W(t)\}_{-\infty < t < \infty}$ be a two-sided Brownian motion starting at 0. Let $\psi$ be a continuous symmetric function defined on $\mathbb{R}$ with $\psi(0) = 0$ such that $\psi(t)/t \rightarrow \infty$ as $t \rightarrow \infty$. Let $A$ be a non-decreasing continuous function defined on $[0,\infty)$ with $A(0) = 0$ that is strictly increasing on some interval and with the property that  $A(t)/t \rightarrow 0$ as $t \rightarrow \infty$. For all $t \in \mathbb{R}$, define $\mathcal{A}(t) = A(t)\,1(t \geq 0) - A(-t)\,1(t < 0)$ and consider the drifted process $Z(t) := W(t) + \psi(t) + \mathcal{A}(t)$. Then $M := \arg \min_t Z(t)$ exists almost surely and is unique, and if $M^{+}$ and $M^{-}$ denote its positive and negative parts respectively, then $M^{+}$ is stochastically strictly smaller than $M^{-}$, i.e., $F_{M^{+}}(x) \geq 
F_{M^{-}}(x)$, with strict inequality for some $x$ (and therefore for all $x$ in an interval, by right continuity). 
\end{theorem}
Invoking this theorem with $\psi(t) = t^2$ and $\mathcal{A} = \mathcal{D}$ (which corresponds to $A(t) = \int_0^{t \wedge 1} B(u)du$ for $t > 0$) and with $X \equiv M := \arg\min_{t}\,W(t) + t^2 + \mathcal{D}(t)$, then immediately shows that $\E X < 0$. 

\subsubsection{The Switching Relationship} 
The proceses $X_{c,d,\mathcal{D}}$ and $X_{c,d}$ where $X_{c,d}(t) := c W(t) + d t^2$ for $t \in \mathbb{R}$ both have a unique minimum with probability 1 (\cite{KP90}, Lemma 2.6) and induce mutually absolutely continuous distributions on $C_{min}(\mathbb{R})$ (the space of continuous functions on the real line that possess a unique minimum equipped with the topology of uniform convergence on compact sets). Therefore, the GCM of $X_{c,d,\mathcal{D}}$ has the same almost sure characterization as that of $X_{c,d}$: it is a piecewise linear function that touches $X_{c,d,\mathcal{D}}$ at finitely many points in every compact interval \cite{groeneboom2011vertices}. For a given $\eta$, let 
$\mbox{Largmin}\,X_{c,d,\mathcal{D},\eta}$ denote the largest minimizer of the process $h \mapsto X_{c,d,\mathcal{D}}(h) - \eta\,h$. Then, the switching relationship originally proposed by Groeneboom, which can be readily verified in this case by drawing a simple diagram says:
\[ \mbox{Largmin}\,X_{c,d,\mathcal{D},\eta} \leq x \Leftrightarrow g_{c,d,\mathcal{D}}(x) \geq \eta \,.\]
Since the minimizer of $X_{c,d,\mathcal{D},\eta}$ is almost surely unique, referring to (a version of) this random variable as $\arg \min X_{c,d,\mathcal{D},\eta}$, we have:
\[ \arg \min\,X_{c,d,\mathcal{D},\eta} \leq x \Leftrightarrow g_{c,d,\mathcal{D}}(x) \geq \eta \;\;\; a.s. \]

\bibliographystyle{apalike}
\bibliographystyle{plain}
\bibliography{AG}
\end{document}